\def\VERSION{25.03.2019}
\def\users{us}\def\friends{no} 
\documentclass[12pt,a4paper,final]{article} 
\usepackage{amsmath}
\usepackage{amsthm}
\usepackage{amssymb}
\usepackage{epsfig}
\usepackage{psfrag}
\usepackage{graphicx}
\usepackage{textcomp}
\usepackage{pgf}
\usepackage{todonotes}
\numberwithin{equation}{section}
\usepackage{upgreek} 
\newtheorem{theorem}{Theorem}[section]
\newtheorem{lemma}[theorem]{Lemma}
\newtheorem{definition}[theorem]{Definition}
\newtheorem{example}[theorem]{Example} 
\newtheorem{proposition}[theorem]{Proposition}
\newtheorem{corollary}[theorem]{Corollary}
\newtheorem{remark}[theorem]{Remark}

\newcommand{\ITEM}[2]{\parbox[t]{.045\textwidth}{\textrm{#1}}%
         \hfill\parbox[t]{.95\textwidth}{#2}\pagebreak[3]\vspace{.5em}}
\usepackage{mathrsfs,cite}
\usepackage{ifthen}
\usepackage[normalem]{ulem}
\usepackage{cancel}\ifthenelse{\equal{\users}{world}}
{
\newcommand{\REM}[1]{}

	\newcommand{\INSERT}[1]{#1}
	\newcommand{\DELETE}[1]{}

        \newcommand{\COMMENT}[1]{}
        \newcommand{\TCOMMENT}[1]{}

}	
{

\usepackage[notcite,notref,color]{showkeys}
\usepackage{showkeys}
\definecolor{brown}{rgb}{0.6,0.2,0.2}
\newcommand{\REM}[1]{\marginpar{\bfseries\tiny{\color{blue}#1}}}

 \newcommand{\INSERT}[1]{{\color{brown}\uwave{#1}\color{black}}}
\newcommand{\COMMENT}[1]{\todo[inline]{#1}\noindent} 
 \newcommand{\DELETE}[1]{{\color{brown}\cancel{#1}\color{black}}}

 \newcommand{\TCOMMENT}[1]{{\color{blue}{ #1}}}

\newcount\hour \newcount\minute
\hour=\time
\divide \hour by 60
\minute=\time
\loop \ifnum \minute > 59 \advance \minute by -60 \repeat
}

\newcommand{\sym}{\mathrm{sym}}
\newcommand{\DDD}[3]{\begin{array}[t]{c}#1\vspace*{-1em}\\_{#2}\vspace*{-.45em}\\_{#3}\end{array}}
\newcommand{\ddd}[3]{\DDD{\begin{array}[t]{c}\underbrace{#1}\vspace*{.6em}\end{array}}{\text{\footnotesize #2}}{\text{\footnotesize #3}}}

\newcommand{\ol}{\overline}
\newcommand{\weak}{\rightharpoonup}

\newcommand{\bfH}{\mathbf{H}}
\newcommand{\mfh}{\mathfrak h}
\newcommand{\mfhel}{\mfh_\mathrm{el}}
\newcommand{\bftheta}{\boldsymbol\theta}

\ifthenelse{\equal{\friends}{for-friends}}
{
}

\newcommand\DT[1]{\mathchoice
                 {{\buildrel{\hspace*{.1em}\text{\LARGE.}}\over{#1}}}
                 {{\buildrel{\hspace*{.1em}\text{\Large.}}\over{#1}}}
                 {{\buildrel{\hspace*{.1em}\text{\large.}}\over{#1}}}
                 {{\buildrel{\hspace*{.1em}\text{\large.}}\over{#1}}}}
\newcommand\DDT[1]{\mathchoice
   {{\buildrel{\hspace*{.13em}\text{\LARGE.\hspace*{-.13em}.}}\over{#1}}}
   {{\buildrel{\hspace*{.1em}\text{\Large.\hspace*{-.1em}.}}\over{#1}}}
   {{\buildrel{\hspace*{.1em}\text{\large.\hspace*{-.1em}.}}\over{#1}}}
   {{\buildrel{\hspace*{.1em}\text{\large.\hspace*{-.1em}.}}\over{#1}}}}

\newcommand{\linesunder}[3]{\LSU{\begin{array}[t]{c}\underbrace{#1}\vspace*{.5em}\end{array}}{\mbox{\footnotesize\rm #2}}{\mbox{\footnotesize\rm#3}}}

\newcommand{\LSU}[3]{\begin{array}[t]{c}#1\vspace*{-1em}\\_{#2}\vspace*{-.3em}\\_{#3}\end{array}}

\newcommand{\DELTA}{\boldsymbol\delta_{\!\tau}}

\newcommand{\Vdots}{\!\;\vdots\!\;}

\renewcommand{\d}{\,\mathrm{d}}  
\newcommand{\dd}{\:\mathrm{d}} 
\newcommand{\rmd}{\mathrm{d}}    
\newcommand{\rmD}{\mathrm{D}}

\newcommand{\Cof}{\operatorname{\mathrm{Cof}}}
\newcommand{\eps}{\varepsilon}
\newcommand{\Epsilon}{\hat{\epsilon}}
\newcommand{\calD}{\mathcal D}
\newcommand{\calE}{\mathcal E}

\newcommand{\calH}{\mathcal H}
\newcommand{\calJ}{\mathcal J}
\newcommand{\calM}{\mathcal M}
\newcommand{\calR}{\mathcal R}
\newcommand{\calW}{\mathcal W}
\newcommand{\calY}{\mathcal Y}
\newcommand{\calYid}{\calY_\mathrm{id}}
\newcommand{\R}{\mathbb{R}}
\newcommand{\N}{\mathbb{N}}

\newcommand{\bbI}{\mathbb{I}}
\newcommand{\bbM}{\mathbb M}
\newcommand{\bbK}{\mathbb K}
\newcommand{\bbD}{\mathbb{D}}
\newcommand{\GDir}{{\varGamma}_\text{\sc D}}
\newcommand{\GNeu}{{\varGamma}_\text{\sc N}}
\newcommand{\SDir}{{\varSigma}_\text{\sc D}}
\newcommand{\SNeu}{{\varSigma}_\text{\sc N}}
\newcommand{\DIV}{\mathop{\mathrm{div}}\nolimits}
\newcommand{\pl}{\partial}
\newcommand{\projS}{P_{\scriptscriptstyle\textrm{\hspace*{-.2em}S}}^{}}
\newcommand{\nablaS}{\nabla_{\scriptscriptstyle\textrm{\hspace*{-.3em}S}}^{}}
\newcommand{\divS}{\mathrm{div}_{\scriptscriptstyle\textrm{\hspace*{-.1em}S}}^{}}

\newcommand{\set}[2]{\{\, #1\,; \: #2\, \} }
\newcommand{\bigset}[2]{ \big\{\, #1\,; \, #2\, \big\} }

\newcommand{\mathfrakw}{\text{\large$\mathfrak{w}$}}

\newcommand{\finite}{{large}}
\newcommand{\STEP}[1]{\medskip\noindent{\emph{#1}}}

\begin{document}
\begin{sloppypar}

\noindent{\LARGE\bfseries
Thermoviscoelasticity in \\[0.25em] Kelvin-Voigt rheology at large strains
}\bigskip\bigskip

\noindent{\Large Alexander Mielke$^{1,2}$ \ and \ Tom\'a\v s Roub\'\i\v cek$^{3,4}$}\bigskip\bigskip

\author{Alexander Mielke}

{\noindent
$^1$ Weierstra\ss-Institut f\"ur Angewandte Analysis und 
Stochastik,\\\hspace*{.8em}Mohrenstr.39, D-10117 Berlin, Germany
\\[0.3em]
$^2$
Institut f\"ur Mathematik,  
Humboldt Universit\"at zu Berlin,\\
\hspace*{.8em}Rudower Chaussee 25, D-12489 Berlin, Germany.\\
\hspace*{.8em}alexander.mielke@wias-berlin.de
\vspace{0.3em}}
\author{Tom\'a\v s Roub\'\i\v cek}\\
%
{
$^3$ Mathematical Institute, Charles University,\\
\hspace*{.8em}Sokolovsk\'a 83, CZ-186~75~Praha~8,  Czech Republic
\\[0.3em]
$^4$ Institute of Thermomechanics, Czech Academy of Sciences,\\
\hspace*{.8em}Dolej\v skova 5, CZ-182~00~Praha~8, Czech Republic\\
\hspace*{.8em}tomas.roubicek@mff.cuni.cz       
}\bigskip\bigskip

\centerline{\itshape\large In memory of Erwin Stein, who advocated
  the importance}\smallskip 
\centerline{\itshape\large of \finite-strain elasticity in engineering practice} 
\bigskip\bigskip

{\it Abstract}: The frame-indifferent thermodynamically-consistent
model of thermoviscoelasticity 
{at large strain} is formulated
in the reference configuration with using the concept of the second-grade
nonsimple materials. We focus on physically correct viscous stresses
that are frame indifferent under time-dependent rotations. Also elastic
stresses are frame indifferent under rotations and respect positivity
of the determinant of the deformation gradient.
The heat transfer is governed by the Fourier law in the actual deformed 
configuration, which leads to a nontrivial description when pulled back into 
the reference configuration. Existence of weak solutions in the
quasistatic setting, i.e.\ inertial forces are ignored, is shown by time 
discretization. 

\bigskip

\noindent
AMS Classification: 
 35K55, 
 35Q74, 
 74A15, 
 74A30, 
 80A17. 

\section{Introduction}

For a long time, thermoviscoelasticity was considered as a quite
difficult problem even at small strains, mainly because of the
nonlinear coupling with the heat-transfer equation which has no
obvious variational structure; hence special techniques had to be
developed. It took about two decades after the pioneering work by
C.M.\ Dafermos \cite{Dafe82GSSI} in one space dimension that first
three-dimensional studies occurred (cf.\ e.g.\
\cite{BlaGui00ESNS,BonBon03EUS3,Roub09TVES}). The basic new ingredient
was the $L^1$-theory for the nonlinear heat equation developed in
\cite{BDGO97NPDE,BocGal89NEPE}.  At large strains, in simple
materials, the problem is still recognized to be very
difficult even for the case of mere viscoelasticity without coupling
with temperature, and only few results are available if the physically
relevant frame-indifference is respected, as articulated by J.M.\,Ball
\cite{Ball77CIET}, see also \cite{Ball02SOPE,Ball10PPNE}.  In
particular, local-in-time existence \cite{LewMuch13LERS} or existence
of measure-valued solutions \cite{Demo00WSCN,DeStTz01VAST} are known
for simple materials. Further examples in this direction
are\cite{Tved08QEVS} for a general three-dimensional theory, but not
respecting frame indifference and the determinant constraints, or
\cite{MiOrSe14ANVM} for a one-dimensional theory using the variation
structure. While the static theory for \finite-strain elasticity
developed rapidly after \cite{Ball77CIET}, there are still only few
result for time-dependent processes respecting frame indifference as
well as the determinant constraint. The first cases were restricted to
rate-independent processes, such as elastoplasticity (cf.\
\cite{MaiMie09GERI,MieRou16RIEF}) or crack growth (cf.\
\cite{DalLaz10QCGF}, see \cite[Sec.\,4.2]{MieRou15RIST} for a
survey. Recently the case of viscoplasticity was treated in
\cite{MiRoSa18GERV}.

The main features of the model discussed in this work can be
summarized in brief as follows: the thermo-visco-elastic continuum is
formulated at large strains in a reference configuration, i.e.\ the
Lagrangian approach. The concepts of 2nd-grade nonsimple material is
used, which gives higher regularity of the deformation. The heat
transfer is modeled by the Fourier law in the actual deformed
configuration, but transformed (pulled back) into the reference
configuration for the analysis. Our model respects both static
frame-indifference of the free energy and
dynamic frame indifference for the dissipation potential. Moreover, 
the local non-selfpenetration 
is realized by imposing a blowup of the free energy if the determinant of the 
deformation gradient approaches 0 from above, however we do not enforce global 
non-selfpenetration. Also, we neglect inertial effect; 
cf.~Remark~\ref{rem-Galerkin} for 
more detailed discussion. 

Let us highlight the important aspects of the presented model 
and their consequences:

\noindent\ITEM{$(\upalpha)$}{The temperature-dependence of the free energy 
creates adiabatic effects involving the rate of the deformation gradient. 
To handle this, the Kelvin-Voigt-type viscosity is used to 
control the rate of the deformation gradient. In addition, we separate 
the purely mechanical part, cf.\ \eqref{ansatz} below, which allows 
us to decouple the singularities of \finite-strain elasticity from 
the heat equation.}

\noindent\ITEM{$(\upbeta)$}{The heat transfer itself (and also the
  viscosity from $(\upalpha)$)
is 
{clearly} rate dependent and the technique of
rate-independent processes supported by variationally efficient
energetic-solution concept cannot be used (which also 
{prevents us from excluding} possible global selfpenetration).
}

\noindent\ITEM{$(\upgamma)$}{The equations for the solid continuum
  need to be formulated and analyzed in the fixed reference
  configuration but transport processes (here only the heat transfer)
  happen rather in the actual configuration and the pull-back
  procedure needs the determinant of the deformation gradient to be
  well away from $0$. To achieve this, we exploit
  the concept of 2nd-grade nonsimple materials together
  with the results of T.J.\,Healey and S.\,Kr\"omer
  \cite{HeaKro09IWSS}, which allow us to show that the
  determinant for the deformation gradient is bounded away from $0$,
  see Section \ref{suu:HealeyK}. 
}  

\noindent\ITEM{$(\updelta$)}{The transport coefficients depend on 
the deformation gradient because of the reasons in point $(\upgamma)$.
For this, measurability in time is needed and thus the concept of global 
quasistatic minimization of deformation (as in
rate-independent systems \cite{MieRou15RIST} or in viscoplasticity in
\cite{MiRoSa18GERV}) would not be satisfactory; 
therefore we rather control the time derivative of the deformation,
which can be done either by inertia (which is neglected in our work) or
by the Kelvin-Voigt-type viscosity from $(\upalpha)$.} 

\noindent\ITEM{$(\upepsilon$)}{The 
  viscosity from $(\upalpha)$ must satisfy time-dependent frame
  indifference as explained in \cite{Antm98PUVS}, thus it is dependent
  on the rate of the right Cauchy-Green tensor rather than on the rate
  of the deformation gradient itself. However, the adiabatic heat
  sources/sinks involve terms where the rate of the deformation
  gradient occurs directly. To control the latter by the former, we
  exploit results of P.\,Neff \cite{Neff02KFIN} in the extension by
  W.\,Pompe \cite{Pomp03KFIV} for generalized Korn's inequalities,
  see Section \ref{suu:Korn}. Here, again 
  the mentioned concept of 2nd-grade nonsimple materials is used to
  control determinant of the deformation gradient, see $(\upgamma)$.}

\medskip

As mentioned above, our model heavily relies on the strain-gradient
theories to describe materials, referred as nonsimple, or also
multipolar or complex. This concept has been introduced long time ago,
cf.\ \cite{Toup62EMCS} or also e.g.\ \cite{FriGur06TBBC, MinEsh68FSTL,
  Podi02CISM, Silh88PTNB, TriAif86GALD, BalCro11LMPI} and in the
thermodynamical concept also \cite{Batr76TNEM}. In the simplest
scenario, which is also used here, the stored-energy density depends
only on the strain $F=\nabla y$ and on the first gradient $\nabla F$
of the strain. This case is called {\it 2nd-grade nonsimple
  material}. Possible generalization using only certain parts of
the 2nd in the spirit of \cite{KrPeSc19?GPEM} still need to be
explored. 

The structure of the paper is as follows. In Section \ref{se:Model} we
present the model in physical and mathematical terms. After the
precise definition of our notion of solution, Theorem
\ref{thm:MainExist} presents the main existence result for
global-in-time solutions for the \finite-strain thermoviscoelastic
system, while Corollary~\ref{Cor:ViscoElast} gives the corresponding
existence result for viscoelasticity at \finite-strain and at constant
temperature, which, to the knowledge of the authors, is also new.
A related result for isothermal \finite-strain viscoelasticity is
derived in \cite{FriKru18PNLV}, but there the limit of small strains
is treated. 

In Section \ref{se:TimeDisc} we start the proof of the main result by
introducing certain regularizations as well as a time-incremental
approach that is particularly constructed in such a split (sometimes
called staggered) way that the deformation is first updated at fixed
temperature and then the temperature is updated, where in some terms
the old and in others the new deformation is used. Another important
step in the analysis is the usage of an energy-like variable
$w=\mathfrakw(\nabla y,\theta)$ instead of temperature $\theta$, which
enables us to exploit the balance-law structure of the heat equation;
cf.\ \cite{Miel13TMER,MieMit18CEER} for arguments for the preference
of energy in favor of temperature. As an intermediate result
Proposition \ref{pr:tau-0} provides the existence of solutions
$(y_\eps,\theta_\eps)$ of the regularized problem.

In Section \ref{se:Limit} we finally show that the limit $\eps_k\to 0$
for $(y_{\eps_k},\theta_{\eps_k}) \to (y,\theta)$ can be controlled in
such a way that $(y,\theta)$ are the desired solutions. We conclude
with a few remarks concerning potential generalizations and further
applications of the methods.

\section{Modeling of thermoviscoelastic materials in the reference
configuration}
\label{se:Model}

We will use the Lagrangian approach and formulate the model in the 
reference (fixed) domain $\varOmega\subset\R^d$
being bounded with smooth boundary $\varGamma$. We assume $d\ge2$ although,
of course, the rather trivial case $d=1$ works too if $p\ge2$ is
assumed additionally to $p>d$ in \eqref{ass} below. 
We will consider a fixed time horizon $T>0$ and use the notation $I:=[0,T]$,
$Q:=I\times\varOmega$, and $\varSigma:=I\times\varGamma$.
For readers' convenience, Table~1 summarizes the main nomenclature used 
throughout the paper.


\begin{table}[ht]
\centering
\fbox{
\begin{minipage}[t]{0.42\linewidth}
\normalsize

$y$ deformation, $y(t,x) \in \R^d$,

$\theta$ absolute temperature,

$(\cdot)\!\DT{^{}}$ time derivative,

$\psi = \phi+\varphi $ free energy,

$\sigma_{\rm el}=\partial_F^{}\psi$ elastic stress,

$\sigma_{\rm vi}=\pl_{\DT F}\zeta$ viscous stress,

$F=\nabla y$ deformation gradient, 

$\mathsf{G}=\nabla F=\nabla^2y$ valued in $\R^{d\times d \times d}$,



$w$ heat part of internal energy,

$\mathfrak{h}_{\rm el}$ elastic hyperstress,


$c_{\rm v}=c_{\rm v}(F,\theta)$ heat capacity,

$\vec{q}$ heat flux,

$\calM=\Phi_\text{el}+\calH$ main mechanical energy,

$\calH$ hyperstress energy,

$\Phi_\text{cpl}$ coupling energy,

$\varPsi=\calM+\Phi_\text{cpl}$ free energy,

$\calW$ thermal energy,

$\calE=\calM+\calW$ total energy,
\end{minipage}
\hfill 
\begin{minipage}[t]{0.50\linewidth}
\normalsize

$\zeta$ potential of dissipative forces,

$\xi$ rate of dissipation (=heat production), 





$\bbK=\bbK(\theta)$ material heat conductivity, 

$\mathcal{K}=\mathcal{K}(F,\theta)$ pulled-back heat conductivity,

$C=F^\top\!F$ right Cauchy-Green tensor,

$\kappa$ heat-transfer coefficient on $\varGamma$,



$g:I{\times}\varOmega\to\R^d$ a time-dependent dead force, 

$f:I{\times}\GNeu\to\R^d$  a boundary traction, 

 $\ell$ an external mechanical loading, 

$\varOmega$ the reference domain,

$\varGamma$ the boundary of $\varOmega$, $\varGamma=\GDir\cap\GNeu$,

$I:=[0,T]$ the fixed time interval,

$Q:=I\times\varOmega$,

$\varSigma:=I\times\varGamma$, 

$\mathscr{H}=\mathscr{H}(\nabla F)$ the potential of $\mathfrak{h}_{\rm el}$,


 $\calY_0,\calYid$ sets of admissimble deformations,

$\mathrm{GL}^+(d):=\{A\in\R^{d\times d};\ \det A>0\}$,

 $\mathrm{SO}(d):=\{A\in\mathrm{GL}^+(d);\ A^\top A=\bbI=AA^\top\}$.

\end{minipage}
} 
\medskip

\centerline{{\normalsize\bf Table~1.}
 {\normalsize\sl Summary of the basic notation used throughout the paper.}}
\label{Tab_Notation}
\end{table}

To introduce our model in a broader context, we may define  
the {\it total free energy} and the {\it total dissipation potential} 
\begin{align}
 \varPsi(y,\theta)=\int_\varOmega \! \big(\psi(\nabla y,\theta) 
         +\mathscr{H}(\nabla^2y)\big) \d x
\ \ \text{ and }\ \ 
\mathcal{R}(y,\DT y,\theta)=\int_\varOmega\zeta(\nabla y,\nabla\DT y,\theta)\,\d x,
\end{align}
respectively.
The mechanical evolution part can then be viewed as an abstract gradient flow 
\begin{align}\label{grad-flow}
\rmD_{\DT y}\mathcal{R}(y,\DT y,\theta)+\rmD_y \varPsi(y,\theta)
=\ell (t)
\ \text{ with }
\langle\ell(t) ,y \rangle=\int_\varOmega\!
g(x,t){\cdot}y(x)
\d x+\int_{\GNeu}\!\!f(x,t){\cdot}y(x)
\dd S,
\end{align}
cf.\ also \cite{Tved08QEVS,MiOrSe14ANVM} for the isothermal case 
and \cite{Miel11FTDM} for the general case. The sum of the
conservative and the dissipative parts corresponds to the {\it
  Kelvin-Voigt rheological model} in the quasistatic variant
(neglecting inertia). The notation ``$\,\pl\,$'' is used for partial
derivatives (here functional or later in Euclidean spaces), while
$(\cdot)'$ will occasionally be used for functions of only one
variable.

Writing \eqref{grad-flow} locally in the classical formulation, one
arrives at the nonlinear parabolic 4th-order partial differential
equation expressing quasistatic {\it momentum equilibrium}
\begin{align}\label{moment-eq}
\DIV \,\sigma+g=0 \qquad \text{ with }\quad
\sigma=\sigma_\mathrm{vi}+\sigma_\mathrm{el}
-\DIV \,\mfhel,
\end{align}
where the viscous 
stress is $\sigma_\mathrm{vi}=\sigma_\mathrm{vi}(F,\DT F,\theta)$ and the
elastic stress is $\sigma_\mathrm{el}=\sigma_\mathrm{el}(F,\theta)$,
 while $\mfhel$ is a so-called hyperstress
arising from the 2nd-grade nonsimple material concept, cf.\ e.g.\ 
\cite{Podi02CISM,Silh88PTNB,Toup62EMCS}.
In view of the local potentials used in \eqref{grad-flow}, we have
\begin{align}\label{potentials}
&\sigma_\mathrm{vi}(F,\DT F,\theta)=\partial_{\DT F}\zeta(F,\DT
  F,\theta),
\quad
\sigma_\mathrm{el}(F,\theta)=\partial_F^{}\psi(F,\theta),
\quad \text{ and} \quad 
\mfhel(\mathsf{G})=\mathscr{H}'(\mathsf{G}),
\end{align}
\mbox{}where $\mathsf{G}\in \R^{d\times d \times d}$ is a
  placeholder for $\nabla F$.

\mbox{}
{An important physical}
requirement is {\it static and dynamic frame indifference}. 
For the elastic stresses, static frame indifference means that 
\begin{subequations}
\begin{align}
  \sigma_\mathrm{el}(RF,\theta)= R\, \sigma_\mathrm{el}(F,\theta)
  \quad \text{and} \quad  
   \mfhel(R\mathsf{G})=R\mfhel(\mathsf{G})
\end{align} 
for all $R\in\mathrm{SO}(d)$, $F$ and $\mathsf{G}$. 
For the viscous stresses\INSERT{,} dynamic frame indifference means that
\begin{align}
\sigma_\mathrm{vi}(RF,\DT RF{+}R\DT
   F,\theta)=R\,\sigma_\mathrm{vi}(F,\DT F,\theta)
\end{align}
\end{subequations}
for all smoothly time-varying $R:t\mapsto R(t)\in\mathrm{SO}(d)$,
cf.\ \cite{Antm98PUVS}.  
{Note that $R$ may depend on $t$ but
  not on $x\in \varOmega$, since frame-indifference relates to
  superimposing time-dependent \emph{rigid-body motions}.}

In terms of the 
thermodynamic potentials 
$\zeta$, $\psi$, and $\mathscr H$, these frame indifferences read as
\begin{subequations}
\label{frame-indif}
\begin{align}
&\psi(RF,\theta)=\psi(F,\theta), \quad 
\mathscr{H}(R\nabla F)=\mathscr H(\nabla F), 
 \quad \text{and}
\\ 
&\zeta(RF,\theta;(RF)\!\DT{^{}}\hspace{.0em})
=\zeta(RF,\theta;\DT{R}F{+}R\DT{F})=\zeta(F,\theta;\DT{F})
\end{align}
\end{subequations}
for $R$, $F$ and $\nabla F$ as above.  These frame indifferences imply
the existence of reduced potentials $\hat\psi$, $\hat\zeta$, and
$\hat{\mathscr H}$ such that
\begin{align}
\label{hat-ansatz}
\zeta(F,\DT F,\theta)=\hat\zeta(C,\DT C,\theta), \quad 
\psi(F,\theta)=\hat\psi(C,\theta),  \quad \text{and}\quad
\mathscr H(\mathsf{G}) = \hat{\mathscr H}(\mathsf B) 
\end{align}
where $\mathsf B=\mathsf{G}^\top\!\cdot \mathsf{G} \in \R^{(d\times
  d)\times (d \times d)} $, and $C\in \R^{d \times d}_\mathrm{sym}$ is
the right Cauchy-Green tensor $C=F^\top F$ with time derivative $\DT
C=\DT F^\top F+F^\top\DT F$.  More specifically, denoting
$\mathsf{G}=[\mathsf{G}_{\alpha ij}]$ the placeholder for
$\frac{\pl}{\pl x_j}F_{\alpha i}$ with $F_{\alpha i}$ the placeholder
for $\frac{\pl}{\pl x_i}y_{\alpha}$, the exact meaning is $ [
\mathsf{G}^\top \! \cdot \mathsf{G}]_{ijkl} :=
\sum_{\alpha=1}^d\mathsf{G}_{\alpha ij}\mathsf{G}_{\alpha kl}$ and
$[F^\top F]_{ij}:=\sum_{\alpha=1}^dF_{\alpha i}F_{\alpha j}$.  The
ansatz \eqref{hat-ansatz} also means that
\begin{subequations}
\label{KV-large-F-vs-C}
\begin{align}
&\sigma_\mathrm{el} (F,\theta):= \pl_F\psi(F;\theta) 
 =2 F\pl_C^{}\hat\psi(F^\top F,\theta) 
 =2 F\pl_C^{}\hat\psi(C,\theta), \\ 
& \mfhel (\mathsf{G}) := \pl_{\mathsf G} \mathscr
H (\mathsf G) = 2\mathsf{G} \pl_{\mathsf B} \hat{\mathscr
  H}(\mathsf{G}^\top\!\!\cdot \mathsf{G}) = 2\mathsf{G} \pl_{\mathsf
  B} \hat{\mathscr H}(\mathsf B), \\
&\sigma_\mathrm{vi}(F,\DT F,\theta)
 :=\pl_{\DT F}\zeta(F,\DT F,\theta) 
= 2 F\pl_{\DT C}\hat\zeta(F^\top \! F,\DT F^\top\! F{+}F^\top\DT F,\theta)
= 2 F\pl_{\DT C}\hat\zeta(C,\DT C,\theta).
\end{align}
\end{subequations}

The simplest choice, which is adopted in this paper for avoiding
unnecessary technicalities, is that the viscosity $\sigma_\mathrm{vi}$
is linear in $\DT C$.  This is the relevant modeling choice for
non-activated dissipative processes with rather moderate rates (in
contrast to activated processes like plasticity having nonsmooth
potentials that are homogeneous of degree 1 in a small-rate
approximation). This linear viscosity leads to a potential which is
quadratic in $\DT C$, viz.
\begin{align}\label{hat-zeta-special}
\hat\zeta(C,\DT C,\theta):=\frac12\,\DT C : \bbD(C,\theta)\DT C\,.
\end{align}
Although for this choice the material viscosity is linear, the
geometrical nonlinearity arising from large strains is still a vital
part of the problem due to the requirement of frame indifference.
Note that $\sigma_\mathrm{vi}(F,\DT F,\theta)$ necessarily depends on
$F$ if we express $\DT C$ in terms of the velocity gradients $\DT F$,
even if $\bbD$ is constant: $\sigma_\mathrm{vi}(F,\DT F,\theta) =
2F \bbD(C,\theta)(\DT F^\top F{+}F^\top \DT F)$. While we will be able
to handle general dependence
on $F$, it will be a crucial restriction that $\DT F \mapsto
\sigma_\mathrm{vi}(F,\DT F,\theta)$ is linear. 

Furthermore, the {\it specific dissipation rate} can be simply identified 
in terms of $\hat\zeta$ as 
\begin{align}\nonumber
\xi(F,\DT F,\theta)&=\sigma_\mathrm{vi}(F,\DT F,\theta){:}\DT F=
2F\pl_{\DT C}\hat\zeta(F^\top F,\DT F^\top F{+}F^\top \DT F,\theta){:}\DT F
\\&=
 \pl_{\DT C}\hat\zeta(F^\top F,\DT F^\top F{+}F^\top \DT F,\theta){:}
(\DT F^\top F{+}F^\top \DT F) =\pl_{\DT C}\hat \zeta(C,\DT C 
  ,\theta){:}\DT C. 
\label{def-of-xi}
\end{align}
For our choice \eqref{hat-zeta-special}, we simply have 
$\xi(F,\DT F,\theta)=\bbD(C,\theta)\DT C{:}\DT C
=2\hat\zeta(C,\DT C,\theta)=2\zeta(F,\DT F,\theta)$. 

In brief, the standard thermodynamical arguments start from the free
energy density $\psi$ and the definition of {\it entropy} via
$s=-\partial_{\theta}^{}\psi$
{(here $\mathscr H$ does play no role as it is chosen to be
  independent of $\theta$)} and the entropy equation
\begin{align}
\label{entropy-eq}
\theta\DT s=\xi-\DIV \,\vec{q}
\end{align}
with 
{the dissipation rate} $\xi$ from \eqref{def-of-xi} and 
the heat flux $\vec{q}$. We further use the 
formula $\DT s=-\pl_{\theta\theta}^2\psi\,\DT\theta-\pl_{F\theta}^2\psi{:}\DT F$ 
and the  {\itshape Fourier law} formulated in the reference configuration 
\begin{align}
\vec{q}=-\mathcal{K}(F,\theta)\nabla\theta, 
\end{align}
which will be specified later in \eqref{K-pull-back}. Altogether, we
arrive at the coupled system
\begin{subequations}
 \label{system}
\begin{align}\nonumber
&\DIV \!\big(\sigma_\mathrm{vi}(\nabla y,\nabla\DT y,\theta)
  +\sigma_\mathrm{el}(\nabla y,\theta) 
  -\DIV \mfhel(\nabla^2y)\big)+g
\\ 
&\qquad\qquad\qquad\qquad \text{with }\ 
\sigma_\mathrm{vi}(F,\DT F,\theta)=  \pl_{\DT F}
\zeta(F,\DT F,\theta) 
\  \text{and }\  \sigma_\mathrm{el}(F,\theta)=\pl_F^{}\psi(F,\theta)
\,,\label{momentum-eq}
\\&\nonumber
c_\mathrm{v}(\nabla y,\theta)\DT\theta
  =\DIV \!\big(\mathcal{K}(\nabla y,\theta)\nabla\theta\big)
   +\xi(\nabla y,\nabla\DT y,\theta)
   +\theta\pl_{F\theta}^2\psi(\nabla y,\theta){:}\nabla\DT y
\\
&\qquad\qquad\qquad\qquad\text{with }\ 
c_\mathrm{v}(F,\theta)=-\theta\pl_{\theta\theta}^2\psi(F,\theta)\
\text{ and } \ \xi
\ \text{ from \eqref{def-of-xi}}
\label{heat-equation}
\end{align}
\end{subequations}
on $Q$.  We complete \eqref{system} by some boundary conditions.  For
simplicity, we only consider a mechanically fixed part $\GDir$ time
independent undeformed (i.e.\ identity) while the whole boundary is
thermally exposed with a phenomenological heat-transfer coefficient
$\kappa \ge 0$:
\begin{subequations}
  \label{BC}
 \begin{align}\label{BC1} 
&\big(\sigma_\mathrm{vi}(\nabla y,\nabla\DT y,\theta)
+\sigma_\mathrm{el}(\nabla y,\theta)\big)\vec{n}
-\divS\big(\mfhel(\nabla^2y)\vec{n}\big) =f&&\text{on }\ \GNeu,
\\
& \label{BC2}
y(x)=x\ \ \ \text{(identity)}&&\text{on }\ \GDir,
\\&\label{BC3} 
\mfhel(\nabla^2y){:}(\vec{n}\otimes\vec{n})=0&&\text{on }\ \varGamma,
\\
& \label{BC4}
\mathcal{K}(\nabla y,\theta)\nabla\theta\cdot\vec{n}
+\kappa\theta=\kappa\theta_\flat&&\text{on }\ \varGamma,
\end{align}
\end{subequations}
where $\vec n$ is the outward pointing normal vector, and
$\theta_\flat$ is a given external temperature.  Moreover, 
following \cite{Beto86KSMC} the surface divergence ``$\divS$'' in
\eqref{BC1} is defined as
$\divS(\cdot)=\mathrm{tr}\big(\nablaS(\cdot)\big)$, where
$\mathrm{tr}(\cdot)$ denotes the trace and $\nablaS$ denotes the
surface gradient given by $\nablaS v=(\mathbb I- \vec n{\otimes}\vec
n)\nabla v= \nabla v-\frac{\partial v}{\partial\vec{n}}\vec{n}$. See
\eqref{momentum-weak+++} for a short mathematical derivation of the
boundary conditions \eqref{BC1} and \eqref{BC3}, and
\cite[pp.\,358-359]{Stei15GFCM} for the mechanical interpretation in
second-order materials. 

 
In order to facilitate the subsequent mathematical analysis,
we assume a rather weak thermal coupling through the free
energy (together with the coupling through the temperature-dependent
viscous dissipation). To distinguish the particular 
coupling thermo-mechanical term from the purely mechanical one, 
we consider the explicit ansatz
\begin{align}\label{ansatz}
\psi(F,\theta)=\varphi(F)+\phi(F,\theta)
\ \ \ \text{ with }\ \ \phi(F,0)=0.
\end{align}
In applications, the internal energy $e$ given by Gibbs' relation
$$
e=e(F,\theta)=\psi(F,\theta)+\theta s
=\psi(F,\theta)-\theta\pl_\theta\psi(F,\theta)=
\psi(F,\theta)-\theta\pl_\theta\phi(F,\theta).
$$
is often balanced. Here, we rather use the 
thermal part of the internal energy $w:=e-\varphi(F)$. In view of the ansatz 
\eqref{ansatz}, we have 
\begin{align}
w=\mathfrakw(F,\theta)=\psi(F,\theta)-\theta\pl_\theta\phi(F,\theta)
-\psi(F,0)=\phi(F,\theta)-\theta\pl_\theta\phi(F,\theta).
\label{def-of-w}\end{align}
Note that $\mathfrakw(F,\cdot)$ is the primitive function of the
specific heat $c_\mathrm{v}(F,\cdot)$ calibrated as
$\mathfrakw(F,0)=0$, so that also $e=\psi(F,0)+w$. The heat-transfer
equation \eqref{heat-equation} simplifies as
\begin{align}
\DT w-\DIV \big(\mathcal{K}(\nabla y,\theta)\nabla\theta\big)
=\xi(\nabla y,\nabla\DT y,\theta)+
\pl_F^{}\phi(\nabla y,\theta){:}\nabla\DT y
\ \text{ with }\ w=\mathfrakw (F,\theta)\,.
\label{heat-equation+}
\end{align}
In particular, the purely mechanical stored energy $\varphi$ does not
occur in \eqref{def-of-w} and does not influence the heat production and 
transfer \eqref{heat-equation+}.

The energetics of the system \eqref{system}--\eqref{BC} can be 
best described by introducing additional energy functionals as
follows:
\begin{subequations}
\label{eq:Energies}
\begin{align}
\label{eq:Energ.H}
&\calH(y)&&:=\int_\varOmega \mathscr H(\nabla^2y) \d x &&\text{hyperstress
  energy},
\\
\label{eq:Energ.M}
&\calM(y)&&:= \calH(y){+}\Phi_\text{el}(y) \text{ with } 
 \Phi_\text{el}(y):=\int_\varOmega \varphi(\nabla y)\d x\  && 
  \text{main mech.\ energy} ,&&
\\
\label{eq:Energ.cpl}
&\Phi_\text{cpl} (y,\theta)\hspace{-2em}&&:=\int_\varOmega \phi(\nabla y,\theta)\d x
  &&\text{coupling energy},
\\
\label{eq:Energ.F}
& \varPsi(y,\theta)&&:= \calM(y) + \Phi_\text{cpl}(y,\theta) 
   && \text{free energy} ,
\\
\label{eq:Energ.W}
&\calW(y,\theta)&&:= \int_\varOmega \mathfrakw(\nabla y,\theta) \d x
     &&\text{thermal energy},
\\
\label{eq:Energ.E}
&\calE(y,\theta)&&:= \calM(y)+ \calW(y,\theta) &&\text{total energy}. 
\end{align}
\end{subequations}
An mechanical energy balance is 
revealed by testing \eqref{momentum-eq} by $\DT y$ and
\eqref{heat-equation} by 1, and using the boundary conditions after
integration over $\varOmega$ and using Green's formula twice together
with another $(d{-}1)$-dimensional Green formula over $\varGamma$ for
\eqref{momentum-eq} and once again Green's formula for
\eqref{heat-equation}.  The last mentioned technique is related with
the concept of nonsimple materials; for the details about how the
boundary conditions are handled see e.g.\
\cite[Sect.\,2.4.4]{Roub13NPDE}.  This test of \eqref{momentum-eq}
gives the mechanical energy balance:
\begin{align}\label{mech-engr}
\int_\varOmega\!\!\!\!\ddd{\xi(\nabla y,\nabla\DT y,\theta)_{_{_{_{}}}}}{dissipation}{rate}\!\!\!\!
+\!\!\!\!\ddd{\sigma_\mathrm{el}{:}\nabla\DT y_{_{_{_{}}}}}{mechanical}{power}\!\!\!\!\d x
+\frac{\rmd}{\rmd t}
 \calH(y)
=\int_\varOmega\!\!\!\!\!\!\!\!\!\!\ddd{g\cdot\DT y_{_{_{_{}}}}}{power of the}{bulk force}
\!\!\!\!\!\!\!\!\!\!\d x
+\int_{\GNeu}\!\!\!\!\!\!\ddd{f\cdot\DT y_{_{_{_{}}}}}{power of}{the traction}
\!\!\!\!\!\d S.
\end{align}
Using $\sigma_\text{el} = \pl_F\varphi + \pl_F \phi$ and
integrating in time leads to the relation  
\begin{equation}
  \label{eq:MechEnergBal}
\calM(y(T)) + \int_0^t \! \int_\varOmega \Big(\int_\varOmega\xi(\nabla
y,\nabla\DT y,\theta) + \pl_F \phi(\nabla y,\theta){:}\nabla \DT y\Big) \d x 
\d t = \calM(y(0)) + \int_0^t \langle \ell,\DT y \rangle \d t.  
\end{equation}
that will be very useful for obtaining a priori estimates in the
following sections. 

Next, we test the heat equation in its simplified form
\eqref{heat-equation+} together with the boundary conditions
\eqref{BC4} by the constant function 1 (i.e.\ we 
merely integrated over  
$\varOmega$) and add the result to \eqref{eq:MechEnergBal}. After
major cancellations we obtain the total energy balance:
\begin{align}\label{total-engr}
\frac{\rmd}{\rmd t} \calE(y, \theta) 
=\int_\varOmega\!\!\!\!\!\!\!\!\!\!\!\!\ddd{g\cdot\DT y_{_{_{_{}}}}}{power of mecha-}{nical bulk load}\!\!\!\!\!\!\!\!\!\!\!\!
\d x
+\!\int_{\GNeu}\!\!\!\!\!\!\!\ddd{f\cdot\DT y_{_{_{_{}}}}}{power of}{the traction}
\!\!\!\!\!\!\d S - \!
\int_\varGamma\!\!\!\!\!\!\!\!\ddd{\kappa(\theta{-}\theta_\flat)_{_{_{_{}}}}}
 {power of the}{external heating}\!\!\!\!\!\! \d S. 
\end{align}
In particular, we see that the total energy is conserved up to the
work induced by the external loadings or the flux of heat through the
boundary.

{}From the entropy equation \eqref{entropy-eq}, we can read the total
entropy balance (the {\it Clausius-Duhem inequality}):
\begin{align}\nonumber
\frac{\rmd}{\rmd t}\int_\varOmega s(t,x)\,\d x&=\int_\varOmega\frac{\xi
+\DIV (\mathcal{K}\nabla\theta)}\theta\,\d x
=
\int_\varOmega\frac{\xi
}\theta-\mathcal{K}\nabla\theta{\cdot}\nabla\frac1\theta\,\d x
+\int_\varGamma\frac{\mathcal{K}\nabla\theta}\theta{\cdot}\vec{n}\,\d S
\\[.1em]&=
\int_\varOmega\!\!\!\!\linesunder{\frac{\xi
}{\theta}
+\frac{\mathcal{K}\nabla\theta{\cdot}\nabla\theta}{\theta^2}}
{entropy-production}{rate}\!\!\!\!\d x
\ +\ \int_\varGamma\frac{\mathcal{K}\nabla\theta}\theta{\cdot}\vec{n}\,\d S
\ge\int_\varGamma\!\!\!\!\!\!\!\!\!\!\!\!\!
\linesunder{\ \frac{-\vec{q}}\theta{\cdot}\vec{n}\ }
{entropy flux}{through boundary}\!\!\!\!\!\!\!\!\!\!\!\!\!\d S.
\label{ent-inquality}\end{align}
This articulates, in particular, the second law of thermodynamics that
the total entropy in the isolated systems (i.e.\ here $\vec{q}=0$ on
$\varGamma$) is nondecreasing with time provided $\mathcal{K} =
\mathcal{K} (\nabla y,\theta)$ is positive semidefinite and the
dissipation rate is non-negative.

It is certainly a very natural modeling choice that Fourier's law is
formulated in the actual (also called the deformed) configuration in a
simple form, namely the \emph{actual heat flux} is given by
\begin{align}
  \label{Fourier}
 \vec{\mathsf q}=-\bbK(\uptheta)\nabla_z\uptheta, \quad
  \text{where }z=y(x) \text{ and }\uptheta(z) = \theta(y^{-1}(z))\
  \text{ for $z\in y(\varOmega)$} 
\end{align}

with the heat-conductivity tensor $\bbK= \bbK(x,\theta)$ considered as
a material parameter possibly dependent on $x\in\varOmega$.  We transform
(i.e.\ {\it pulled-back}) this {\it Fourier law} into the reference
configuration via the heat flux $\vec q(x)=\mathcal K(x) \nabla
\theta = \mathcal K (\nabla y(x))^\top \nabla_z\uptheta (y(x))$ and
$\vec q = (\Cof F^\top) \vec{\mathsf q}$, because fluxes should be
considered as $(d{-}1)$-forms. With \eqref{Fourier} the usual
transformation rule for 2nd-order contra-variant tensors yields
the heat-conductivity tensor
\begin{align}
 \mathcal{K}(x,F,\theta)
 & =(\Cof F^\top)\bbK(x,\theta)F^{-\top}
  =\frac{(\Cof F^\top)\bbK(x,\theta)\Cof F}{\det F}
  =(\mathrm{det}F)F^{-1}\bbK(x,\theta)F^{-\top}
 \label{K-pull-back}
\end{align}
{if $\det F>0$, whereas } the case $\det F\le0$ is considered
nonphysical, so $\mathcal{K}$ is then not defined. Here we used the
standard shorthand notation $F^{-\top}=[F^{-1}]^\top=[F^\top]^{-1}$
and also the algebraic formula $F^{-1}=(\Cof F^\top)/\det F$.  In what
follows, we omit explicit $x$-dependence for notational simplicity.
Let us emphasize that in our 
  formulations $\nabla\theta$ is \emph{not} treated as a vector, but a
  contravariant 1-form. Starting from $\theta(x)=\uptheta(y(x))$ the
  chain-rule gives $\nabla (x)= \nabla y(x)^\top
  \nabla_Y\bftheta(y(x))$.
It should be noted that \eqref{Fourier} is rather formal argumentation,
assuming injectivity of the deformation $y$ and thus existence of $y^{-1}$,
which is however not guaranteed in our model; anyhow, handling only local 
non-selfpenetration while ignoring possible global selfpenetration  
is our modeling approach often accepted in engineering, too.

For the isotropic case $\bbK(\theta)=\varkappa(\theta)\bbI$, relation
\eqref{K-pull-back} can also be written by using the right
Cauchy-Green tensor $C=F^\top F$ as $\mathcal
K=\det(F)\varkappa(\theta)C^{-1}$, cf.\ e.g.\ \cite[Formula
(67)]{DuSoFi10TSMF} or \cite[Formula (3.19)]{GovSim93CSD2} for the
mass instead of the heat transport.  In principle, $\bbK$ in
\eqref{Fourier} itself may also depend on $C=F^\top F$, which we
omitted to emphasize that $\mathcal{K}$ in \eqref{K-pull-back} will
depend on $F$ anyhow.

In what follows, we will use the (standard) notation for the Lebesgue
$L^p$-spaces and $W^{k,p}$ for Sobolev spaces whose $k$-th
distributional derivatives are in $L^p$-spaces and the abbreviation
$H^k=W^{k,2}$.  The notation $W^{1,p}_{\rmD}$ will indicate the closed
subspace of $W^{1,p}$ with zero traces on $\GDir$.  Moreover, we will
use the standard notation $p'=p/(p{-}1)$.  In the vectorial case, we
will write $L^p(\varOmega;\R^n)\cong L^p(\varOmega)^n$ and
$W^{1,p}(\varOmega;\R^n)\cong W^{1,p}(\varOmega)^n$.  Thus, for example,
\begin{align}
H^{1}_{\rmD}(\varOmega;\R^{d}):=\big\{v\in L^2(\varOmega;\R^d);\ \nabla v\in 
 L^2(\varOmega;\R^{d\times d}),\ \ v|_{\GDir}^{}=0\,\big\}.
\end{align}
For the fixed time interval $I=[0,T]$, we denote by $L^p(I;X)$ the
standard Bochner space of Bochner-measurable mappings $I\to X$ with
$X$ a Banach space. Also, $W^{k,p}(I;X)$ denotes the Banach space of
mappings from $L^p(I;X)$ whose $k$-th distributional derivative in
time is also in $L^p(I;X)$. The dual space to $X$ will be denoted by
$X^*$.  Moreover, $C_{\rm w}(I;X)$ denotes the Banach space of weakly
continuous functions $I\to X$.  The scalar product between vectors,
matrices, or 3rd-order tensors will be denoted by ``$\,\cdot\,$'',
``$\,:\,$'', or ``$\,\Vdots\,$'', respectively. Finally, in what
follows, $K$ denotes a positive, possibly large constant.

We consider an initial-value problem, imposing the initial conditions
\begin{align}\label{IC}
y(0,\cdot)=y_0\quad  \text{ and } \quad \theta(0,\cdot)=\theta_0 \quad
\text{ on }\ \varOmega.
\end{align}

Having in mind the form \eqref{heat-equation+} of the heat equation, we can now
state the following definition for a weak solution: 

\begin{definition}[Weak solution]
\label{def}
  A couple $(y,\theta):Q=[0,T]{\times}\varOmega \to \R^d\times\R$ is
  called a \emph{weak solution} to the initial-boundary-value problem
  \eqref{system}\,\&\,\eqref{BC}\,\&\,\eqref{IC} 
  if $(y,\theta)\in C_{\rm w}(I;W^{2,p}(\varOmega;\R^d))\times L^1(I;W^{1,1}(\varOmega))$
  with $\nabla\DT y\in L^2(Q;\R^{d\times d})$, if $\min_Q\det\nabla y>0$
  and $y|_{\SDir}= \text{\rm identity}$, and if it satisfies 
  the integral identity
\begin{subequations}
\begin{align}\nonumber&
\int_0^T\bigg(\int_\varOmega
 \Big( \big(\sigma_\mathrm{vi}(\nabla y ,\nabla\DT y,\theta)
+  \sigma_\mathrm{el}  (\nabla y,\theta)\big) {:} \nabla z + \mfhel
   (\nabla^2 y)\vdots \nabla^2 z \bigg) \,\d x \d t
\\[-.5em]
&\hspace*{18em}
= \int_Q g{\cdot}z\,\d x\d t  +\int_{\SNeu}\!\!f{\cdot}z\,\d S\d t 
\label{momentum-weak+}
\intertext{for all smooth $z:Q\to\R^d$ with $z=0$ on $\SDir$ together
  with $y(0,\cdot)=y_0$, and if}
&\nonumber
\int_Q 
\mathcal{K}(\nabla y,\theta)\nabla\theta{\cdot}\nabla v
-\big(\xi(\nabla y,\nabla\DT y,\theta) {+}
\pl_F^{}\phi(\nabla y,\theta)
{:}\nabla\DT y\big)v
-\mathfrakw(\nabla y,\theta)\DT v\, \d x\d t
\\[-.5em]
&\qquad\qquad\qquad\qquad
+\int_\varSigma\kappa\theta v\,\d S\d t
=\int_\varSigma\kappa\theta_\flat v\,\d S\d t+
\int_\varOmega \mathfrakw(\nabla y_0,\theta_0)v(0)\,\d x
\label{heat-weak}
\end{align}
\end{subequations}
for all smooth $v:Q\to\R$ with $v(T)=0$, where $\mathfrakw$ is
defined in \eqref{def-of-w}. 
\end{definition}

At first sight,  it seems that \eqref{momentum-weak+} is not
suited to apply the test function $z = \DT y$, which is the natural
and necessary choice for deriving energy bounds. Obviously, we will
not be able to obtain enough control on $\nabla^2\DT
y$. However, using the abstract chain rules provides in Section
\ref{suu:ChainRule} this problem can be handled by extending 
$\calH(y)=\int_\varOmega \mathscr H(\nabla^2y)\dd x $ to a lower
semicontinuous and convex functional on $H^1(\varOmega;\R^d)$ by setting
it $\infty$ outside $W^{2,p}(\varOmega; \R^d)$, see the
rigorous proof of \eqref{eq:BalMechEnerg} in Step 3 of the proof of Proposition
\ref{pr:tau-0}.   

It will be somewhat technical to see that the weak
formulation \eqref{momentum-weak+} is indeed selective enough, in the
sense that for sufficiently smooth solutions one can indeed obtain 
the classical formulation \eqref{system} together with the boundary
conditions \eqref{BC}, cf.\ also 
\cite[Sect.\,2.4.4]{Roub13NPDE}.  In particular, abbreviating
$\sigma=\sigma_\mathrm{vi}(\nabla y,\nabla\DT y,\theta)
+\sigma_\mathrm{el}(\nabla y,\theta)$, integrating by part once, and 
 using the boundary conditions (\ref{BC}a,c) yields 
\begin{align}
 \label{momentum-weak++} 
&\int_Q\!\!\Big(\big(\sigma{-}\DIV\mfhel(\nabla^2y)\big){:}\nabla z
 -g{\cdot}z \Big)\d x\d t=\int_{\SNeu}\!\!f{\cdot}z \d S\d t-
 \int_\varSigma\mfhel(\nabla^2y)\Vdots(\nabla z {\otimes} \vec{n} )\d S\d t.
\end{align}

We now want to show how the strong form
\eqref{momentum-eq} and the associated boundary conditions
(\ref{BC}a,c) follow from \eqref{momentum-weak++}. For this goal,
we apply Green's formula in the opposite direction to remove 
$\nabla$ in front of the test function $z$. Using also the orthogonal
decomposition of $\nabla z=\nablaS z + \frac{\pl}{\pl\vec{n}}z\otimes \vec{n} $
involving the surface gradient $\nablaS z$ and writing shortly
$\mfh$ for $\mfhel(\nabla^2y)\in \R^{d\times d \times d}$, relation
\eqref{momentum-weak++} leads to the identity 
\begin{align}
 \nonumber
&\int_Q\big({-}\DIV\,\sigma+\DIV^2\mfh -g\big){\cdot}z\,\d x\d t
\\[-.3em]&\nonumber=
\int_\varSigma\Big(\big(\sigma  {-}\DIV \mfh \big) :
( z {\otimes} \vec{n} ) - \mfh \Vdots (\nabla z {\otimes} \vec{n}
   ) \Big) \d x \d t + \int_{\SNeu}f{\cdot}z\,\d S\d t
\\[-.3em]&\nonumber=
\int_\varSigma \!\Big((\sigma {-}\DIV \mfh )\vec{n}{\cdot}z + 
 \big(\mfh:(\vec{n} {\otimes} \vec{n})\big)\cdot \frac{\pl z}{\pl\vec{n}}
        +\mfh\vec n :\nablaS z) \Big)\d S \d t
-\int_{\SNeu}f{\cdot}z\,\d S\d t
\end{align}
Using the surface divergence $\divS$ and the projection
$\projS:A\mapsto A-A\vec{n} \otimes \vec{n} $ to the tangential part, we
obtain the integration by parts formula (cf.\ \cite{Beto86KSMC} or
\cite[pp.\,358-359]{Stei15GFCM}) 
\[
\int_{\pl\varOmega} A:\nablaS z\d S 
 =\int_{\pl\varOmega} (\projS A):\nablaS z\d S 
   = -\int_{\pl\varOmega} \divS(\projS A)\cdot z \d S,
\]
where the surface $\varGamma$ is now assumed to be sufficiently smooth. 
Using this with $A=\mfh\vec{n}$ for the previous relation we find 
\begin{align} 
\nonumber
&\int_Q\big({-}\DIV\,\sigma+\DIV^2\mfh -g\big){\cdot}z\,\d x\d t
\\&
=\int_{\SNeu} \!\!\! \Big((\sigma{-}\DIV \mfh)\vec{n}
  -\divS\big(\projS (\mfh \vec{n})\big)  
   -f\Big){\cdot}z \d S\d t
+\int_\varSigma \!\!\big(\mfh{:}(\vec{n} {\otimes} \vec{n})\big)\cdot \frac{\pl
  z}{\pl\vec{n}} \d S\d t,
\label{momentum-weak+++}
\end{align}
where we have used $z=0$ on $\SDir =\varSigma\setminus \SNeu$. Now,
taking $z$'s with a compact support in $Q$, we obtain the equilibrium
\eqref{momentum-eq} in the bulk. Next taking taking $z$'s with zero
traces on $\varSigma$ but general $ \frac{\pl
  z}{\pl\vec{n}}$, we obtain \eqref{BC3}. Note that the latter
condition implies $\projS (\mfh \vec{n})=\mfh \vec{n} -
\big(\mfh:(\vec{n}{\otimes}\vec{n})\big)\otimes\vec{n} = \mfh
\vec{n}$. Hence, taking finally general $z$'s, we obtain
\eqref{BC1}, as $\projS$ can be dropped because of \eqref{BC3}. 

Moreover, also note that, from the integral identity
\eqref{heat-weak}, one can read $ \mathfrakw ( \nabla y(0), \theta(0))
= \mathfrakw(\nabla y_0,\theta_0)$ from which
$\theta(0)=\theta_0$ follows when taken the invertibility of
$\mathfrakw(F,\cdot)$ and $y(0)=y_0$ into account.

Now we exploit the decomposition \eqref{ansatz} of $\psi$ into
$\phi$ and $\varphi$, which allows us to impose coercivity
assumptions for the purely elastic part $\phi$ that are
independent of those for $\varphi$,  namely
%

%
\begin{subequations}
\label{ass}
\begin{align}
&\nonumber \hspace*{-0.7em}\exists\, 
p\in {]d,\infty[}\cap {[2,\infty[},\ s>0,\ q\ge pd/(p{-}d)\ \ \exists\,
\alpha,K,\Epsilon>0: 
\\
\nonumber&\varphi:\mathrm{GL}^+(d)\to\R^+\ \text{ twice continuously 
differentiable},\ \forall\, F\in\mathrm{GL}^+(d):
\\
&\label{ass-psi}
\qquad
\varphi(F)\ge\Epsilon|F|^s+\Epsilon/(\det F)^q,
\\
\nonumber&\phi:\mathrm{GL}^+(d){\times} \R^+ 
 \to\R^+\ \text{twice continuously differentiable},\ \forall\, F,\tilde F\in
\mathrm{GL}^+(d),\ \theta \geq 0:\\
&\label{ass-phi}
 \qquad\big|\phi(F,\theta){-}\phi(\tilde F,\theta)\big|
  \leq K\big(1 {+}|F|^{s/2}{+}|\tilde F|^{s/2}\big)|F{-}\tilde F|,
 \quad 
\\[0.3em]&\label{ass-phi+}\qquad
\pl_{FF}^2\phi(F,\theta)\le K,\ \ \ 
|\theta\pl_{F\theta}^2\phi(F,\theta)|\le K,
\ \ \
\Epsilon\le-\theta\pl_{\theta\theta}^2\phi(F,\theta)\le K,
\\[0.4em]
&\nonumber\mathscr{H}:\R^{d\times d\times d}\to\R^+
\text{ convex, continuously differentiable}, \forall\, 
  G \in\R^{d\times d\times d}:
\\
&\label{ass-H1}
\qquad
\Epsilon|G|^p\le\mathscr{H}(G) \leq K(1{+}|G|^p),
\\[0.4em]
 &\nonumber
\hat\zeta:\R^{d\times d}_\mathrm{sym} {\times}
  \R^{d\times d}_\mathrm{sym}{\times}\R  \to\R^+
 \text{ is continuous and }\forall\: (C,\DT C,\theta)\in\R^{d\times
  d}_\mathrm{sym}{\times}\R{\times}\R^{d\times d}_\mathrm{sym}:
\\
&\label{ass-hat-zeta}
 \qquad\hat\zeta(C,\cdot,\theta):\R^{d\times d}_\mathrm{sym}\to\R^+
\text{ quadratic (cf.\ \eqref{hat-zeta-special})},\ \
\alpha|\DT C|^2\le\hat\zeta(C,\DT C,\theta)\le K|\DT C|^2,
\\[0.4em]
&\label{ass-K}
\bbK:\R\to\R^{d\times d}\ \text{ is continuous, uniformly positive
  definite, and bounded},
\\[0.4em]
&\label{ass-g}
  \, g \in L^2(Q;\R^d),\quad f\in L^2( \SNeu ;\R^d),
  \quad \kappa>0,
\\[0.3em]
&\label{ass-IC}
y_0\in \calYid :=
\set{ y\in W^{2,p}(\varOmega;\R^d)}{y|_{\GDir}=\text{identity} },
 \quad \mathrm{det}(\nabla y_0)\ge\Epsilon,
\\[0.3em]
&\label{ass-BC}
\theta_\flat \in  L^1(\varSigma),\quad \theta_\flat\geq 0, \quad
\theta_0 \in L^1(\varOmega), \quad \theta_0\geq 0, \quad
\psi(\nabla y_0,\theta_0) \in L^1(\varOmega),
\end{align}
\end{subequations}
where $\mathrm{GL}^+(d)$ denotes the set of matrices in $\R^{d\times
  d}$ with positive determinant.  The last assumption in  
\eqref{ass-phi+} means that $c_\mathrm{v}$ together with $c_\mathrm{v}^{-1}$ are
bounded, which is a major restriction. However, it allows for a
 rather simple estimation in Lemma~\ref{est++}; for alternative,
more general 
situations dealing with 
increasing $c_\mathrm{v}(\cdot)$ we refer to \cite[Sec.\,8.3]{KruRou19MMCM}. 

The function $w=\mathfrakw(F,\theta)$ defined in \eqref{def-of-w}
satisfies $\mathfrakw(F,0)=0$ by \eqref{ansatz}. Moreover, we have
$\pl_\theta \mathfrakw(F,\theta)=-\theta
\pl_\theta^2\phi(F,\theta)$. Hence assumption \eqref{ass-phi+}
implies, for all $ F\in \mathrm{GL}^+(\R^d) $, 
the two-sided estimates
\begin{equation}
\label{eq:mfw.estim}
 \begin{aligned} 
  &\Epsilon \theta \leq \mathfrakw(F,\theta) \leq K\theta \quad 
  \text{for all } \theta\geq 0. 
 \\
  &\Epsilon |\theta_1{-}\theta_2| \leq |\mathfrakw(F,\theta_1) {-}
 \mathfrakw(F,\theta_2)| \leq K |\theta_1{-}\theta_2| \quad 
  \text{for all } \theta_1,\theta_2\geq 0.  
 \end{aligned}
\end{equation}
 
The assumptions (\ref{ass}b,c) make the thermomechanical coupling
through $\phi$ rather weak 
in order to allow for a simple handling of the mechanical part
independently of the temperature. These restrictive assumptions are
needed for our specific and simple way of approximation method
rather than with the problem itself. E.g.\ the assumption in
\eqref{ass-phi} is used to facilitate the estimate
\eqref{est-of-dtphi/dt}, which allows us to control the difference
between $\int_\varOmega(\nabla y^k,\theta) \d x $ and $\int_\varOmega(\nabla
y^{k-1},\theta) \d x $ in terms of $\calM(y^k)$, $\calM(y^{k-1})$, and
$\| \nabla y^k{-}\nabla y^{k-1}\|_{L^2}^2$. Moreover, after having
derived uniform bounds on $|\nabla y^k|$ it will be exploited to show
that the thermo-coupling stress $\pl_F\phi$ is bounded. 
Finally, (\ref{ass}d,h) makes the stored energy finite at time $t=0$.

It will be important that $\pl_F^{}\phi(F,\theta)$ vanishes for
$\theta=0$ (which follows from \eqref{ansatz}), so that
temperature stays non-negative if $\theta_0\ge0$ and
$\theta_\flat\ge0$, as assumed.

We now state our main existence results, which will be
proved in the following Sections~\ref{se:TimeDisc} to 
\ref{se:Limit}. The method will be constructive, avoiding
non-constructive Schauder fixed-point arguments, however some
non-constructive attributes such as selections of converging
subsequences will remain. More specifically, the proof is obtained by
first making the a priori estimate for time-discretized solutions in,
see Proposition~\ref{pr:FirstAprEstim}, and then deriving an existence
result for time-continuous solutions of an $\eps$-regularized problem,
see Proposition~\ref{pr:tau-0}. Finally, Proposition \ref{prop2}
provides convergence for $\eps \to 0$.

\begin{theorem}
   [Existence of energy-conserving weak solutions]
\label{thm:MainExist}
Assume that the conditions \eqref{ass} hold.  The original
initial-boundary-value problem \eqref{system}--\eqref{BC}--\eqref{IC}
with $\mathcal{K}$ from \eqref{K-pull-back} possesses at least one
weak solution $(y,\theta)$ in the sense of Definition~\ref{def}.  In
addition, these solutions satisfy $\nabla\theta\in L^r(Q;\R^d)$ for
all $1\le r<(d{+}2)/(d{+}1)$, the mechanical energy balance
\eqref{mech-engr}, and the total energy balance \eqref{total-engr}.
\end{theorem}

As mentioned in the introduction, a lot of publications are devoted
to the simpler isothermal viscoelasticity at \finite strain, yet, in
the multi-dimensional case, they do not satisfy all the necessary
physical requirements. It is therefore worthwhile to present a
version of our existence result by restricting it to this simpler
case, for which a lot of assumptions are irrelevant or simplify. In
particular, \eqref{ansatz} simplifies as $\psi(F,\theta)=\varphi(F)$.
Of course, our theory only works because we are using a non-degenerate
second-grade material, where $\calH(y):=\int_\varOmega \mathscr H(\nabla^2 y) \,\d
x $ generates enough regularity to handle the geometric and physical
nonlinearities. To the best of the authors knowledge, even the
following result for isothermal viscoelasticity is new.

A similar regularization approach to isothermal \finite-strain
viscoelasticity was considered in \cite{FriKru18PNLV}, where the
$\calH(y)$ is multiplied with a small parameter that vanishes slower
than the loading. Hence, the authors are able to show that their
solutions are sufficiently close to the identity which allows them to
exploit a simpler Korn's inequality obtained by a perturbation
argument.  Hence, to the best of the author's knowledge the following
result is the first that allows for truly \finite strains.

\begin{corollary}[Viscoelasticity at constant temperature]
   \label{Cor:ViscoElast}
   Let $\varphi$ satisfy \eqref{ass-psi}, and let
   \textrm{(\ref{ass}d-e,g-h)} be satisfied with
   $\hat\zeta=\hat\zeta(C,\DT C)$ and with $\psi=\varphi$. Then, the
   initial-boundary-value problem
   \eqref{momentum-eq}--\eqref{BC1}--\eqref{IC} (with $\theta$
   ignored) possesses at least one weak solution $y$ in the sense that
   the integral identity \eqref{momentum-weak+} holds. In addition,
   the mechanical energy balance \eqref{eq:MechEnergBal} holds with
   $\xi=\xi(F,\DT F)$ and without the last term involving
   $\pl_F^{}\phi$.
\end{corollary}

Before going into the proof of our main result, we show that our
conditions are general enough for a series of nontrivial
applications:

\begin{example}[Classical thermomechanical coupling]
\label{ex:Classical}\upshape 
The classical example of a free energy in
thermomechanical coupling is given in the form 
\begin{align}
\psi(F,\theta) = \varphi(F)- a(\theta)\,\varphi_1 (F)
   +c\theta(1{-}\log\theta),
\label{ansatz-special}\end{align}
i.e.\ $\phi(F,\theta)$ involves a term in the product form
$-a(\theta)\varphi_1(F)$.  For the purely mechanical part we may take
the polyconvex energy $\varphi(F)= c_1 |F|^s + c_2/(\det F)^q$ for
$\det F>0$ and $\infty $ otherwise. For the thermomechanical coupling
we obtain $c_\text{v}(F,\theta) = - \theta \pl_{\theta\theta}^2
\psi(F,\theta)= c + a''(\theta)\varphi_1(F)$, thus to have positivity
of the heat capacity $c_\text{v}$, we assume $a''(\theta)\geq 0$ and
$\varphi_1(F)\geq 0$. Moreover, we have
\[
w=\mathfrakw(F,\theta) = c\,\theta + \big(\theta
a'(\theta){-}a(\theta)\big) \varphi_1(F) \quad \text{and}\quad 
\pl_F\phi(F,\theta)=a(\theta)\varphi_1'(F)
. 
\]
Thus, we see that all assumptions in \eqref{ass} can easily be
satisfied, e.g.\ by 
choosing $a(\theta)=(1{+}\theta)^{-\alpha}$ with $\alpha>0$, which is
smooth bounded and convex, and taking any $\phi_1 \in
C^2_\text{c}(\R^{d\times d})$.  
\end{example}

\begin{example}[{\it Phase transformation in shape-memory alloys}]
\label{ex:PhaseTrafo}
\upshape 
An interesting example of a free energy $\psi$ occurs in modeling
of austenite-martensite transformation in so-called shape-memory 
alloys: 
\begin{align*}
\psi(F,\theta)=(1{-}a(\theta))\varphi_{_\mathrm{A}}(F)
+a(\theta)\varphi_{_\mathrm{M}}(F)+\psi_0(\theta).
\end{align*}
cf.\ e.g.\ \cite{Roub04MMES} and references therein.
Here $a$ denotes the volume fraction of the austenite versus martensite 
which is supposed to depend only on temperature. Of course, this is 
only a rather simplified model.
%
For, $\psi_0(\theta)=
c\theta(1{-}\log\theta)$
it complies with the ansatz \eqref{ansatz-special} with $\varphi(F)=\varphi_{_\mathrm{A}}(F)$
and $\varphi_1(F)=
\varphi_{_\mathrm{M}}(F){-}\varphi_{_\mathrm{A}}(F)$.
The heat capacity then reads as
\begin{align*}c_\mathrm{v}(F,\theta)=\theta a''(\theta)
[\varphi_{_\mathrm{A}}{-}\varphi_{_\mathrm{M}}](F)-\theta\psi_0''(\theta).
\end{align*}
To ensure its positivity, $\psi_0$ is to be strictly concave in
such a way that $\psi_0''(\theta)\le K/\theta$ and 
then
$\inf_{(F,\theta)}\theta a''(\theta)[\varphi_{_\mathrm{A}}{-}\varphi_{_\mathrm{M}}](F)+K>0$
is to (and can) be ensured by suitable modeling assumptions. 
\end{example}

\begin{example}[{\it Thermal expansion}]
  \upshape 
  Multiplicative decomposition $F=F_\mathrm{el}F_\mathrm{th}$ with the
  ``thermal strain'' $F_\mathrm{th}=\bbI/\mu(\theta)$ and the elastic
  strain $F_\mathrm{el}$ which enters the elastic part of the stored
  energy $\varphi$. This leads to
\begin{align}
 \label{psi-expansion}
  \psi(F,\theta)=\beta(\theta)\varphi(F_\mathrm{el})+\phi(\theta)
  =\beta(\theta)\,\varphi\big(\mu(\theta)F\big)-\phi(\theta).
\end{align}
Unfortunately, \eqref{psi-expansion} is inconsistent with the ansatz
\eqref{ansatz} because the contribution $\varphi$ which has been
important for our analysis due to uniform coercivity, cannot be
identified in \eqref{psi-expansion}.
\end{example}

\section{A few auxiliary results}
\label{subs:AuxRes}

In this subsection we provide a series of auxiliary results that are
crucial to tackle the difficulties arising from \finite-strain theory.
First we show how the theory developed by Healey and Kr\"omer
\cite{HeaKro09IWSS} which allows us to show that a bound for the
elastic energy $\calM(y,\theta)$ provides lower bounds on the $\det
\nabla y$. This can then be used to establish the validity of the
Euler-Lagrange equations and useful $\lambda$-convexity result, which
is needed for obtaining optimal energy estimates. Second we provide a
version of Korn's inequality from Pompe \cite{Pomp03KFIV} that allows
us to obtain dissipation estimates via $\calD(y,\DT y,\theta) \geq c_0
\| \DT y\|^2_{H^1(\varOmega)} $. Finally, in Section
\ref{suu:ChainRule} we provide abstract chain rules as derived in
\cite[Sec.\,2.2]{MiRoSa13NADN} that allows us to derive energy
balances like \eqref{eq:MechEnergBal} from the corresponding weak
equations.

\subsection{Local invertibility  and Euler-Lagrange equations}
\label{suu:HealeyK}

A crucial point in \finite-strain theory is the blow-up of the energy
density $\psi(F,\theta)$ for $\det F\searrow 0$. Thus, it is desirable 
to find a suitable positive lower bound for $\det \nabla y(t,x)$. 
The following theorem is an adaptation of the result in
\cite[Thm.\,3.1]{HeaKro09IWSS}. 

\begin{theorem}[Positivity of determinant] 
\label{th:HealeyKromer}
Assume that the functional $\calM:W^{2,p}(\varOmega;\R^d)
  \to \R_\infty$ satisfies the assumption \eqref{ass-psi} and
  \eqref{ass-H1}.  Then, for each $C_M>0$ there exists a
  $C_\mathrm{HK}>0$ such that all $y\in \calYid$ with $\calM(y)\leq C_M$
  satisfy
    \begin{equation}
      \label{eq:HealKroEstim}
    \| y\|_{W^{2,p}} \leq C_\mathrm{HK},\ \  \| y\|_{C^{1,1-d/p}} \leq
    C_\mathrm{HK},\ \   \det \nabla y(x) \geq \frac1{C_\mathrm{HK}}, 
    \ \  
     \|(\nabla y)^{-1}\|_{C^{1-d/p}} \leq C_\mathrm{HK}. 
    \end{equation}
\end{theorem}
\begin{proof} We give the full proof, since our mixed boundary
  conditions are not covered in \cite{HeaKro09IWSS}. From
  $\calM(y)\leq C_M$ and the coercivities of $\varphi$ and $\mathscr
  H$ we obtain $\det\nabla y\geq 0$ a.e.\ in $\varOmega$ and the a priori
  bounds
\[
\| \nabla y\|_{L^s} + \| \big(\det(\nabla y)\big)^{-1}\|_{L^q} + \|
   \nabla^2 y\|_{L^p} \leq C_M^{(1)}.
\]
Together with the Dirichlet boundary conditions in $\calYid$ we obtain an a priori
bound for $y$ in $W^{2,p}(\varOmega;\R^d)$ and hence also in 
$C^{1,\lambda}(\varOmega;\R^d)$, where $\lambda = 1-d/p>0$. This proves
the first two assertions.

In particular, the function $\delta:x \mapsto \det(\nabla y(x))$ is
H\"older continuous as well with $\|\delta\|_{C^\alpha}\leq
C_M^{(2)}$.  Since $\varOmega$ is a bounded Lipschitz domain, there exist
a radius $r_*>0$ and a constant $\alpha_*>0$ such that for all $x\in
\ol\varOmega$ the sets $B_{r_*}(x)\cap \ol\varOmega$ contains an interior
cone $C_{x}=\big\{ x{+}z \:\big|\: 0<|z|<r_*,\ \frac1{|z|}z \in A(x)
\big\}$ where the set $A(x)\subset \mathbb S^{d-1}$ of cone directions
has a surface measure $\int_{A(x)} 1 \d S \geq \alpha_*$. Thus, using the
H\"older continuity
\[
\delta(y) \leq \delta (x) + C^{(2)}_M |x {-}y|^\lambda \quad \text{
  for all } x,y\in \ol\varOmega,
\]
we can estimate as follows:
\begin{align*}
\big(C_M^{(1)}\big)^q &\geq \int_\varOmega \frac1{\delta(y)^q} \d y \geq
\int_{\varOmega\cap B_{r_*}(x)} 
\frac1{\big(\delta(x)+C^{(2)}_M|x{-}y|^\alpha \big) ^q} \d y\\
&\geq  \int_{\omega\in A(x)} \int_{r=0}^{r_*} 
\frac{r^{d-1}\:\d r}{\big(\delta(x)+C^{(2)}_M r^\alpha \big) ^q}  \d
\omega \ \geq  \ \frac{\alpha_*}{2^q} \int_{r=0}^{r_*} 
\frac{r^{d-1}\:  \d r}{\max\{\delta(x)^q,(C^{(2)}_M r^\alpha)^q\} }
\\
& \geq c^{(3)}_M
\min\{\delta(x)^{-q},\delta(x)^{-(q{-}d/\lambda)}  \} \ = \ 
\frac{c^{(3)}_M}{\max \{\delta(x)^q ,\delta(x)^{q{-}d/\lambda} \}} ,
\end{align*}
where in the last estimate we crucially used the assumption
$q>pd/(p{-}d)$ which implies $\lambda q>d$. Since in the last
expression both
exponents of $\delta(x)$ are positive, we obtain the explicit lower
bound 
\[
\det\nabla y(x)=\delta(x) \geq \min\Big\{ \big(c^{(3)}_M\big)^{1/q}/C^{(1)}_M,
 \big(c^{(3)}_M/(C^{(1)}_M)^q\big)^{\lambda/(\lambda q-d)}   \Big\},
\]
which gives the third assertion in \eqref{eq:HealKroEstim}.

The last assertion follows via the implicit function theorem.
\end{proof}

The most important part of the above result is that the determinant of
$\nabla y$ is bounded away from $0$. Hence, the function $f \mapsto
\varphi(F)$, which is blows up for $\det F \searrow  0$, is evaluated
only in a compact subset of $\mathrm{GL}^+(d) \subset \R^{d \times d}$ such that
$\pl_F\phi$ and $\pl^2 \varphi$ exist.  Again following
\cite[Cor.\,3.3]{HeaKro09IWSS} we obtain the G\^ateaux
differentiability of $\calM$ and as well as a useful
$\Lambda$-semiconvexity result. 

\begin{proposition}[G\^ateaux derivative and $\Lambda$-semiconvexity] 
\label{pr:Gat.Cvx}
Assume that $\calM$ satisfies \eqref{ass-psi} and
\eqref{ass-H1}. Then, in each point $y \in \calYid$ with
$\calM(y)<\infty$ the G\^ateaux derivative in all directions $h\in
\calY_0:=\bigset{v\in W^{2,p}(\varOmega)}{v|_{\GDir}}$ exists and has the form
\begin{equation}
    \label{eq:M.Gateaux}
    \rmD \calM(y)[h]= 
    \int_\varOmega \Big( \rmD \mathscr H(\nabla^2y)\vdots\nabla^2 h +
    \pl_F \varphi(\nabla y) : \nabla h \Big) \d x  
\end{equation}
Moreover, for each $C_M>0$ there exists $\Lambda(C_M)>0$ such that for all
$y^{(1)},y^{(2)} \in \calYid$ with $\calM(y^{(j)}) \leq C_M $ and
$\|\nabla y^{(1)} - \nabla y^{(2)}\|_{L^\infty} \leq 1/C_G$ we have
$C_G$ convexity 
\begin{equation}
  \label{eq:M.LambdaCvx}
  \calM(y^{(2)}) \geq \calM(y^{(1)}) 
        + \rmD \calM(y^{(1)}) [y^{(2)} {-}y^{(1)}] 
      - \Lambda(C_M) \| \nabla y^{(2)} {-} \nabla y^{(1)} \|_{L^2}^2. 
\end{equation}
\end{proposition} 
\begin{proof} We decompose $\calM=\calH+\Phi_\mathrm{el}$, see \eqref{eq:Energ.M}.
The differentiability of the convex functional $y \mapsto \calH(y)$ on
$W^{2,p}(\varOmega;\R^d)$ is standard and
follows from \eqref{ass-H1}. For treating $\Phi_\mathrm{el}$  
we use the embedding $W^{2,p}(\varOmega) \subset C^{1,\lambda}(\varOmega)$
and exploit the result $\det \nabla y(x)\geq 1/C_\mathrm{HK}$ from Theorem
\ref{th:HealeyKromer}. For all $h\in W^{2,p}_{\GDir}(\varOmega;\R^d)$ we
find a $t_*>0$ such that $\det\big( \nabla(y{+}th)(x)\big) >
1/(2C_\mathrm{HK}) $  for all $t\in [-t_*,t_*]$ and all $x \in
\varOmega$. Hence, 
\[
\rmD \Phi_\mathrm{el}(y) [h] = \lim_{t\to 0} \frac1t \big(\Phi_\mathrm{el}(y{+}t h)
-\Phi_\mathrm{el}(y)\big) = \lim_{t\to 0} \int_\varOmega  \frac1t \big( 
\varphi(\nabla y{+}t\nabla h) - \varphi(\nabla y)\big)
\d x, 
\]
and the limit passage is trivial as the convergence in the integrand
is uniform. 

To derive \eqref{eq:M.LambdaCvx} that the convexity of $\mathscr H$
implies  
\[
\calH(y^{(2)}) \geq 
\calH(y^{(1)}) + \int_\varOmega \rmD\mathscr H(\nabla^2
y^{(1)})\vdots\big (\nabla^2 y^{(2)} - \nabla^2 y^{(1)}\big) \,\d x.
\]  
To treat the functional $\Phi_\mathrm{el}$ 
we apply Theorem
\ref{th:HealeyKromer} to $y^{(1)}$ and $y^{(2)}$, which implies the
pointwise bounds
\[
| \nabla y^{(j)}(x)| \leq C_\mathrm{HK} \quad \text{ and } \quad 
\det \nabla y^{(j)}(x) \geq 1/C_\mathrm{HK}. 
\]
Clearly there is a $\delta >0$ such that all 
\begin{align*}
\forall\, F_1,F_2 \in &\R^{d \times d}\ \forall\, s\in [0,1]: \\ 
 &\left. \begin{array}{c} |F_1|,|F_2|\leq C_\mathrm{HK}, \\
   \det F_1, \det F_2 \geq 1/C_\mathrm{HK}  \end{array} \right\} \ 
    \Longrightarrow  \ \det\big((1{-}s)F_1 + s F_2\big)\geq
     1/(2C_\mathrm{HK}). 
\end{align*}
This we denote by $-\Lambda_*$ the minimum of smallest eigenvalue of
of the matrices $\pl_F^2 \varphi(F)$ where $F \in \R^{d \times d} $
runs through the 
compact set given by $ |F|\leq C_\mathrm{HK}$ and $\det F\geq
1/(2C_\mathrm{HK})$. Hence, assuming $\|\nabla y^{(2)}{-}\nabla
y^{(2)}\|_{L^\infty} \leq \delta $  we find 
\begin{align*}
&\Phi_\mathrm{el}(y^{(2)}) - \Phi_\mathrm{el}(y^{(1)})-\rmD \Phi_\mathrm{el}(y^{(1)})[   y^{(2)}
{-}y^{(1)}] \\
&= \int_\varOmega \Big( \varphi(\nabla y^{(2)}) -
  \varphi(\nabla  y^{(1)}) - \pl_\varphi(y^{(1)}):(\nabla y^{(2)}{-}
  \nabla y^{(1)})  \Big) \d x \\
& =\int_\varOmega \frac12 \int_{s=0}^1 \pl_F^2\varphi\big((1{-}s)\nabla y^{(1)} 
          {+} s \nabla y^{(2)}\big) 
  \big[\nabla y^{(2)}{-} \nabla y^{(1)}, 
       \nabla y^{(2)}{-} \nabla y^{(1)}\big]  \d s \d x \\
& \geq -\frac{\Lambda_*}2 \int_\varOmega |  \nabla y^{(2)}{-} 
                   \nabla y^{(1)} |^2 \d x .  
\end{align*}
This establishes the result with $\Lambda(C_M):= \max\{ C_\mathrm{CK},
1/\delta, \Lambda_*/2\}$.  
\end{proof}

\subsection{A generalized Korn's inequality}
\label{suu:Korn}

The following result will be crucial to show that the nonlinear
viscosity depending on $ F=\nabla y$ really controls the $H^1$ norm of
of the rate $\DT y$. It relies on Neff's 
generalization \cite{Neff02KFIN} of the Korn inequality, in the
essential improvement obtained by Pompe \cite{Pomp03KFIV}.

\begin{theorem}[Generalized Korn's inequality]
\label{th:Korn}
For a fixed $\lambda\in {]0,1[}$ and positive constants $K>1$ 
define the set 
\[
\mathsf F_K:= \bigset{ F\in C^\lambda(\varOmega;\R^{d \times d})} 
  { \| F\|_{C^\lambda} \leq K,\  \min_{x\in \varOmega} \det F(x) \geq 1/K }.
\]
Then, for all $K>1$ there exists a constant $c_K>0$ such that for all
$F\in \mathsf F_K$ we have 
\begin{equation}
  \label{eq:KornIneq}
 \forall\, v \in H^1(\varOmega;\R^d) \text{ with } v|_{\GDir} =0: \  
 \int_\varOmega \big| F^\top\nabla v{+}(\nabla v)^\top F\big|^2 \d x \geq c_K
 \| v\|^2_{H^1}. 
\end{equation}
\end{theorem}
\begin{proof} 
In  \cite[Thm.\,2.3]{Pomp03KFIV} it is shown that \eqref{eq:KornIneq}
holds for any given $F\in \mathsf F_K$. Let us denote by $c(F)>0$ the
supremum of all possible such constants for the given $F$. 
By a
perturbation argument it is easy to see that the mapping $F\mapsto
c(F)$ is continuous with respect to the $L^\infty$ norm in
$C^0(\ol\varOmega;\R^{d \times d})$. Since $\mathsf F_K$ is a compact
subset of $C^0(\ol\varOmega;\R^{d \times d})$ the infimum of $c$ on
$\mathsf F_K$ is attained at some $F_*\in \mathsf F_K$ 
by Weierstra\ss{}' extremum principle. Because of $c(F)\geq c(F_*)$,  
we conclude that \eqref{eq:KornIneq} holds with $c_K=c(F_*)$.  
\end{proof}

We emphasize that estimate \eqref{eq:KornIneq} is not valid if $F$ is not
continuous, see \cite[Thm.\,4.2]{Pomp03KFIV}. This shows that without
the in $W^{2,p}$ is crucial to control the rate of the strain $\nabla
\DT y$, which is necessary to handle the thermomechanical coupling.
The following corollary combines Theorems \ref{th:HealeyKromer} and
\ref{th:Korn}, by using the compact embedding $W^{2,p}(\varOmega;\R^d) \subset
C^{1,\lambda}(\varOmega;\R^d)$.  

\begin{corollary}[Uniform generalized Korn's inequality on sublevels]
  \label{cor:GenKorn} Given any $C_M>0$ there exists a $\ol c_K>0$ such
  that for all $y \in \calYid$ with $\calM(y)\leq C_M$ we have the
  generalized Korn inequality 
  \begin{equation}
    \label{eq:GenKorn}
    \forall\, v \in H^1(\varOmega;\R^d) \text{ with } v|_{\GDir} =0: \  
 \int_\varOmega \big| (\nabla y)^\top\nabla v{+}(\nabla v)^\top \nabla
 y\big|^2 \d x \geq  \ol c_K
 \| v\|^2_{H^1}.  
  \end{equation}
\end{corollary}

\subsection{Chain rules for energy functionals} 
\label{suu:ChainRule}         

Abstract chain rules for energy functionals $\calJ : X\to
\R_\infty:=\R{\cup} \{\infty\}$ on a Banach
space concern the question under what conditions for an absolutely
continuous curve $z:[0,T]\to X$ the composition $t \mapsto 
\calJ(z(t))$ is absolutely continuous and satisfies 
$\frac\d{\d t} \calJ(z(t))= \langle \varXi(t), \DT z(t)\rangle $ for
$\varXi\in \overline\pl \calJ(z(t))$, where $\overline\pl$ denotes a
suitable subdifferential. In particular, this implies 
\[
\calJ(z(t_1)) = \calJ(z(t_0)) + \int_{t_0}^{t_1}  \langle \varXi(t), \DT
z(t)\rangle \d t \quad \text{for } 0\leq t_0<t_1\leq T.
\]

The case that $X$ is a Hilbert space and $\calJ$ is convex and
lower semicontinuous goes back to 
\cite[Lem.\,3.3]{Brez73OMMS}, see also \cite[Lemma 4.4]{Barb10NDEM}:

\begin{proposition}
  [Chain rule for convex functionals in a Hilbert space] 
\label{pr:ChainRule1}
Let $X$ be a Hilbert space and $\calJ : X\to
\R_\infty:=\R{\cup} \{\infty\}$ a lower semicontinuous and convex
functionals. If the functions $z:[0,T]\to X$ and $\varXi:[0,T]\to X^*$
satisfy
\[
    z \in H^1([0,T]; X), \quad \varXi\in L^2([0,T];X^*), \text{and} \quad  
\varXi(t)\in \pl\calJ(z(t)) \text{ a.e.\ in }[0,T],
\]
where $\pl\calJ$ denotes the convex subdifferential, then 
\[
t\mapsto \calJ(z(t)) \text{ lies in }W^{1,1}(0,T) \quad \text{and}
\quad \frac\d{\d t} \calJ(z(t))= \langle \varXi(t), \DT z(t)\rangle 
 \text{ a.e.\ in }[0,T].
\]
\end{proposition}

A first generalization to Banach spaces $X$ with separable dual $X^*$ is given
in \cite[Prop.XI.4.11]{Visi96MPT}. We provide a slight generalization
of the results in \cite[Sec.\,2.2]{MiRoSa13NADN} that work for
arbitrary reflexive Banach spaces and include also certain nonconvex
functionals. The functional $\calJ$ is called \emph{locally
semiconvex}, if for all $z$ with $\calJ(z)< \infty$ there
exists a $\Lambda=\hat\Lambda(z)\geq 0$ and a balls $B_r(z)=\set{\hat z\in
  X}{\|\hat z{-}z\|_X\leq r }$ with $r=\hat r(z)$ the restriction
$\calJ|_{B_r(z)} $ is $\Lambda$-semiconvex, viz.
\[
\forall\, z_0,z_1 \!\in\! B_r(z)\ \forall \,s \!\in\! [0,1]:\ 
\calJ\big((1{-}s)z_0+sz_1\big) \leq (1{-}s) \calJ(z_0) + s \calJ(z_1)
+\frac{\Lambda}2 (s{-}s^2) \|z_1{-}z_0\|_X^2.
\] 
By $\ol\pl\calJ$ we denote the Fr\'echet subdifferential which is
defined by 
\[
\ol\pl\calJ(z)=\bigset{\varXi \in X^*}{ \calJ(\hat z) \geq \calJ(z)
  +\langle \varXi,\hat z{-}z\rangle -2\hat\Lambda(z)\|\hat z{-}z\|_X^2
  \text{ for }\hat z \in B_{\hat r(z)}(z) }. 
\]
The next results follows by a simple adaptation of the proof of
\cite[Prop.\,2,4]{MiRoSa13NADN}.

\begin{proposition}
 [Chain rule for locally semiconvex functionals] 
\label{pr:ChainRule2} 
Consider a separable reflexive Banach space, a $q\in {]1,\infty[}$ with
$q'=q/(q{-}1)$,  and $\calJ:X\to \R_\infty$
a lower semicontinuous and locally semiconvex functional. If the
functions $z\in W^{1,q}([0,T]; X)$ and $\varXi\in L^{q'}([0,T];X^*) $ 
satisfy 
\[
 \sup\bigset{\calJ(z(t))}{t\in [0,T]} < \infty 
\quad \text{and} \quad  
\varXi(t)\in \ol\pl\calJ(z(t)) \text{ a.e.\ in }[0,T],
\]
then 
\[
t\mapsto \calJ(z(t)) \text{ lies in }W^{1,1}(0,T) \quad \text{and}
\quad \frac\d{\d t} \calJ(z(t))= \langle \varXi(t), \DT z(t)\rangle 
 \text{ a.e.\ in }[0,T].
\]
\end{proposition}
\begin{proof} The result follows by the fact that the image of $z$
  lies in $\mathrm{dom}\calJ=\set{z\in X}{\calJ(z)<\infty}$ and is compact in
  $Z$. Hence there is one $\Lambda_*<\infty$ and one $r_*>0$ such that
  provides $\Lambda_*$ semiconvexity on $B_{r_*}(z(t))$ for all $t\in
  [0,T]$. Hence, the results in the proof of
  \cite[Prop.\,2,4]{MiRoSa13NADN} can be applied when choosing
  $\omega^R(\hat z,z)=\Lambda_*\|\hat z{-}z\|_X$ and using that fact
  that all needed arguments are local and use only information of
  $\calJ$ in a neighborhood of the image of $z$.
\end{proof}

\section{Time discretization of a regularized problem}
\label{se:TimeDisc}

Before we construct solution by a suitable time-discretization, we
introduce regularizations in two points. Firstly, we add a linear viscous
damping which allows us to obtain simple a priori bounds for the
strain rate $\nabla \DT y$, because in the first steps of the 
construction we are not yet in the position to 
exploiting the generalized Korn inequality of Theorem \ref{th:Korn}. 
Secondly, we modify the creation of heat through the viscous damping,
which in the physically correct form leads to an $L^1$ source term
which can only be handled in the first steps of the construction
either. 

Hence, introducing the regularization parameter $\eps >0$ 
we consider the coupled system%
\begin{subequations}
  \label{system-reg}
\begin{align}
 \label{momentum-eq-reg}
 &\DIV \big(\sigma_\mathrm{vi}(\nabla y,\nabla\DT y,\theta) 
   +\eps \nabla\DT y 
  +\sigma_\mathrm{el}(\nabla y,\theta)-\DIV \mfhel(\nabla^2y)
  \big)+g  = 0, 
\\[0.4em]
 \label{heat-equation-reg}
 &\DT w -\DIV (\mathcal{K}(\nabla y,\theta)\nabla\theta)
  =  \xi^\mathrm{reg}_\eps(\nabla y,\nabla\DT y,\theta) 
     + \pl_F^{}\phi(\nabla y,\theta){:}\nabla\DT y
 \\
  \label{eq:syst.reg.w}
  & w=\mathfrakw(\nabla y,\theta),
 \\
  &\nonumber
   \text{with } \xi^\mathrm{reg}_\eps(F,\DT F,
\theta):=\frac{\xi(F,\DT F, \theta)}{1{+}\eps \, \xi(F,\DT F,
\theta) }, 
\end{align}
\end{subequations}
where $\mathfrakw$ is from \eqref{def-of-w} and
$\mathcal{K}$ from \eqref{K-pull-back}. This system is defined on $Q$
and is complemented with regularized boundary and initial conditions
\begin{subequations}
 \label{BC-reg}
\begin{align}
 \label{BC1-reg} \hspace*{2em} 
& 
\big(\sigma_\mathrm{vi}(\nabla y,\nabla\DT y,\theta)
  {+}  \eps \nabla\DT y
  {+}\sigma_\mathrm{el}(\nabla y,\theta) \big)\vec{n} 
   -\divS\big( \mfhel(\nabla^2y)\vec{n}\big) =f & 
 \text{on }&\varSigma_\mathrm{N}\hspace*{2em} 
\\[0.4em]
 \label{BC1A-reg}
 &y=\text{identity} \quad \text{on }\varSigma_\rmD, \hspace{6.8em}
   \mfhel(\nabla^2y) : (\vec{n}{\otimes}\vec{n})=0
 &\text{on }&\varSigma,
\\
& \label{BC2-reg} 
\mathcal{K}(\nabla y,\theta)\nabla\theta\cdot\vec{n}+\kappa\theta
=\kappa\theta_{\flat,\eps }\ \ 
 \text{ with }\ \ \theta_{\flat,\eps }
:=\frac{\theta_\flat}{1{+}\eps \theta_\flat}, 
 &\text{on } &\varSigma,
\\
& \label{IC-reg}
y(0,\cdot)=y_0\ \ \ \ \text{ and }\ \ \ \theta(0,\cdot)=\theta_{0,\eps }:=\frac{\theta_0}{1{+}\eps \theta_0}\ & \text{on }& \varOmega.
\end{align}
\end{subequations}

This system is solved by time discretization. For this we consider a
constant time step $\tau>0$ such that $T/\tau$ is an integer, 
leading to an equidistant partition of the considered time interval $[0,T]$.
(Let us emphasize, however, that a varying time-step and non-equidistant 
partitions  can be easily implemented because we will always 
consider only first-order time differences and one-step formulas.)

For time discretization of the regularized system
\eqref{system-reg}--\eqref{BC-reg}
we use the difference notation 
\[
\DELTA f^k =\frac1\tau\big(f^k- f^{k-1}\big)
\]
and define a staggered scheme, where first $y^{k-1}_{\eps\tau}$ is
updated to $y^{k}_{\eps\tau}$ while keeping $\theta^{k-1}_{\eps\tau}$
fixed, and then $\theta$ is updated implicitly by updating
$w^{k-1}_{\eps\tau} $ to $w^{k}_{\eps\tau} = \mathfrakw(\nabla
y^{k}_{\eps\tau}, \theta^k_{\eps\tau})$. More precisely, in the domain
$\varOmega$ we ask for 
\begin{subequations}
  \label{system-disc}
\begin{align}\nonumber
&-\DIV \bigg(\sigma_\mathrm{vi}\Big(\nabla y_{\eps\tau}^{k-1},
 \DELTA \nabla y_{\eps\tau}^k\, ,\theta_{\eps\tau}^{k-1}\Big)
+\eps  \DELTA \nabla  y^k
\\[-.3em]
& \label{momentum-eq-disc}
\qquad\qquad
+\sigma_\mathrm{el}(\nabla y_{\eps\tau}^k,\theta_{\eps\tau}^{k-1})
-\DIV \mfhel(\nabla^2y_{\eps\tau}^k) \bigg)
=g_\tau^k:=\frac1\tau\int_{(k-1)\tau}^{k\tau}\!\!g(t)\,\d t,
\\[0.3em]
&\nonumber
\DELTA w^k_{\eps\tau}
-\DIV (\mathcal{K}(\nabla
y_{\eps\tau}^{k-1},\theta_{\eps\tau}^{k-1})\nabla\theta_{\eps\tau}^k) 
=\xi^\mathrm{reg}_\eps(\nabla y_{\eps\tau}^{k-1},\nabla
\delta_\tau y_{\eps\tau}^k ,\theta_{\eps\tau}^{k-1}) 
\\[.3em]&\hspace*{15.6em}
+ \pl_F^{}\phi(\nabla y_{\eps\tau}^k,\theta_{\eps\tau}^k){:}
\DELTA \nabla y_{\eps\tau}^k
\label{heat-equation-disc}
\end{align}\end{subequations}
together with the discrete variant of the boundary conditions 
\eqref{BC-reg} as
\begin{subequations}
\label{BC-reg-disc}
\begin{align}\nonumber \hspace*{2em}
&\Big(\sigma_\mathrm{vi}\Big(\nabla y_{\eps\tau}^{k-1}, 
  \DELTA \nabla y_{\eps\tau}^k\, ,\theta_{\eps\tau}^{k-1}\Big)
  +\eps  \DELTA \nabla y_{\eps\tau}^k
  +\sigma_\mathrm{el}(\nabla y_{\eps\tau}^k,\theta_{\eps\tau}^{k-1})
 \Big)\vec{n}
\\ \label{BC1-reg-disc}
&\quad\qquad\qquad\qquad
  -\divS\big(\mfhel(\nabla^2y_{\eps\tau}^k) \vec n\big) =
  f_\tau^k:=\frac1\tau\int_{(k-1)\tau}^{k\tau}\!\!\!f(t)\,\d t
 &&\text{on }\GNeu, \hspace*{2em} 
\\
&y_{\eps\tau}^k=\text{identity  \qquad on } \GDir,
 \hspace*{5em}\mfhel(\nabla^2y_{\eps\tau}^k) : 
 (\vec{n} {\otimes} \vec{n})=0
&&\text{on }\ \varGamma,
\label{BC2-reg-disc}
\\
\label{BC3-reg-disc}
&\mathcal{K}(\nabla
y_{\eps\tau}^{k-1},\theta_{\eps\tau}^{k-1}) \nabla 
 \theta_{\eps\tau}^k\cdot\vec{n}+\kappa\theta_{\eps\tau}^k
 =\kappa\theta_{\flat,\eps ,\tau}^k  :=\frac\kappa\tau 
 \int_{(k-1)\tau}^{k\tau}\!\! \theta_{\flat,\eps }(t)\,\d t 
&&\text{on }\ \varGamma.
\end{align}
\end{subequations}
The main advantage is that the  boundary-value problem 
\eqref{momentum-eq-disc}, \eqref{BC1-reg-disc}, and
\eqref{BC2-reg-disc} for $y^k_{\eps\tau} $ are the Euler-Lagrange
equation of a functional, so that solutions 
can be obtained by solving the global minimization problem 
\begin{align}
\nonumber
y^k_{\eps\tau} \in \text{Arg\;\!Min} 
    \Big\{ \:\frac1\tau \calR(y_{\eps\tau}^{k-1},
    y{-}y_{\eps\tau}^{k-1}, \theta_{\eps\tau}^{k-1}) + 
  \frac{\eps}{2\tau}\|\nabla y{-}\nabla y_{\eps 
    \tau}^{k-1}\|_{L^2}^2 &
\\
\label{minimize-y}
     {}+   \varPsi(y,\theta_{\eps\tau}^{k-1})- 
 \langle\ell^k_\tau, y\rangle &\:\Big|\: y \in \calYid \:\Big\}  ,
\end{align}
where $\langle\ell^k_\tau, y\rangle=\int_\varOmega g^k_\tau{\cdot} y\d x
+ \int_{\GNeu} f^k_\tau{\cdot} y \d S$. 
Clearly, the Euler-Lagrange equation may have more solutions, however
for deriving suitable a priori bounds, we will exploit the minimizing
properties. 

Similarly, the boundary value problem \eqref{heat-equation-disc} and
\eqref{BC3-reg-disc} for $\theta^k_{\eps\tau}$, where
$y^{k-1}_{\eps\tau}$ and $ y^k_{\eps\tau}$ are given, has a
variational structure. For this we define the functions 
$\phi_\mathrm{C}(F,\theta):=\int_0^\theta \phi(F,\hat\theta)\d \hat\theta$ and
$W(F,\theta)=2\phi_\mathrm{C}(F,\theta) - \theta \phi(F,\theta)$ to obtain the
relation
\begin{equation}
  \label{eq:W.mfw}
  \pl_\theta W(F,\theta)=\mathfrakw(F,\theta) = \phi(F,\theta) -
  \theta \pl_\theta \phi(F,\theta) \ \text{ and } \ 
  \pl_\theta\pl_F \phi_\mathrm{C}(F,\theta)= \pl_F\phi(F,\theta).  
\end{equation}
With $\pl_\theta^2 W(F,\theta)=\pl_\theta \mathfrakw(F,\theta)=
-\theta\pl_\theta^2 \phi(F,\theta)  \geq \Epsilon$ we see that
$W(F,\cdot)$ is uniformly convex by assumption \eqref{ass-phi+}. 
Thus, we can obtain solutions $\theta^k_{\eps\tau}$ of
\eqref{heat-equation-disc} and \eqref{BC3-reg-disc}  via the
minimization problem 
\begin{align}\nonumber
\theta^k_{\eps\tau} \in \text{Arg\;\!Min}\Big\{& \int_\varOmega \Big(
     \frac1\tau \big( W(\nabla y^k_{\eps\tau}, \theta )-
      w^{k-1}_{\eps\tau}\theta \big)+ \frac12\nabla\theta{\cdot}
      \mathcal{K}(\nabla y_{\eps\tau}^{k-1}, 
      \theta_{\eps\tau}^{k-1})\nabla\theta \Big) \d x  \\
&\nonumber
+ \int_\varOmega \Big( {-}\xi^\text{reg}_\eps 
  (\nabla y_{\eps\tau}^{k-1},\DELTA y^k_\eps,\theta_{\eps\tau}^{k-1}) \theta 
-\pl_F \phi_\mathrm{C}( \nabla y_{\eps\tau}^k,\theta):\DELTA \nabla
y^k_{\eps\tau} \Big) \d x \\
\label{minimize-theta}
&
\hspace*{8em}
+\int_\varGamma\frac\kappa2\big(\theta{-}\theta^k_{\flat,\eps,\tau}\big)^2\,\d
S \; \Big| \; \theta \in H^1(\varOmega),\ \theta\geq 0 \:\Big\} .
\end{align}
We emphasize that this staggered scheme is constructed in a very
specific way by taking $\theta=\theta^{k-1}_{\eps\tau}$ from the
previous time step in the mechanics problem for $y^k_{\eps\tau}$, 
see \eqref{minimize-y}. For the construction of
$\theta=\theta^k_{\eps\tau}$ from the heat equation we have to use
sometimes the explicit (backward) approximations
$\theta^{k-1}_{\eps\tau} $ and sometimes the implicit (forward)
approximation $\theta^k_{\eps\tau}$.  Clearly, the former is simpler
and it is used in the heat conduction tensor $\mathcal{K}( \nabla
y_{\eps\tau}^{k-1}, \theta_{\eps\tau}^{k-1})$ and in the heat
production $\xi^\text{reg}_\eps$.  It is tempting to use the
explicit choice $\theta^{k-1}_{\eps\tau}$ also in the
thermo-mechanical coupling term $\pl_F\phi(\nabla y^{k}_{\eps\tau}
,\theta){:}\nabla \DELTA y^k_\eps$ (last term in
\eqref{heat-equation-disc}) as it would simplify the energy balance,
see Remark \ref{rm:DiscreteEnergyBal}. However, as this term does
not have a sign, we would not be able to guarantee positivity of
$\theta^k_{\eps\tau}$. Thus we are forced to use the more involved
implicit term $\theta \mapsto \pl_F\phi_\mathrm{C}(\nabla y^{k}_\eps
,\theta ) {:} \nabla \DELTA y^k_\eps$ in \eqref{minimize-theta}
instead of the simpler, linear choice $\theta \mapsto \theta \pl_F
\phi (\nabla y^{k}_{\eps\tau}, \theta^{k-1}_{\eps\tau} ){:} \nabla
\DELTA y^k_\eps$.  This choice may introduce a nonconvexity, so that
$\theta^k_{\eps\tau}$ may not be unique.

The following result states that we can obtain solutions
$(y^k_{\eps\tau},\theta^k_{\eps\tau})$ of \eqref{system-disc} and
\eqref{BC-reg-disc} by solving the minimization problems
\eqref{minimize-y} and \eqref{minimize-theta}, alternatingly. 
For notational simplicity we have written  the minimization problem
\eqref{minimize-theta} for $\theta$ with the constraint $\theta\geq
0$, however, for establishing the Euler-Lagrange
\eqref{heat-equation-disc} and \eqref{BC3-reg-disc} we need to  show
that non-negativity of $\theta$ comes even without imposing the
constraint.  This will be achieved by minimization over $\theta \in
H^1(\varOmega)$ after extending all functionals suitably for $\theta
<0$. 

\begin{proposition}[Time-discretized solutions via minimization]
  \label{pr:IncrMinim} 
  Let our assumptions \eqref{ass} be satisfied. For $N\in \N$ set
  $\tau=T/N$ and $ (y^0_{\eps\tau} , \theta^0_{\eps\tau} ) = (y_0,
  \theta_{0,\eps} ) $ as in \eqref{IC-reg}. Then, for $k=1,\ldots,N$
  we can iteratively find $(y^k_{\eps\tau},\theta^k_{\eps\tau})\in
  \calYid \times H^1_+(\varOmega)$ by solving first the incremental
  global minimization problem \eqref{minimize-y} and then
  \eqref{minimize-theta}. The global minimizers satisfy the
  time-discretized problem \eqref{system-disc}.
\end{proposition}
\begin{proof} 
       {\emph{Mechanical step:}} We first show that the
minimization problem in \eqref{minimize-y} has 
a solution for any $\theta^{k-1}_{\eps\tau}\in  H^1(\varOmega)$ with
$\theta^{k-1}_{\eps\tau} \geq 0$. 
By assumption we have $\phi(F,\theta)\geq 0 $ which implies
$ \varPsi(y,\theta)\geq \calM(y)$. Thus, the functional in the
minimization problem is coercive on $\calYid\subset W^{2,p}
(\varOmega;\R^d)$. By lower semicontinuity in $ W^{2,p}
(\varOmega;\R^d)$ we obtain the desired minimizer $y^k_{\eps\tau} \in
\calYid$ with $\calM(y^k_{\eps\tau})<\infty$. 
Hence, Theorem \ref{th:HealeyKromer} shows that the minimizer
satisfies $\det\nabla y(x) \geq \delta>0$. As in Proposition
\ref{pr:Gat.Cvx} we conclude that $y^k_{\eps\tau}$ satisfies the
Euler-Lagrange equation 
\begin{align*}
\int_\varOmega \!\Big( \pl_{\DT F}\zeta(\nabla y^{k-1}_{\eps\tau},
\nabla \DELTA y^k_\eps,\theta^{k-1}_{\eps\tau}):\nabla z + \eps \nabla\DELTA
y^k_\eps:\nabla z +\pl_F \psi(\nabla
y^k_\eps,\theta^{k-1}_{\eps\tau}):\nabla z  \Big) \d x \quad \\
{} + \rmD\calH(y_\eps^k)[z] - \langle
\ell^k_{\tau} , z\rangle \quad \text{for all }z \in \calY_0.
\end{align*}
But this gives exactly \eqref{momentum-eq-disc}, \eqref{BC1-reg-disc}, and
\eqref{BC2-reg-disc}. 

\medskip

     {\emph{Energy step:}} We now assume that 
$\theta^{k-1}_{\eps\tau}\in H^1(\varOmega)$ and $y^{k-1}_{\eps\tau},
y^k_{\eps\tau}\in \calYid$ are given with 
$\theta^{k-1}_{\eps\tau} \geq 0$ and $\calM(y^{k-1}_{\eps\tau}),
\, \calM(y^{k}_{\eps\tau}) < \infty$. With this, we show that a
variant of the 
minimization problem \eqref{minimize-theta} has a minimizer
$\theta^k_{\eps\tau}$. For this we extend the function $\phi$, which
satisfies $\phi(F,0)=0$ by assumption \eqref{ansatz}, continuously by
$\phi(F,\theta)=0$ whenever $ \theta < 0 $. As the functions 
$\mathfrakw$, $\phi_\mathrm{C}$, and $W$ are defined through $\phi$ they all
extend continuously differentiable for $\theta<0$ to the constant
value $0$. Thus, the 
integrands in \eqref{minimize-theta} are defined for all $\theta\in
\R$ and we can minimize over $\theta \in H^1(\varOmega)$, i.e.\ without
the constraint $\theta\geq 0$.  

Clearly, the extended functional is lower weakly semicontinuous on
$H^1(\varOmega)$ because of $\mathcal K \geq 0$.  To show coercivity of
the functional, we use that $\calM(y^{k-1}_{\eps\tau})<\infty$
implies $ \nabla y^{k-1}_{\eps\tau} \in L^\infty$ and $\det\nabla
y^{k-1}_{\eps\tau}(x) \geq \delta>0$. Hence, $\mathcal K$ given in
\eqref{K-pull-back} satisfies $\nabla \theta \cdot \mathcal K(\nabla
y^{k-1}_{\eps\tau} , \theta^{k-1}_{\eps\tau}) \nabla \theta \geq
\alpha_*|\nabla \theta|^2$ for some $\alpha_*>0$. Together with
the boundary integral, where $\kappa>0$ due to \eqref{ass-g}, we have
two terms that generate a lower bound $c_0\|\theta\|_{H^1(\varOmega)}^2
-C$. 

For the remaining term we observe $W(F,\theta)\geq 0$ by
construction, while $\frac1\tau w^{k-1}_{\eps\tau} $ and
$\xi^\text{reg}_\eps$ are given functions in $L^2(\varOmega)$.  Finally,
the last bulk term involving $\pl_F \phi_\mathrm{C}$ we use
\eqref{ass-phi} giving $|\pl_F \phi(F,\theta)|\leq K(1+|F|^{s/2})$ and
hence, because of $\nabla y^k_{\eps\tau} \in L^\infty(\varOmega; \R^{d
  \times d}) $, we have
\[
\big|\pl_F\phi_\mathrm{C} (\nabla y^{k}_{\eps\tau},\theta)\big| = 
\Big| \int_0^\theta
\pl_F\phi(\nabla y^{k}_{\eps\tau},\hat\theta) \dd \hat\theta \Big|
\leq C_* |\theta|.   
\]
Together with $\DELTA \nabla y^k_{\eps} \in L^2(\varOmega;\R^{d\times
  d})$ we have show that all remaining terms can be estimated from
below by $-C\|\theta\|_{L^2(\varOmega)}$. 

In summary, we conclude that the extended functional in
\eqref{minimize-theta} is weakly lower semicontinuous and and
coercive. Hence, a global minimizers $\theta_*$ exist and moreover these
minimizers solve the associated Euler-Lagrange equation as $\pl_\theta
W (F,\theta) =\mathfrakw (F,\theta)$ and $\pl_\theta \phi_\mathrm{C}
(F,\theta) = \phi  (F,\theta)$ depend continuously on $ \theta$.  

To show that all global minimizers are non-negative we test the
Euler-Lagrange equation by the negative part
$\theta_*^-:=\min\{\theta_*,0\}$ of $\theta_*$, which is still an $H^1$
function:
\begin{align*}
0&= \int_\varOmega \Big(\frac1\tau\mathfrakw(\nabla
y^k_{\eps\tau},\theta_*)\theta_*^- - \frac1\tau w^{k-1}_{\eps\tau}
\theta_*^- + \nabla\theta_*{\cdot}
      \mathcal{K}(\nabla y_{\eps\tau}^{k-1}, 
      \theta_{\eps\tau}^{k-1})\nabla\theta_*^- \Big) \d x  \\
&\quad
+ \int_\varOmega \Big( {-}\xi^\text{reg}_\eps 
  (\nabla y_{\eps\tau}^{k-1},\DELTA y^k_\eps,\theta_{\eps\tau}^{k-1}) \theta_*^- 
- \theta_*^- \pl_F \phi( \nabla y_{\eps\tau}^k,\theta_*):\DELTA \nabla
y^k_{\eps\tau} \Big) \d x 
\\[-.3em]&\hspace{19em} 
+ \int_\varGamma \big( \kappa \theta_* \theta_*^- 
 - \theta^k_{\flat,\eps,\tau}\theta_*^- \big) \d S \\
& \geq \int_\varOmega \big( 0 + p_2 + \alpha_* 
  | \nabla \theta_*^-|^2 + p_4+ 0 \big) \d x 
  + \int_\varGamma \big( \kappa (\theta_*^-)^2 + p_7\big) \d S
  \geq  \ c_0 \| \theta_*^- \|_{H^1(\varOmega)}^2 .  
\end{align*}  
In the first estimate we have used $w^{k-1}_{\eps\tau}
= \mathfrakw (\nabla y^{k-1}_{\eps\tau}, \theta^{k-1}_{\eps\tau}) \geq
\Epsilon \theta^{k-1}_{\eps\tau}\geq 0$, $\xi^\text{reg}_\eps \geq
0$, and $\theta^k_{\flat,\eps,\tau}\geq 0$ which gives the non-negativity
of $p_2$, $p_4$, and $p_7$, while the first and fifth term vanish
identically since for $\theta_*>0$ we have $\theta_*^-=0$ while for
$\theta_*<0$ we have $\mathfrakw(F,\theta_*)=0$ and
$\pl_F\phi(F,\theta_*)=0$ (here we crucially use the implicit
structure). Thus, we conclude $\theta_*^-=0$ which is
equivalent to $\theta_*\geq 0$. 

Thus, choosing $\theta^k_{\eps\tau}=\theta_*$ for any global minimizer
of the extended functional we see that it is also a global minimizer
of \eqref{minimize-theta} and that the Euler-Lagrange equations hold. 
\end{proof}

Considering discrete approximations $\big(y_{\eps
  \tau}^k\big)_{k=0,...,T/\tau}$, we introduce a notation for the
piecewise-constant and the piecewise affine interpolants defined
respectively by
\begin{align}
\nonumber
&&& \left.
\begin{aligned}
&\overline{y}_{\eps\tau}(t)= y_{\eps\tau}^k\,,\quad 
\underline y_{\eps\tau}(t)= y_{\eps\tau}^{k-1},\quad \text{and}
\\[0.3em]
&y_{\eps\tau}(t)=\frac{t-(k{-}1)\tau\!\!}\tau\  y_{\eps\tau}^k
+\frac{k\tau-t}\tau y_{\eps\tau}^{k-1} 
\end{aligned}
\  \right\} && \text{for }(k{-}1)\tau< t < k\tau,
\\[0.5em]
&&&\ \ \underline y_{\eps\tau}(k\tau) = \ol y_{\eps\tau}(k\tau) = 
y_{\eps\tau}(k\tau) = y^{k}_{\eps\tau} &&  \text{for }k=0,1,\ldots,
T/\tau. 
 \label{def-of-interpolants}
\end{align}
The notations $\theta_{\eps\tau}$, $\overline\theta_{\eps\tau}$, and
$\underline\theta_{\eps\tau}$ or $w_{\eps\tau}$ have analogous
meanings. However, with $\overline g_{\tau}(t)$ we refer to the
locally averaged loadings $ \ol g_\tau(t)= g_\tau^k$ for $t\in
{]k\tau{-}\tau, k\tau]}$ (cf.\ \eqref{momentum-eq-disc}), and
similarly for $\ol f_\tau$, $\ol \ell_\tau$ and $\ol
\theta_{\flat,\eps,\tau}$. 

The following result provides the basic energy estimates where we will
crucially use the carefully chosen semi-implicit scheme defined
through the staggered minimization problems \eqref{minimize-y} and
\eqref{minimize-theta}. Here also we will essentially rely regularizing
viscous term $\eps \Delta \DT y$, as $\calR$ cannot be used because of
the missing a priori bound for $y^k_{\eps\tau}$ in
$W^{2,p}(\varOmega;\R^d)$. 
Moreover, we will exploit the fact that we
have global minimizers in  \eqref{minimize-y} rather than arbitrary
solutions of the Euler-Lagrange equations
\eqref{momentum-eq-disc}. This latter argument works because we have
neglected inertial terms in the momentum balance
\eqref{momentum-weak+} and hence in \eqref{momentum-eq-disc}. We refer
to \cite{KruRou19MMCM} to cases where inertial effects are treated but
in the isothermal case.

\begin{proposition}[First a-priori estimates]
  \label{pr:FirstAprEstim}
  Let \eqref{ass} be satisfied, then for all $\eps >0$ there exists a
  $K_\eps>0$ such that the following holds.  For $\tau <
  1/K_\eps$ the interpolants constructed from the discrete solutions
  $(y_{\eps\tau}^k,\theta_{\eps\tau}^k) \in W^{2,p}(\varOmega;\R^d)
  \times H^1(\varOmega)$, $k=1,...,T/\tau$, obtained in Proposition
  \ref{pr:IncrMinim} satisfy the following estimates:
\begin{subequations}
  \label{est}
 \begin{align}
  \label{est-u}
  & \big\| y_{\eps\tau} \big\|_{L^\infty(I;W^{2,p} (\varOmega;\R^d))
      \,\cap\, H^1(I;H^1(\varOmega;\R^d)) }\leq K_\eps ,
 \\
  \label{est-det.nabla.y}
  & \det\big(\nabla y_{\eps\tau}(t,x)\big)\geq 1/K_\eps \quad 
    \text{ a.e.\ on } Q, 
 \\ 
  \label{est-theta}
  &\big\| \ol\theta_{\eps\tau} \big\|_{ L^2(I;H^1(\varOmega)) 
 \,\cap\, L^\infty(I;L^2(\varOmega))}\le K_\eps ,
 \\
  \label{est-w}
  &\big\| \ol w_{\eps\tau} \big\|_{ L^2(I;H^1(\varOmega))\,\cap\, 
         L^\infty (I;L^2(\varOmega)^*)}\le K_\eps ,
\\
  \label{est-w'}
  & \big\|   w_{\eps\tau}  \big\|_{C(I;L^2(\varOmega)) \,\cap\,
    L^2([\tau,T],H^1(\varOmega) )\,\cap\, 
    H^1(I;H^1(\varOmega)^*)}  \leq K_\eps , 
 \\
  \label{est-theta'}
  & \big\|   \theta_{\eps\tau}  \big\|_{C(I;L^2(\varOmega)) \,\cap\,
    L^2([\tau,T],H^1(\varOmega) )\,\cap\, 
    H^1(I;H^1(\varOmega)^*)}  \leq K_\eps , 
\end{align}
\end{subequations}
\end{proposition}
We emphasize that we did not make any smoothness assumptions for
$\theta_0$, hence the regularized initial values
$\theta^0_{\eps\tau}:=\theta_{0,\eps}$ and
$w^0_{\eps\tau}:=\mathfrakw(\nabla y_0,\theta_{0,\eps})$ are not
smooth. This explains, why we have to use the left-continuous
interpolants in \eqref{est-theta}  and \eqref{est-w} and why in
\eqref{est-w'} we have
to exclude the interval $[0,\tau]$ in $L^2([\tau,T];H^1(\varOmega))$. 
%
%
\begin{proof}
As $y^k_{\eps\tau}$ is a global minimizer, we can insert
$y=y^{k-1}_{\eps\tau}$ as testfunction in \eqref{minimize-y} to
obtain the estimate (recall $\DELTA y^k_\eps=
\frac1\tau(y^k_{\eps\tau}{-}y^{k-1}_{\eps\tau})$)  
\begin{align}
 \varPsi(y^k_{\eps\tau},\theta^{k-1}_{\eps\tau})
- \varPsi(y^{k-1}_{\eps\tau},\theta^{k-1}_{\eps\tau}) + \tau \calR(
y^{k-1}_{\eps\tau}, \DELTA y^k_\eps,\theta^{k-1}_{\eps\tau}) +
\frac{\eps\tau}2 \| \DELTA y^k_\eps\|_{L^2}^2 \leq \tau\langle
\ell^k_\tau , \DELTA y^k_\eps\rangle. 
\label{test-y-discrete+}
\end{align}
The proof will be divided into three steps.

\STEP{Step 1: Uniform energy bound.} 
Using the decomposition
$ \varPsi(y,\theta)=\calM(y)+\Phi_\mathrm{cpl}(y,\theta)$, see
\eqref{eq:Energ.M},  we can write equivalently  
\begin{align}\nonumber
&\calM(y^k_{\eps\tau}) -\calM(y^{k-1}_{\eps\tau}) + \tau \calR(
    y^{k-1}_{\eps\tau}, \DELTA y^k_\eps,\theta^{k-1}_{\eps\tau})
   + \frac{\eps\tau}2 \| \nabla \DELTA y^k_\eps\|_{L^2}^2 \\
&\leq  \tau\langle \ell^k_\tau , \DELTA y^k_\eps\rangle 
  + \int_\varOmega \big( \phi(\nabla  y^{k-1}_{\eps\tau}, \theta^{k-1}_{\eps\tau}) 
         - \phi(\nabla y^k_{\eps\tau},\theta^{k-1}_{\eps\tau})  \big) \d x . 
\label{test-y-discrete}
\end{align}
To estimate the last term use the assumption \eqref{ass-phi} on
$|\pl_F\phi(F,\theta)|$ as follows 
\begin{align}
  \nonumber
 \phi(F_1,\theta)-\phi(F_2,\theta) &\leq
   K(1{+}|F_1|+|F_2|)^{s/2}\,|F_1{-}F_2| \\
 &\leq \frac{K^2}{2\rho}
  (1{+}|F_1|+|F_2|)^{s} + \frac\rho2 |F_1{-}F_2|^2, 
\label{est-of-dtphi/dt}
\end{align}
where  $\rho>0$ is arbitrary. Choosing $\rho=\eps/(4\tau)$ and
$F_j=\nabla y^{k+j-2}_{\eps\tau}$ we can insert this into the estimate
\eqref{test-y-discrete}. Moreover we can use $\calR\geq 0$ and 
$\langle \ell^k_\tau, \DELTA y^k_\eps\rangle \leq
\|\ell^k_\tau\|_{H^{-1}} \| \DELTA y^k_\eps\|_{H^1} \leq 
\|\ell^k_\tau\|_{H^{-1}}c_\mathrm{P}\| \nabla\DELTA y^k_\eps\|_{L^2}  $ as $\DELTA
y^k_\eps \in \calY_0$. This 
leads to  
\begin{align*}
& \calM(y^k_{\eps\tau}) -\calM(y^{k-1}_{\eps\tau}) 
    + \frac{\eps\tau}2 \| \nabla \DELTA y^k_\eps\|_{L^2}^2 
 \\
& \leq  
  \frac{2\tau c_P^2}\eps \|\ell^k_\tau\|_{H^{-1}}^2 
  + \frac{\eps\tau}8\|\nabla\DELTA y^k_\eps\|_{L^2}^2 
  + \frac{2\tau K^2}\eps \int_\varOmega  \big(1{+}|\nabla y^k_{\eps\tau}| 
        {+} |\nabla y^{k-1}_{\eps\tau}|)^s \d x 
  + \frac{\eps\tau}8\|\nabla\DELTA y^k_\eps\|_{L^2}^2 .    
\end{align*}
Using the coercivity assumption \eqref{ass-phi} for $\phi$ the
second-last term can be estimated by $\calM$ again and setting
$m_k:= \calM(y^k_{\eps\tau})$ we obtain the recursive estimate 
\begin{equation}
  \label{eq:mk.mk-1}
  m_k - m_{k-1}  +\frac{\eps\tau}4 \| \nabla \DELTA y^k_\eps\|_{L^2}^2 
 \leq \tau \ol c_\eps  \|\ell^k_\tau\|_{H^{-1}}^2  +
\tau C_\eps (|\varOmega|{+}m_k{+}m_{k-1}) 
\end{equation}
with $C_\eps =2{\cdot}3^s K^2/\eps$ and $\ol c_\eps= 2
c_\mathrm{P}^2/\eps$ . In a first step we neglect the
last term on the left-hand side and obtain 
\[
\big(1{-}\tau C_\eps\big) m_k \leq \big(1{+}\tau C_\eps\big)m_{k-1} +
\ol c_\eps \tau \|\ell^k_\tau\|_{H^{-1}}^2 + \tau C_\eps |\varOmega|.
\]
We now restrict $\tau>0$ via $\tau <1/(2C_\eps)$ by choosing
$K_\eps\geq 2C_\eps$, so we can iterate
the above estimate. With \eqref{ass-IC} we have
$m_0:= \varPsi (y_0,\theta_0)<\infty$ and a
simple induction yields the discrete Gronwall-type estimate  (with
$Q_\eps= (1{+}\tau C_\eps)/(1{-}\tau C_\eps)$) 
\begin{align}
\nonumber
m_k &\leq Q_\eps^k m_0 + \frac{\tau}{1{-}\tau C_\eps} \sum_{j=1}^k Q^{k-j}
     \big(\ol c_\eps \|\ell^j_\tau\|_{H^{-1}}^2 {+} C_\eps|\varOmega|\big) \\
\nonumber
  & \leq Q^k \Big( m_0 + 2 \ol c_\eps \big(\sum_{j=1}^k \tau
  \|\ell^j_\tau\|_{H^{-1}}^2\big) + k\tau \,2C_\eps|\varOmega| \Big)\\
\label{eq:calM.bound}
&\leq 4\mathrm{e}^{2C_\eps T} \Big(  \varPsi (y_0,\theta_0) + 2 \ol c_\eps  
  \int_0^T \|\ell(s)\|^2_{H^{-1}} \d s + 2T C_\eps
  |\varOmega|\Big):=\widetilde K_\eps .
\end{align}
Using Theorem \ref{th:HealeyKromer} we obtain the desired uniform 
upper bound in \eqref{est-u} for the interpolant
$y_{\eps\tau}:I=[0,T]\to \calYid$ in
$L^\infty\big(I;W^{2,p}(\varOmega;\R^d)\big)$ as well as the lower bound
\eqref{est-det.nabla.y} for the determinant.  
%
%

    \STEP{Step 2: Dissipation bound.} 
We return to \eqref{eq:mk.mk-1}
and add all estimates from $k=1$ to $N_\tau:=T/\tau \in \N$ to obtain 
\begin{align*}
\frac\eps4 \int_Q |\nabla \DT y_{\eps\tau}|^2 \d x \d t &
 =\frac{\eps\tau}4 \sum_{k=1}^{N_\tau}  \| \nabla \DELTA y^k_\eps\|_{L^2}^2 \\
&\leq m_0-m_{N_\tau} + \tau \sum_{k=1}^{N_\tau}
   \Big(\ol c_\eps\|\ell^k_\tau\|_{H^{-1}}^2 + C_\eps(|\varOmega|{+}m_{k-1}
   {+} m_k) \Big) 
\\
&\leq  \varPsi (y_0,\theta_0) + \ol c_\eps \| \ell\|_{L^2(I;H^{-1})}^2 + C_\eps T
(|\varOmega|{+}2\widetilde K_\eps) =:\widehat K_\eps. 
\end{align*}
This provides the uniform bound for $y_{\eps\tau}$ in
$H^1(I;H^1(\varOmega;\R^d))$, and \eqref{est-u} is established. 
\\[0.3em]
%

\STEP{Step 3: Temperature bounds.} 
Testing the Euler-Lagrange equations \eqref{heat-equation-disc} and
\eqref{BC3-reg-disc} by $w_{\eps\tau}^k$ yields the identity
\begin{align}
\nonumber
 \int_\varOmega \!\Big( \frac{w_{\eps\tau}^k{-}w_{\eps\tau}^{k-1}\!\!} 
    \tau\;w_{\eps\tau}^k + \nabla w_{\eps\tau}^k {\cdot} 
   \mathcal{K} (\nabla y_{\eps\tau}^{k-1} , \theta_{\eps\tau}^{k-1} ) 
   \nabla \theta_{\eps\tau}^k  \Big) \d x 
  + \int_\varGamma \kappa\theta_{\eps\tau}^k w_{\eps\tau}^k \dd S 
\qquad 
\\[-.1em]
=\int_\varOmega h_{\eps\tau}^kw_{\eps\tau}^k\,\d x+\int_\varGamma
\kappa\theta_{\flat,\eps ,\tau}^kw_{\eps\tau}^k\,\d S 
\label{heat-L2-test}
\\[.2em]
\nonumber
\text{with } h_{\eps\tau}^k:= \xi^\mathrm{reg}_\eps(\nabla y_{\eps
  \tau}^{k-1}, \nabla \DELTA y_\eps ^k ,\theta_{\eps\tau}^{k-1}) 
   +\pl_F^{}\phi(\nabla y_{\eps\tau}^k,\theta_{\eps\tau}^k){:}
   \nabla \DELTA y_\eps ^k. \hspace*{5em}
\end{align}
Recalling $c_\mathrm{v}(F,\theta)=\pl_\theta
\mathfrakw(F,\theta)$ we obtain the chain rule 
\begin{align}
 \label{nabla-w}
 \nabla w_{\eps\tau}^k=\nabla \mathfrakw (\nabla y_{\eps\tau}^k,
 \theta_{\eps\tau}^k) =\pl_F^{} \mathfrakw (\nabla y_{\eps\tau}^k,
 \theta_{\eps\tau}^k) {:} \nabla^2y_{\eps\tau}^k
 +c_\mathrm{v}(\nabla y_{\eps\tau}^k, \theta_{\eps\tau}^k ) \nabla
 \theta_{\eps\tau}^k.
\end{align}
Moreover, we have the elementary estimate $
\frac1\tau(w^k_{\eps\tau}{-}w^{k-1}_{\eps\tau}) w^k_{\eps\tau} \leq \frac1{2\tau}\big(
(w^k_{\eps\tau})^2-(w^{k-1}_{\eps\tau})^2\big)$, and  $\theta
w=\theta\mathfrakw(F,\theta)\ge0$ by the definition of
$\mathfrakw$. Using additionally $c_\mathrm{v}(F,\theta)=-\theta
\pl_\theta^2 \phi(F,\theta) \geq \Epsilon $ (see \eqref{ass-phi+}, the
above identity \eqref{heat-L2-test} leads to 
\begin{align}
\nonumber
&\int_\varOmega \Big( \frac1{2\tau}(w_{\eps\tau}^k)^2-\frac1{2\tau}(w_{\eps\tau}^{k-1})^2
+ \Epsilon \nabla\theta_{\eps\tau}^k \cdot \mathcal{K}^{k}_{\eps\tau}
\nabla\theta_{\eps\tau}^k \Big) \d x 
\\&\ \ 
\le\int_\varOmega h_{\eps\tau}^kw_{\eps\tau}^k
-\nabla\theta_{\eps\tau}^k \cdot \mathcal{K}^{k}_{\eps\tau} b^k_{\eps\tau} 
\,\d x+\int_\varGamma \kappa\theta_{\flat,\eps ,\tau}^kw_{\eps\tau}^k\,\d S.
\label{est-heat-eq-disc}
\\
\nonumber
&\text{where } \mathcal K^{k}_{\eps\tau}= \mathcal{K} (\nabla
     y_{\eps\tau}^{k-1} ,\theta_{\eps\tau}^{k-1}) \ \text{ and } \ 
 b^k_{\eps\tau}:= \pl_F^{}\mathfrakw(\nabla y_{\eps
   \tau}^k,\theta_{\eps\tau}^k){:}\nabla^2y_{\eps\tau}^k. 
\end{align}
Using uniform bounds for $\nabla y_{\eps\tau}$ and $\det \nabla
y_{\eps\tau}$ from Step~1, the assumption \eqref{ass-K} on $\mathbb
K$,  as well as formula \eqref{K-pull-back} we find a $\varkappa_\eps$ such that
\begin{equation}
  \label{eq:calK}
  |\mathcal K^k_{\eps\tau}| \leq \varkappa_\eps \quad \text{ and } \quad 
  a\cdot \mathcal K^k_{\eps\tau} a \geq \frac1{\varkappa_\eps} |a|^2 \ 
  \text{ for all }a\in \R^d.
\end{equation}
Moreover, using $\pl_F^{}\mathfrakw= \pl_F^{} \phi - \theta
\pl_{F\theta}^2 \phi $ the assumptions \eqref{ass-phi} and
\eqref{ass-phi+} together with the uniform $L^\infty$ bound for
$\nabla y_{\eps\tau}$ we find $\|\pl_F^{}\mathfrakw( \nabla y_{\eps
   \tau}^k,\theta_{\eps\tau}^k)\|_{L^\infty} \leq C_\eps$. 
Realizing also that we have $\nabla^2y_{\eps\tau}^k$ already estimated
in $L^p(\varOmega;\R^{d\times d\times d})$ with $p\ge2$ we obtain
$\|b^k_{\eps\tau}\|_{L^2} \leq \ol C_\eps$. 
For the right-hand side  $h^k_{\eps\tau}$ of
\eqref{heat-L2-test}
we have 
\[
\| h_{\eps\tau} \|_{L^2} \leq \| \xi^\mathrm{reg}_\eps\|_{L^2} + \|
\pl_F \phi(\nabla y^k_{\eps\tau}, \theta^k_{\eps\tau})\|_{L^\infty} \|
\nabla \DELTA y^k_\eps\|_{L^2} \leq \ol C_\eps \big( 1 + \|  \DELTA
y^k_\eps\|_{H^1} \big),
\] 
where we again used the $L^\infty$ bounds for $\nabla
y^k_{\eps\tau}$. 
Finally, by definition we have  $\theta_{\flat,\eps} \in [0,1/\eps]$,
and \eqref{eq:mfw.estim} allows us to estimate $w$ by $\theta$, which
yields the boundary estimate 
\[
\Big| \int_\varGamma \theta^k_{\flat,\eps,\tau} w^k_{\eps\tau} \d S \Big|
\leq \frac1\eps \int_\varGamma K |\theta^k_{\eps\tau}| \d S \leq C_\eps  \|
\theta^k_{\eps\tau}\|_{H^1} \leq \ol C_\eps\big(
\|w^k_{\eps\tau}\|_{L^2} + \|\nabla \theta^k_{\eps\tau}\|_{L^2} \big).  
\]
Based on the above estimates and introducing the abbreviations 
\[
\gamma_k:=\|w^k_{\eps\tau}\|_{L^2}, \quad 
\Theta_k:=\|\nabla \theta^k_{\eps\tau}\|_{L^2} ,
\text{ and } \nu_k:=\|\DELTA y^k_\eps\|_{H^1}  
\]
we can estimate the right-hand side in \eqref{est-heat-eq-disc} via
\begin{align*}
\text{RHS} &\leq \ol C_\eps (1{+} \nu_k) \gamma_k + C_\eps \Theta_k +
\ol C_\eps ( \gamma_k{+}\Theta_k) \leq \ol c_\eps   \Big( 
\frac1\alpha +\nu_k^2 + \gamma_k^2 + \alpha \Theta^2\Big),  
\end{align*}
where $\alpha>0$ is arbitrary. Estimating the last term on the
left-hand side in \eqref{est-heat-eq-disc} from below by
$\frac\Epsilon\varkappa \Theta_k^2$ we may choose $\alpha=
\Epsilon/(2\varkappa \ol c_\eps)$. After multiplying
\eqref{est-heat-eq-disc} by $2\tau$ we obtain  
\begin{equation}
  \label{eq:heat-rekursion}
  \gamma_k^2 -\gamma_{k-1}^2 + \frac{\Epsilon}{2\varkappa} \Theta^2_k
  \leq \tau \hat c_\eps\big(1 +  \nu_k^2 + \gamma_k^2\big).  
\end{equation}
Arguing as in Steps~1 and 2 for \eqref{eq:mk.mk-1} and using
$\gamma_0^2 = \int_\varOmega w^0_{\eps\tau}\d x \leq K^2\int_\varOmega
\theta_{0,\eps}^2 \dd x \leq K^2|\varOmega|/\eps^2 < \infty$ (cf.\
\eqref{IC-reg}) the left-continuous interpolants $\ol\theta_{\eps\tau}$
and $\ol w_{\eps\tau}$ satisfy the a priori estimates 
\[
\Epsilon \|\ol\theta_{\eps\tau} \|_{L^\infty(I;L^2(\varOmega))} \leq 
\| \ol w_{\eps\tau} \|_{L^\infty(I;L^2(\varOmega))} = \sup_{k=0,...,N_\tau}
\gamma_k \leq K_\eps 
 \text{ and } 
\|\nabla \ol \theta_{\eps\tau}\|_{L^2(Q)}^2 = \tau \sum_{k=1}^{N_\tau} 
\Theta_k^2 \leq K_\eps.
\] 
With $\theta \leq \mathfrakw(F,\theta)/\Epsilon$ we immediately find
\eqref{est-theta} for $\ol\theta_{\eps\tau}$. 
The estimate \eqref{est-w} follows by using \eqref{nabla-w} once again. 

The uniform estimate the piecewise affine interpolant $w_{\eps\tau}$ 
in the spaces $ C(I;L^2(\varOmega)) \cap L^2( [\tau,T], H^1(\varOmega) ) $ 
follows from the previous estimates for $\ol w_{\eps\tau}$. 
Finally, we note that the time derivative
interpolant $w_{\eps\tau}$  is equal to $\DELTA w^k_\eps$ on the
intervals ${](k{-}1)\tau,k\tau[}$. We now use the Euler-Lagrange
equations \eqref{heat-equation-disc} and \eqref{BC3-reg-disc}, which
provides for $\DELTA w^k_\eps =  \frac1\tau(w^k_\eps{-}
w^{k-1}_\eps) $  the estimate
\[
\| \DELTA w^k_\eps\|_{(H^1)^*} 
 \leq  C^{\mathcal K}_\eps   \|\nabla \theta^k_{\eps\tau}\|_{L^2} +
 C^\xi_\eps + C_\eps^{\pl_F\phi}         \| \DELTA y^k_\eps\|_{H^1} +
 C^\kappa_\eps \big(\| \theta^k_{\eps\tau}\|_{H^1} +
 |\varGamma|/\eps\big). 
\]           
Squaring and summation over $k=1,\dots, N_\tau$ gives the remaining
uniform bound in \eqref{est-w'} for $\pl_t w_{\eps\tau}$ in
$L^2 \big( I; H^1(\varOmega)^* \big) $. 

Using \eqref{eq:mfw.estim} once again, we bound the
increments $\DELTA \theta^k_\eps$ via the pointwise estimate 
\begin{align*}
\Epsilon |\DELTA \theta^k_\eps| &= \frac{\Epsilon}\tau |
\theta^k_{\eps\tau} {-}\theta^{k-1}_{\eps\tau}| \leq 
 \frac1\tau|\mathfrakw(\nabla y^{k-1}_{\eps\tau}, \theta^k_{\eps\tau}) 
  {-}\mathfrakw(\nabla y^{k-1}_{\eps\tau}, \theta^{k-1}_{\eps\tau})| 
\\
 & \leq \frac1\tau |w^k_{\eps\tau} {-} w^{k-1}_{\eps\tau}| + \frac1\tau | 
   \mathfrakw(\nabla y^k_{\eps\tau}, \theta^k_{\eps\tau}) 
  {-}\mathfrakw(\nabla y^{k-1}_{\eps\tau}, \theta^k_{\eps\tau})|
 \leq |\DELTA w^k_\eps| + C_\eps |\nabla \DELTA y^k_\eps|.   
\end{align*}
Taking the $H^1(\varOmega)^*$ norm we obtain $\| \DELTA
\theta^k_\eps\|_{H^1(\varOmega)^*} \leq K_\eps \big( \| \DELTA
w^k_\eps\|_{H^1(\varOmega)^*} + \| \DELTA
y^k_\eps\|_{H^1(\varOmega)}\big) $, such that \eqref{est-theta'} follows from
\eqref{est-w'}, \eqref{est-u}, and \eqref{est-theta}. 

This finishes the proof of Proposition \ref{pr:FirstAprEstim}.  
\end{proof}

\section{The limit $\tau\to 0$  in the regularized problem}
\label{su:Limit.tau0} 

Using the above a priori estimates for the interpolants we will be
able to extract convergent subsequences. First we will observe that
the three different types of interpolants have to converge to the same
limit. Next we want to pass to the limit in the discretized weak forms
of the momentum balance and the heat equation. While most terms can be
handled by compactness arguments or 
weak-convergence methods, there is one term that needs special
attention namely the heat-source term $\xi^\mathrm{reg}_\eps$ that is
quadratic in $\nabla \DT y_\eps $. Thus, it will be a crucial step to
show strong convergence of $\DT y_{\eps\tau}$ in $L^2(I;
H^1(\varOmega))$, which can be done by passing to the limit in a suitable
discretized version of the mechanical energy balance
\eqref{eq:MechEnergBal}. In this argument 
we will use the $\Lambda$-convexity derived in Proposition
\ref{pr:Gat.Cvx} to relate the mechanical energies $\calM(y^{k-1}_{\eps\tau}) $
and $\calM(y^k_{\eps\tau})$.

With the definition \eqref{def-of-interpolants} for the three
types of interpolants, we see that the following discretized version
\eqref{system-disc+} of the momentum balance and heat equations
\eqref{system-reg} and \eqref{BC-reg} holds for the discrete
solutions constructed in Proposition \ref{pr:IncrMinim}: 
\begin{subequations}
  \label{system-disc+}
 \begin{align}
  \nonumber
 &-\DIV \! \big(\sigma_\mathrm{vi}(\nabla\underline y_{\eps\tau},
   \nabla\DT y_{\eps\tau},\underline\theta_{\eps\tau})
  +\eps \nabla\DT y_{\eps\tau}    +\sigma_\mathrm{el}(\nabla 
          \overline y_{\eps\tau},\underline\theta_{\eps\tau})
 \\
  \label{momentum-eq-disc+}
 &\hspace{16em}
      -\DIV\mfhel(\nabla^2\overline y_{\eps\tau})
       \big)    =\overline g_{\tau} \,,
\\&
\DT w_{\eps\tau} {-}\DIV\!\big(\mathcal{K}(\nabla\underline y_{\eps\tau}, 
    \underline\theta_{\eps\tau})\nabla\overline\theta_{\eps\tau}\big)
  = \xi^\mathrm{reg}_\eps(\nabla\underline y_{\eps\tau}, 
      \!\nabla\DT y_{\eps\tau},\underline\theta_{\eps\tau})
    {+}\pl_F^{}\phi(\nabla\overline y_{\eps\tau}, 
      \overline\theta_{\eps\tau}){:}\nabla\DT y_{\eps\tau},
\label{heat-equation-disc+}
\\&
\overline w_{\eps\tau}=\mathfrakw(\nabla\overline y_{\eps
  \tau},\overline\theta_{\eps\tau}), 
 \label{heat-equation-disc++}
\end{align}
\end{subequations}
to hold on $Q=[0,T]\times \varOmega$,
 while the regularized boundary conditions \eqref{BC-reg-disc} read
\begin{subequations}
  \label{BC-reg-disc+}
 \begin{align}
  \nonumber
&&&\big(\sigma_\mathrm{vi}\big(\nabla\underline y_{\eps\tau},
  \nabla\DT y_{\eps\tau},\underline\theta_{\eps\tau}\big) 
  +\eps \nabla\DT y_{\eps\tau} +\sigma_\mathrm{el}( 
  \nabla\overline y_{\eps\tau}, \underline\theta_{\eps\tau})\big)\vec{n} \\
&&&\hspace*{13em}-\divS\big(\mfhel(\nabla^2\overline y_{\eps\tau}) \vec{n}\big) 
   =\overline f_{\tau}\!\!&&\text{on }\SNeu,
 \\
  \label{BC1-reg-disc+} 
&&&\overline y_{\eps\tau}=\text{\,identity on }\ \SDir, \hspace{8em}
   \mfhel(\nabla^2\overline y_{\eps\tau}){:}(\vec{n}\otimes\vec{n})=0
   &&\text{on }\varSigma,&&
 \\
&&&\label{BC2-reg-disc+} 
\mathcal{K}(\nabla\underline y_{\eps\tau},\underline\theta_{\eps\tau})
\nabla\overline\theta_{\eps\tau}\cdot\vec{n}
+\kappa\overline\theta_{\eps\tau}=\kappa\overline\theta_{\flat,\eps ,\tau}
&&\text{on }\varSigma.
\end{align}
\end{subequations}
Here it is essential that we have to use all three types of 
interpolants, e.g.\ $\ol y_{\eps\tau},\ \underline y_{\eps\tau}$, and
$y_{\eps\tau}$. In particular, we emphasize that $t\mapsto w_{\eps\tau}(t)$ is the
piecewise affine interpolant of $\{w^k_{\eps\tau}\}_{k=0,...,N_\tau}$, 
which does not coincide with 
$t\mapsto \mathfrakw(\nabla y_{\eps\tau}(t),\theta_{\eps\tau}(t))$
except at the nodal points $t=k\tau$.

\begin{proposition}[Convergence for $\tau\to0$]\label{pr:tau-0}
Let \eqref{ass} hold, and let $\eps >0$ be fixed. Then, considering 
a sequence of time steps $\tau\to0$, there is a subsequence 
(not relabeled)  and limit functions $(y_\eps ,\theta_\eps )$ such that
\begin{subequations}
\label{conv-}
\begin{align}\label{conv-y-}
&y_{\eps\tau}\to y_\eps &&\text{weakly* in }\ L^\infty(I;W^{2,p}(\varOmega;\R^d))\,\cap\,H^1(I;H^1(\varOmega;\R^d)),&&
\\&\label{conv-theta-}
\theta_{\eps\tau}\to\theta_\eps &&\text{weakly\, in }\ \,
L^2(I;H^1(\varOmega))\,\cap\,H^1(I;H^1(\varOmega)^*).&&
\end{align}\end{subequations}
Moreover, any couple $(y_\eps ,\theta_\eps )$ obtained by this way is a 
weak solution to the regularized initial-boundary-value problem 
\eqref{system-reg}--\eqref{BC-reg}.
\end{proposition}

\begin{proof}  The proof consists of five steps. 


  \STEP{Step 1: Extraction of convergent subsequences.}  As $\eps
  >0$ is still fixed, we can exploit the a priori estimates obtained
  in Proposition \ref{pr:FirstAprEstim}, namely \eqref{est-u} and
  \eqref{est-theta'}. By Banach's selection principle, we choose
  a subsequence and some $(y_\eps ,\theta_\eps )$ such that
  \eqref{conv-} holds.  By the Aubin-Lions theorem combined with an
  interpolation, as $p>d$, we have also
\begin{subequations}
  \label{conv-strong}
 \begin{align}
  \label{conv-y-strong}
  \nabla y_{\eps\tau}&\to\nabla  y_\eps && \text{uniformly in }\ 
         L^\infty(Q;\R^{d\times d}),&&
 \\
  \label{conv-theta-strong}
  \theta_{\eps\tau}&\to\theta_\eps &&\text{strongly in }\ L^s(Q)\ \ 
   \text{ for all } s \in {[1, \min\{4,2+4/d\}[}.&&
\end{align}
\end{subequations}
Indeed, for the first result we use the continuous embedding 
 $W^{1,p}(\varOmega) \subset C^\alpha( \ol \varOmega)$ with $\alpha
 =1-d/p\in {]0,1[}$ and thus $\| \nabla y_{\eps\tau}\|_{C^\alpha} \leq
 K_0$. Moreover, \eqref{est-u} yields the H\"older estimate 
\begin{align}
\big\|\nabla y_{\eps\tau}(t_1)-\nabla y_{\eps\tau}(t_2)\big\|_{L^2(\varOmega;\R^d)}
\le K_1|t_1-t_2|^{1/2} \quad \text{for all } 
t_1,t_2\in I. 
\end{align}
While the first part of 
\eqref{est-u}
yields just 
$\|\nabla y_{\eps\tau}(t_1)-\nabla y_{\eps\tau}(t_2)\|_{W^{1,p}(\varOmega;\R^d)}
\le K_0$. By interpolation, we find $\beta\in {]0,\alpha[}$ and
$\lambda \in {]0,1[}$ such that 
  we have the interpolation 
$\|\cdot\|_{C^\beta}\le C
\|\cdot\|_{C^\alpha}^{1-\lambda}\|\cdot\|_{L^2)}^\lambda$ and conclude 
\begin{align}
\big\|\nabla y_{\eps\tau}(t_1)-\nabla y_{\eps\tau}(t_2)\big\|_{C^\beta(\bar\varOmega;\R^d)}
\le C K_0^{1-\lambda}K_1^\lambda|t_1-t_2|^{\lambda/2}.
\end{align}
Thus, the sequence $\{\nabla y_{\eps\tau}\}$ is uniformly bounded
in $C^\gamma (\ol Q)$ for $\gamma=\min\{\beta,\lambda/2\}$, and
uniform convergence follows by the Arzel\`a-Ascoli theorem. 

The convergence \eqref{conv-theta-strong} follows from
\eqref{conv-theta-} by the Aubin-Lions theorem when interpolated with 
the estimate in $L^\infty(I;L^2(\varOmega))$ which is contained implicitly in
\eqref{conv-theta-}.

Moreover, both convergences in \eqref{conv-strong} hold also for the
piecewise constant interpolants because of the estimates $\|\nabla
y_{\eps\tau}-\nabla\underline y_{\eps\tau} \|_{L^\infty (I;
  L^2(\varOmega; \R^{d\times d}))}\le K\tau^{1/2}$ (and the same also for
$\nabla\overline y_{\eps\tau}$) and $ \| \nabla \theta_{\eps\tau} -
\nabla\ol\theta_{\eps\tau} \|_{ L^\infty (I;H^1(\varOmega;\R^d)^*)} \leq
K \tau^{1/2}$.

Similarly, using the a priory estimates \eqref{est-w} and
\eqref{est-w'} for $w_{\eps\tau}$ and $\ol w_{\eps\tau}$ yields 
\begin{align}
\nonumber
& w_{\eps\tau} \weak w_\eps \quad \text{weakly in }
     L^2(I;H^1(\varOmega)) \cap H^1(I;H^1(\varOmega)^*)\\
\label{eq:w.convg}
  \ol w_{\eps\tau},\, & w_{\eps\tau}  \to w_\eps \quad \text{strongly in } 
   L^s(Q) \text{ for all } s\in {[1,\min\{4,2+4/d\} [}.
\end{align}

    \STEP{Step 2: Convergence in the mechanical equation.} 
Now the convergence in the discretized momentum balance
\eqref{momentum-eq-disc+} can be done by the above weak convergences
\eqref{conv-} because $\sigma_\mathrm{vi}$ is linear in terms of $\DT
F$ and by Minty's trick for the monotone operator induced by
$\mfhel=\mathscr{H}'$. For a reflexive Banach space $X$ and  a
hemi-continuous, monotone operator $\bfH:X\to X^*$ Minty's
trick means the implication 
\begin{equation}
  \label{eq:Minty}
\left.\begin{aligned} 
\bfH(u_\tau)=b_\tau,\ \ \quad & u_\tau\weak u \text{ in }X,\\ 
 b_\tau\weak b \text{ in }X^*, \quad & \langle b_\tau,u_\tau\rangle \to \langle
 b,u\rangle  \ \ 
\end{aligned} \right\}\quad \Longrightarrow\quad \bfH (u)=b. 
\end{equation}
We apply this for $\bfH$ defined by $\langle \bfH(y),z\rangle = \int_Q
\mfhel(\nabla^2y(t,x))\Vdots \nabla^2 z(t,x)\d x \d t$, where
$X=W^{2,p}(Q)$. Clearly, $\bfH$ is hemi-continuous and
monotone. Choosing $u_\tau =\ol y_{\eps\tau}$ the weak equations
\eqref{momentum-eq-disc+} and \eqref{BC-reg-disc+} are interpreted as
$\bfH(\ol y_{\eps\tau})=b_\tau$ with $b_\tau$ defined via
\[
\langle b_\tau,z\rangle = - \! \int_\varOmega \! \big( \sigma_\text{vi}
(\nabla\underline y_{\eps\tau} , \nabla\DT y_{\eps\tau}, 
 \underline\theta_{\eps\tau}) {+} \eps \nabla\DT y_{\eps\tau} {+} 
 \sigma_\text{el} (\nabla\ol y_{\eps\tau},\underline\theta_{\eps\tau}) \big):
 \nabla z \d x \d t + \int_0^T \! \langle \ol\ell_\tau, z\rangle \d t .
\]
We obtain $b_\tau \weak b$ with $b$ defined by 
\[
\langle b,z\rangle = - \! \int_\varOmega \! \big( \sigma_\text{vi}
(\nabla y_{\eps } , \nabla\DT y_\eps , 
 \theta_\eps ) {+} \eps \nabla\DT y_\eps  {+} 
 \sigma_\text{el} (\nabla y_\eps ,\theta_{\eps }) \big):
 \nabla z \d x \d t + \int_0^T \! \langle \ell, z\rangle \d t ,
\]
because we can pass to the limit $\tau$ in all four terms
separately. For the first term we applying the lower semicontinuity
result \cite[Thm.\,7.5]{FonLeo07MMCV} twice, namely for the integrands
$f_\pm(x,(F,\theta),G)=\pm \sigma_\mathrm{vi}(F,G,\theta){:}\nabla
z(x)$ which both are convex in $G$. The limit passage in the second
term is simple weak convergence, and the fourth term converges because
of $\ol\ell_\tau \to \ell $ in $L^2\big(I;H^1_\mathrm{D}(\varOmega)^*\big)$. 
In the third term we exploit 
\[
\nabla \ol y_{\eps\tau} \in \mathbb F(K_\eps):= 
\bigset{F\in \R^{d \times d} } { |F|\leq K_\eps,\ \det F\geq 1/K_\eps}
\] 
(see \eqref{est-u}
and \eqref{est-det.nabla.y} from Proposition \ref{pr:FirstAprEstim}),
such that using \eqref{ass-psi} and \eqref{ass-phi} the map $(F,\theta)\mapsto
\sigma_\mathrm{el}(F,\theta)=\pl_F\varphi(F)+\pl_F \phi(F,\theta)$ is
continuous and bounded on $\mathbb F(K_\eps)\times \R^+$. Hence, with
\eqref{conv-strong} and Lebesgue's dominated convergence theorem we
obtain the desired convergence.  

To use Minty's trick \eqref{eq:Minty} we still need to check $\langle
b_\tau, \ol y_{\eps\tau}\rangle \to \langle b,y_\eps\rangle$. However,
as we have shown above $b_\tau$ is bounded (and hence weakly
converging to $b$) in
$L^2\big(I;H^1_\mathrm{D}(\varOmega)^*\big)$ and $ \ol
y_{\eps\tau} \to  y_\eps $ in $L^2\big(I;H^1_\mathrm{D}(\varOmega)\big)$
strongly (by \eqref{conv-y-strong}, the result follows immediately. 
Hence, we conclude $\bfH(y_\eps)=b$, which is nothing else than the
regularized momentum balance \eqref{momentum-eq-reg}, \eqref{BC1-reg},
and \eqref{BC1A-reg}. 


\STEP{Step 3: Balance of mechanical energy.} 
For the limit passage in the heat equation we need strong $L^2$-convergence of 
$\nabla\DT y_{\eps\tau}$ due to the viscous dissipation
$\xi^\mathrm{reg}_\eps(F,\DT F,\theta)$ that is nonlinear in $\DT F$. 
The strategy is to use the balance of mechanical energy as
follows. Rewriting the regularized momentum balance
\eqref{momentum-eq-reg}, \eqref{BC1-reg}, and \eqref{BC1A-reg}
in the form 
\[
\rmD_{\DT y}\calR(y_\eps,\DT y_\eps, \theta_\eps) +\eps\nabla \DT y_\eps + \rmD
\calM (y_\eps) + \rmD_y \Phi_\mathrm{cpl}(y_\eps,\theta_\eps) = \ell(t)
\]
with $\calM$ and $\Phi_\mathrm{cpl}$ defined in \eqref{eq:Energies}. 
We can now test with $\DT y_\eps \in L^2(I;H^1_\mathrm{D}(\varOmega))$ and
use (after decomposing $\calM=\calH +\Phi_\mathrm{el}$, see
\eqref{eq:Energies}) the chain rule in 
Proposition \ref{pr:ChainRule2} to obtain the balance of mechanical
energy  in the form 
\begin{align}
 &\nonumber 
   \calM(y_\eps(T))+ \int_0^T\!\! \big(2 \calR(y_\eps,\DT
  y_\eps, \theta_\eps) {+}\eps\|\nabla \DT y_\eps\|_{L^2}^2\big) \d t \\
&  = \calM(y_0) +\int_0^T  \langle \ell,\DT
  y_\eps \rangle \d t - \int_Q \pl_F\phi(\nabla y_\eps,\theta_\eps) {:} \nabla
  \DT y_\eps \d x \d t . 
  \label{eq:BalMechEnerg}
\end{align}
Indeed, by Proposition \ref{pr:Gat.Cvx} we know that $\calM$ satisfies
the assumptions of Proposition \ref{pr:ChainRule2} with space
$X=H^1_{\GDir}(\varOmega;\R^d)$. Clearly, $y_\eps \in H^1(I;X)$ and
$\calM(y_\eps(t))\leq \widetilde K_\eps$, see \eqref{eq:calM.bound}. 
Moreover, for 
\[
\varXi = \ell(t)-\rmD_{\DT y}\calR(y_\eps,\DT y_\eps,\theta_\eps)  - 
\eps\nabla \DT y_\eps  - \rmD_y \Phi_\mathrm{cpl}(y_\eps,\theta_\eps) 
\]
we have $\varXi(t)= \rmD\calM(y(t))$ a.e.\ in $[0,T]$ and our a priori
estimates provide $\varXi \in L^2([0,T];H^1_{\GDir}(\varOmega)^*)$. Thus,
\eqref{eq:BalMechEnerg} follows from Proposition \ref{pr:ChainRule2}.


\STEP{Step 4: Strong convergence of strain rates.}
The next step is now to derive a similar mechanical energy balance for the
time-discretized solutions, which is better than the previously
used estimate \eqref{test-y-discrete}. Passing to the limit $\tau\to
0$ from the latter estimate we would arrive at an estimate like
\eqref{eq:BalMechEnerg}, but with $2\calR$ and $\eps$ replaced by
$\calR$ and $\eps/2$, respectively. 

To improve the discrete bounds used in Proposition
\ref{pr:FirstAprEstim} we can exploit the a priori estimates
$\calM(y^k_{\eps\tau})\leq K_\eps$, which allow us to use the geodesic
$\Lambda$-convexity result in Proposition \ref{pr:Gat.Cvx}. Instead of
using the minimization property of $y^k_{\eps\tau}$ in
\eqref{minimize-y} we test the Euler-Lagrange equation
\eqref{momentum-eq-disc} with boundary conditions \eqref{BC1-reg-disc} 
and \eqref{BC2-reg-disc} by $y^k_{\eps\tau}{-}y^{k-1}_{\eps\tau}$ to
obtain 
\begin{align*}
& \tau 2\calR(y^{k-1}_{\eps\tau}, \DELTA y^k_\eps ,
\theta^{k-1}_{\eps\tau}) + \tau \eps \|\nabla \DELTA
y^k_\eps\|_{L^2}^2 +\rmD_y
\calM(y^k_{\eps\tau})[y^k_{\eps\tau}{-}y^{k-1}_{\eps\tau}] \\
&\hspace{10em} = \langle \ell^k_\tau,
 y^k_{\eps\tau}{-}y^{k-1}_{\eps\tau}\rangle - \rmD_y
 \Phi_\mathrm{cpl}(y^k_{\eps\tau},\theta^{k-1}_{\eps\tau}) 
  [y^k_{\eps\tau}{-}y^{k-1}_{\eps\tau}], 
\end{align*}
where we have the correct factors $2\calR$ and $\eps$. 
To recover the energy values $\calM(y^{j}_{\eps\tau})$ we now
eliminate the term involving $\rmD\calM$ using the
$\Lambda$-convexity estimate \eqref{eq:M.LambdaCvx} with
$y^{(1)}=y^k_{\eps\tau}$ and $y^{(2)}=y^{k-1}_{\eps\tau}$, which yields 
\begin{align*}
&\calM(y^k_{\eps\tau}) + \tau 2 \calR (y^{k-1}_{\eps\tau}, \DELTA y^k_\eps ,
 \theta^{k-1}_{\eps\tau}) + \big(\tau \eps - \tau^2\Lambda(K_\eps) 
 \big) \| \nabla  \DELTA y^k_\eps\|_{L^2}^2 \\
& \qquad \leq \calM(y^{k-1}_{\eps\tau}) + \tau \langle \ell_\tau^k, \DELTA
  y^k_\eps\rangle - \rmD_y \Phi_\mathrm{cpl} (
  y^k_{\eps\tau},\theta^{k-1}_{\eps\tau}) [\DELTA y^k_\eps]. 
\end{align*}
We now sum this inequality over $k=1,,\ldots, N_\tau$ and using the
interpolants we obtain the integral estimate 
\begin{align}
\nonumber
 &\calM(\ol y_{\eps\tau}(T)) + \int_0^T\!\!2\calR(
   \underline y_{\eps\tau}, \DT y_{\eps\tau},
   \underline\theta_{\eps\tau}) \d t + 
 (\eps{-}\tau \Lambda(K_\eps))\int_Q |\nabla\DT y_{\eps\tau}|^2 \d x \d t 
\\
\label{eq:ImprovedMEB}
& \qquad \leq \calM(y_0) + \int_0^T\!\!\bigg( \langle \ol\ell_\tau,
\DT y_{\eps\tau} \rangle - \int_\varOmega \pl_F\phi(\nabla \ol
y_{\eps\tau} ,\underline\theta_{\eps\tau}) \d x  \bigg) \d t .
\end{align}
Using the the convergences \eqref{conv-} and \eqref{conv-strong} it is
immediate to see that the all the terms on the right-hand side
converge to the corresponding terms on the right-hand side in
\eqref{eq:BalMechEnerg}. Now denote the three terms on the
left-hand side by $I^{(j)}_{\eps\tau}$ and set $I^{(j)}_\eps =
\liminf_{\tau\to 0^+} I^{(j)}_{\eps\tau}$. Using lower semicontinuity
arguments (use \cite[Thm.\,7.5]{FonLeo07MMCV} once again for
$I^{(2)}_{\eps\tau}$) we find 
\begin{align}
\nonumber 
&&\ol y_{\eps\tau}(T)\weak y_\eps(T) \text{ in } 
  W^{2,p}(\varOmega;\R^d)&&\Longrightarrow \qquad  & I^{(1)}_\eps  
  \geq \calM(y_\eps(T)) ,
\\ 
\nonumber
&&\nabla\DT y_{\eps\tau}\weak \nabla\DT y_\eps  \text{ in }
  L^2(Q;\R^{d\times d})   &&\Longrightarrow \qquad & I^{(2)}_\eps  \geq
  \int_0^T\!\! 2 \calR(y_\eps, \DT y_\eps,\theta)\d t, &&
\\
\label{eq:liminf.Ij}
&&\nabla\DT y_{\eps\tau}\weak \nabla\DT y_\eps  \text{ in }
  L^2(Q;\R^{d\times d})   &&\Longrightarrow \qquad  & I^{(3)}_\eps  
  \geq \eps \| \nabla \DT y_\eps\|_{L^2(Q)}^2.  
\end{align}
Thus, passing to the liminf on the left-hand side and to the limit on
the right-hand side in \eqref{eq:ImprovedMEB} and comparing with
\eqref{eq:BalMechEnerg} we obtain 
\begin{align*} 
  I^{(1)}_\eps   {+}  I^{(2)}_\eps  {+}  I^{(3)}_\eps  \leq \text{RHS} = 
   \calM(y_\eps(T))+ \int_0^T\!\! \big(2 \calR(y_\eps,\DT
  y_\eps, \theta_\eps) {+}\eps\|\nabla \DT y_\eps\|_{L^2}^2\big) \d t .
\end{align*}
Together with \eqref{eq:liminf.Ij} we conclude that we must have
equality in all three cases after ``$\Longrightarrow$''. However, 
$\nabla\DT y_{\eps\tau}\weak \nabla\DT y_\eps  $  in 
$ L^2(Q;\R^{d\times d}) $ and 
\[
I^{(3)}_\eps  = \liminf_{\tau\to 0} (\eps{-}\tau \Lambda(K_\eps)) \|
\nabla \DT  y_{\eps\tau}\|_{L^2(Q)}^2  
  = \eps \| \nabla \DT y_\eps\|_{L^2(Q)}^2
\]
imply the desired strong convergence $\nabla\DT y_{\eps\tau} \to 
\nabla\DT y_\eps  $ in $ L^2(Q;\R^{d\times d})$. 


\STEP{Step 5: Limit in the heat equation.} We first pass to the
limit $\tau\to 0$ in the constitutive relation
\eqref{heat-equation-disc+}, namely $\ol
w_{\eps\tau}=\mathfrakw(\nabla\ol y_{\eps\tau},
\ol\theta_{\eps\tau})$.  The left-hand side converges to $w_\eps$ by  
\eqref{eq:w.convg}, while the right-hand side converges to $\mathfrakw
(\nabla y_\eps,\theta_\eps)$ by the continuity of $\mathfrakw $, the
bound \eqref{eq:mfw.estim} and the
convergences  \eqref{conv-strong}. Thus, $w_\eps=\mathfrakw
(\nabla y_\eps,\theta_\eps)$ is established, i.e.\
\eqref{eq:syst.reg.w} holds. 
 
We write the heat equation \eqref{heat-equation-disc+} with boundary
conditions \eqref{BC2-reg-disc+} in the weak form 
\begin{align}
\nonumber
&\int_Q\Big(\DT w_{\eps\tau} z +\nabla
  \ol\theta_{\eps\tau}\cdot\mathcal K(\nabla \underline y_{\eps\tau},
    \underline\theta_{\eps\tau} ) \nabla z \Big) \d x \d t 
  + \int_\varSigma  \kappa
      \big(\ol\theta_{\eps\tau}{-}\ol\theta_{\flat,\eps,\tau}\big) z
      \d S \d t 
\\
\label{eq:heat.epstau}
& = \int_Q\big(\xi^\mathrm{reg}_\eps(\nabla\underline y_{\eps\tau}, 
      \!\nabla\DT y_{\eps\tau},\underline\theta_{\eps\tau})
    {+}\pl_F^{}\phi(\nabla\overline y_{\eps\tau}, 
      \overline\theta_{\eps\tau}){:}\nabla\DT y_{\eps\tau}\big) z\d
      x \d t   
\end{align}
for all $z \in L^\infty(I;H^1(\varOmega))$. While we only have the weak
convergences $\DT w_{\eps\tau} \weak \DT w_\eps$ in $L^2 \big( I; H^1
(\varOmega)^*\big)$ (see \eqref{eq:w.convg}) and $\nabla
\ol\theta_{\eps\tau} \weak \nabla \theta_\eps$ in $L^2(Q)$ (see
\eqref{conv-theta-}), all other functions in \eqref{eq:heat.epstau}
converge strongly.  In particular, using the strong convergences
$\nabla\DT y_{\eps\tau} \to \nabla\DT y_\eps $ in $ L^2(Q;\R^{d\times
  d})$ and $0\leq \xi^\mathrm{reg}_\eps(\nabla\underline y_{\eps\tau},
\!\nabla\DT y_{\eps\tau},\underline\theta_{\eps\tau}) \leq K_\eps$ we
obtain
\begin{align}
&\xi^\mathrm{reg}_\eps(\nabla\underline y_{\eps\tau}, 
      \!\nabla\DT y_{\eps\tau},\underline\theta_{\eps\tau})
 \to    \xi^\mathrm{reg}_\eps(\nabla y_{\eps },\nabla\DT y_\eps, 
         \theta_\eps) \quad \text{strongly in } L^p(Q) \text{ for all
         }p\in {]1,\infty[}.
\label{conv-dissip-heat}
\end{align}
Thus, passing to the limit $\tau\to 0$ in \eqref{eq:heat.epstau} leads
exactly to the weak form to the regularized heat equation 
 \eqref{heat-equation-reg} with boundary condition \eqref{BC2-reg}.

This conclude the proof of Proposition \ref{pr:tau-0}. 
\end{proof}

\section{Limit passage $\eps \to 0$}
\label{se:Limit}
 
In this final step of the proof of Theorem \ref{thm:MainExist} we
have to pass to the limit with the regularization parameter $\eps\to
0$. As we are already in the time-continuous setting we are now able
to make the formally derived total energy balance \eqref{total-engr}
for $\calE$ rigorous for all $\eps>0$. From this we will be able to
derive a priori bounds for $(y_\eps, \theta_\eps)$ that are
independent of $\eps$. 

\begin{remark}[Missing discrete estimate for the total energy] 
\label{rm:DiscreteEnergyBal} The derivation of the total energy
balance is achieved by testing the momentum balance by
$\DT y$ and the heat equation by the constant function $1$. The
corresponding step on the time-discrete level would be the test  
\eqref{momentum-eq-disc} by $\DELTA y^k$ and
\eqref{heat-equation-disc} by $1$. We would be able to use the desirable
cancellation of the dissipation, namely  $\xi^\text{reg}_\eps -\xi
\leq 0$; however for the coupling terms 
\[
\pl_F \phi(\nabla y^k_{\eps\tau},\theta^{k-1}_{\eps\tau}) : \DELTA
\nabla y^k_\eps  \quad \text{ and } \quad \pl_F \phi(\nabla 
y^k_{\eps\tau},\theta^{k}_{\eps\tau}) : \DELTA \nabla y^k_\eps ,
\]
which arise from \eqref{momentum-eq-disc} and
\eqref{heat-equation-disc} respectively, we do not have any way to
estimate the first against the second. Recall that we were forced to
use the explicit/forward value $ \theta^{k}_{\eps\tau}$ to maintain
positivity of the temperature.
\end{remark}

To exploit the balance of the total energy we have to strengthen the
assumption on the leading $\ell(t)$, i.e.\ the functions $g$, and $f$,
in \eqref{ass-g}, namely 
\begin{equation}
  \label{eq:g.f.strong}
  g\in W^{1,1}(I; L^2(\varOmega;\R^d)), \qquad 
  f\in W^{1,1}(I; L^2(\GNeu;\R^d)).
\end{equation}
This implies that $t \mapsto \ell(t)$ lies in
$W^{1,1}\big(I;H^1_{\GDir}(\varOmega;\R^d)^*\big)$, which is what we
will only need. 

The new $\eps$-independent estimates on $\nabla \DT y_\eps$ in
$L^2(Q)$ will be obtain by exploiting the Pompe's generalized Korn's
inequality (cf.\ \cite{Pomp03KFIV}) as prepared in Theorem
\ref{th:Korn} above.

\begin{lemma}[A-priori estimates for $y_\eps $]
 \label{lem2}
  Let the assumptions \eqref{ass} and \eqref{eq:g.f.strong} 
  hold. Then there exists a constant $K$ such that for all $\eps \in
  {]0,1[}$ and all weak
  solutions $(y_\eps ,\theta_\eps )$ of the regularized problem
  \eqref{system-reg}-\eqref{BC-reg} obtained in
  Proposition~\ref{pr:tau-0} we have the a priori estimates 

 Then $\det(\nabla y_{\eps })>0$ on $Q$
  and the following estimates hold with $K$ independent of $\eps >0$:
\begin{subequations}\label{est+}
\begin{align}\label{est-u+}
 &\big\|y_{\eps }\big\|_{L^\infty(I;W^{2,p}(\varOmega;\R^d))}\le K,
\\
 \label{est-det.nabla.y+}
 & \det \big( \nabla y_\eps(t,x)\big) \geq 1/K \ 
     \text{ for all } (t,x)\in Q,
\\
 \label{est-theta+}
 &\big\|\theta_{\eps }\big\|_{L^\infty(I;L^1(\varOmega))}\le K,
\\
 \label{est-F-dot+} &\big\|\nabla\DT y_\eps
 \big\|_{L^2(Q;\R^{d\times d})}\le K,
\\ 
 \label{est-xi+} 
 &\int_Q \xi(\nabla y_\eps,\nabla\DT y_\eps, \theta_\eps) \d x \d t \le K,
\end{align}
\end{subequations}
with $q$ from \eqref{ass-psi}, where again $\sym(\cdot)$ denotes
the symmetric part of a $(d{\times}d)$-matrix. 
\end{lemma}
\begin{proof} We proceed in
  two steps that are 
  close to estimates we have done in the time-discrete setting.


  \STEP{Step 1: Estimate for $\calE(y_\eps,\theta_\eps)$.} 
Using the derived regularity for the solution $(y_\eps,\theta_\eps)$
we see that a suitable variant of the total energy balance
\eqref{total-engr} holds. To be specific, we start from
\eqref{eq:BalMechEnerg}, which is also valid for arbitrary $t\in
{]0,T]} $ in place of $T$,  and add the time-integrated version of
\eqref{heat-equation-reg} tested with the constant function $z\equiv
1$. Using $\calE=\calM+\calW$ with
$\calW(y_\eps,\theta_\eps)=\int_\varOmega w_\eps \d x $ we find 
\begin{align*}
 & \calE(y_\eps(t),\theta_\eps(t)) +  \int_0^t\int_\varOmega \Big( 
 \xi(\nabla y_\eps, \nabla \DT y_\eps,\theta_\eps) + \eps |\nabla \DT
 y_\eps|^2 - \xi^\text{reg}_\eps(\nabla y_\eps, \nabla \DT
 y_\eps,\theta_\eps) \Big) \d x  \d s  \\
& = \calE(y_\eps(0),\theta_\eps(0)) + \int_0^t \langle
 \ell(s),\DT y_\eps(s) \rangle \d s + \int_0^t \int_\varGamma
 \kappa(\theta^\eps_\flat{-}\theta_\eps) \d S \d s .
\end{align*} 
The importance is the cancellation of the term $\pl_F \phi : \nabla
\DT y_\eps$ and that the difference of the dissipation integrals has a
sign. 

Defining the auxiliary variable $E_\eps(t):=
\calE(y_\eps(t),\theta_\eps(t)) - \langle \ell(t),y_\eps(t)\rangle $
and using $0\leq \theta^\eps_\flat \leq \theta_\flat$ and $\theta_\eps\geq
0$ gives 
\[
E_\eps(t) \leq E_\eps(0) + \int_0^t \Big( \int_\varGamma \kappa
\theta_\flat \d S - \langle \dot\ell(s),y_\eps(s)\rangle \Big) \d s,
\]
where we have integrated by parts the power of the external loadings,
which was possible by the strengthened assumption
\eqref{eq:g.f.strong}. 

With $\calE\geq \calM\geq \calH$ and the coercivity of $\calH$ we
have $\|y\|_{H^1} \leq c_1 + c_2 \calE(y,\theta)$ and obtain 
\[
E_\eps(t) \leq E_\eps(0) + \int_0^t \big( a(s) + b(s)E_\eps(s)\big) \d
s \ \text{ with }a,\,b \geq 0
\]
and $a,b\in L^1(0,T)$, which follows from  \eqref{eq:g.f.strong} for
$\ell$ and \eqref{ass-BC} for $\theta_\flat$. With $B(t)=\int_0^t b(s)
\d s $ and $A(t) = \int_0^t a(s) \d s$ the Gronwall estimate yields the
a priori estimate
\[
\calE_\eps(t) \leq \mathrm e^{B(t)}\big( E_\eps(0)+ A(t)\big) \leq
\mathrm e^{B(T)} \big( E^0 + A(T)\big) :=M_1,
\]
where we used $\calE_\eps(0)=\calE(y_\eps(0),\theta_\eps(0)) \leq 
\calE(y_0,\theta_0)-\langle \ell(0), y_0\rangle =:E^0<\infty$ by
\eqref{ass-IC}, \eqref{ass-BC}, and \eqref{eq:mfw.estim}. 
This immediately implies 
\[
\calM(y_\eps(t)) + \Epsilon \| \theta_\eps(t)\|_{L^1(\varOmega)}  \leq
\calE_\eps(y_\eps(t),\theta_\eps(t)) \leq M_2 .
\]
Hence, \eqref{est-theta+} is established, whereas \eqref{est-u+} and
\eqref{est-det.nabla.y+} follow by applying Theorem
\ref{th:HealeyKromer}. 


\STEP{Step 2: Estimate for the strain rate $\nabla \DT y_\eps$.}
We return to the mechanical energy balance \eqref{eq:BalMechEnerg} on
the interval $I=[0,T]$. We recall that the dissipation function
$\xi(F,\DT F,\theta)$ is assumed to control the symmetric part of
$F^\top \DT F$ only, namely 
\[
\xi(F,\dot F,\theta) = 2\hat\zeta(F^\top
F, F^T\DT F{+}\DT F^\top F,\theta) \geq \alpha|F^T\DT F{+}\DT F^\top
F|^2.
\]
Using our a priori bounds on $\calM(y_\eps(t))$, we can apply the
generalized Korn's inequality a prepared in Corollary
\ref{cor:GenKorn} with $y=y_\eps(t,\cdot)$ and  $v= \DT y_\eps(t) \in
H^1_{\GDir}(\varOmega; \R^d)$ to
obtain
\begin{align*}
 & \alpha \ol c_K \int_0^T\| y_\eps(t)\|_{H^1}^2 \d t \leq \int_Q
 \alpha \big| 
   \nabla y_\eps^\top \nabla \DT y_\eps {+} \nabla \DT y_\eps^\top
   \nabla y_\eps  \big|^2 \d x \d t   
\leq \int_Q \xi(\nabla y_\eps ,\nabla \DT
y_\eps, \theta_\eps) \d x \d t  
\\
&\leq \calM(y_0)-\calM(y_\eps(T)) +
\int_0^T \big( \|\ell(t)\|_{(H^1)^*} + \|\pl_F\phi(\nabla
y_\eps,\theta_\eps)\|_{L^\infty(Q)}\big) \| \DT y_\eps(t) \|_{H^1} \d t ,
\end{align*}
where we used $|\pl_F \phi(F,\theta)|C(1{+}|F|)^s$ and $|\nabla
y_\eps(t,x)|\leq K$, which follows from \eqref{est-u+}. From this,
\eqref{est-F-dot+} and \eqref{est-xi+}  follow immediately. 
\end{proof}

For the deformation $y_\eps$ we have all the estimates we need
for passing to the limit. But we still need good a priori estimates
for the temperature. Here the problem arises that the heating arising
through the viscous dissipation $\xi(\nabla y_\eps, \nabla \DT
y_\eps,\theta_\eps)$ is only bounded in $L^1(Q)$. So, obtaining
improved estimates we have to invoke special test functions developed
by Boccardo and Gallou\"et \cite{BocGal89NEPE} for parabolic equations
with measure-valued right-hand sides. 

\begin{proposition}[A priori estimates for $\theta_\eps$ and $w_\eps$]
\label{pr:Add.Est.eps}
Under the conditions of Lemma~\ref{lem2}, also the following estimates hold:
\begin{subequations}
\label{est++}
\begin{align}
\label{eq:est.w.theta}
&\forall\, p\in \big[1,\text{\large$\tfrac{d{+}2}d$}
   \big[\ \exists\, C_p>0 \ \forall\,
\eps \in {]0,1]}: \quad 
\|\theta_\eps\|_{L^p(Q)} + \|w_\eps\|_{L^p(Q)} \leq C_p , 
\\
\label{eq:nabla.w.theta}
&\forall\, r\in \big[1,\text{\large$\tfrac{d{+}2}{d{+}1}$} 
 \big[\ \exists\, K_r>0 \ \forall\,
\eps \in {]0,1]}: \quad 
\|\nabla\theta_\eps\|_{L^p(Q)} + \|\nabla w_\eps\|_{L^p(Q)} \leq K_r ,
\\
\label{est-dot-w}
&\exists \,K>0 \ \forall\, \eps \in {]0,1[}: \quad 
  \big\|\DT w_{\eps }\big\|_{L^1(I;H^{(d+3)/2}(\varOmega)^*)}\le K.
\end{align}
\end{subequations}
\end{proposition}
\begin{proof}  We follow the recipe in \cite{BocGal89NEPE} in the
  simplified variant  of \cite{FeiMal06NSET},  see also
  \cite{MieNau18?EWSK}. For $\eta\in {]0,1[}$ we define the function
  $\chi_\eta :\R^+\to \R^+$ via 
\begin{align*}
\chi_\eta(0)=0 \quad \text{and} \quad 
\chi'_\eta(w):=1-\frac1{(1{+}w)^\eta} \in[0,1]. 
\end{align*} 
Clearly, $\chi_\eta$ satisfies $\min\{0,w/2{-}C_\eta \} \leq
\chi_\eta(w) \leq w$ and $\chi''_\eta(w)=
\dfrac\eta{(1{+}w)^{1+\eta}}>0$. 

Now testing \eqref{heat-equation-reg} with the test function $z=
\chi'_\eta\circ w_\eps$ amounts to applying the chain rule in
Proposition \ref{pr:ChainRule1} to the convex functional
$\calJ(w)=\int_\varOmega \chi_\eta(w(x)) \d x$ on the space
$X=H^1(\varOmega)^*$. Indeed, from \eqref{conv-} and $w_\eps = \mathfrakw
(\nabla y_\eps,\theta_\eps)$ we have $w_\eps \in L^2 (I; H^1(\varOmega))
\cap H^1(I;H^1(\varOmega)^*)) $, and the chain rule gives the first
identity in the following calculation: 
\begin{align*}
&\frac\rmd{\rmd t} \int_\varOmega \chi_\eta(w_\eps)\d x = \int_\varOmega
  \chi'_\eta(w_\eps) \DT w_\eps \d x \\
&= - \int_\varOmega \chi''_\eta(w_\eps)
\nabla w_\eps\cdot \mathcal K(\nabla y_\eps,\theta_\eps)\nabla \theta_\eps
\d x   + \int_\varGamma
\kappa (\theta_\flat^\eps{-}\theta_\eps)\d S \\
&\quad + \int_\varOmega \chi'_\eta(w_\eps) \Big( \xi^\text{reg}_\eps(\nabla y_\eps,
\nabla \DT y_\eps, \theta_\eps) + \pl_F \phi(\nabla
y_\eps,\theta_\eps) {:} \nabla \DT y_\eps\Big) \d x .  
\end{align*}
Integration over $t\in I=[0,T]$ and using $\chi'_\eta(w) \in [0,1]$ and
$\|\nabla y_\eps\|_{L^\infty(Q)} \leq K_\infty$ yield
\begin{align}
\nonumber
&\int_Q \chi''_\eta(w_\eps)
\nabla w_\eps\cdot \mathcal K(\nabla y_\eps,\theta_\eps)\nabla \theta_\eps
\d x \d t \\
\label{eq:nabla.weps}
&\leq \int_\varOmega \chi_\eta(w_0)\d x + \int_\varSigma \kappa
\theta_\flat \d S \d t + \int_Q \!\Big( \xi(\cdots)+ C(1{+}K_\infty)^s|\nabla
\DT y_\eps|\Big) \d x \d t \leq C,
\end{align}
where we used \eqref{ass-IC}, \eqref{ass-BC}, \eqref{est-F-dot+}, and
\eqref{est-xi+}. 

From this, we derive an a priori bound on $\nabla w_\eps$ by setting
$\widetilde{\mathcal K}_\eps =\mathcal K(\nabla y_\eps,\theta_\eps)$
and estimate it as in \eqref{eq:calK} (see Step~3 of the proof of Proposition
\ref{pr:FirstAprEstim}) by 
\[
|\widetilde{\mathcal K}_\eps(t,x)| \leq \varkappa  \quad \text{and} \quad 
a\cdot \widetilde{\mathcal K}_\eps(t,x) a  \geq \frac1\varkappa |a|^2 ,
\]
where $\varkappa$ is now independent of $\eps$ because of the
$\eps$-independent bound in \eqref{est-u+} and
\eqref{est-det.nabla.y+}. Moreover, $\nabla w_\eps$ and $\nabla
\theta_\eps$ are related by
\begin{equation}
  \label{eq:na.w.na.theta}
  \nabla w_\eps = \pl_\theta \mathfrakw(\nabla y_\eps,\theta_\eps)
\nabla \theta_\eps + \pl_F \mathfrakw(\nabla y_\eps,\theta_\eps):
\nabla^2 y_\eps.
\end{equation}
With $\pl_\theta w(F,\theta)=-\theta \pl_\theta^2 \phi(F,\theta) \leq
K$ we obtain
\begin{align*}
\frac1{\varkappa}|\nabla w_\eps|^2 &\leq  
  \nabla w_\eps \cdot \widetilde{\mathcal K}_\eps \nabla w_\eps 
\\
& =  \pl_\theta\mathfrakw(\nabla y_\eps,\theta_\eps)\:  \nabla
     w_\eps\cdot \widetilde{\mathcal K}_\eps \nabla \theta_\eps 
  + \nabla w_\eps\cdot \widetilde{\mathcal K}_\eps \pl_F \mathfrakw 
   (\nabla y_\eps,\theta_\eps): \nabla^2 y_\eps\\
&\leq K 
  \nabla w_\eps \cdot \widetilde{\mathcal K}_\eps \nabla \theta_\eps
  +  \varkappa |\nabla w_\eps|\, C(1{+}K_\infty)^s
  |\nabla^2 y_\eps| \\
& \leq K 
  \nabla w_\eps \cdot \widetilde{\mathcal K}_\eps \nabla \theta_\eps
  + \frac1{2\varkappa} |\nabla w_\eps|^2
   + C_* |\nabla^2 y_\eps|^2.
\end{align*}  
Canceling $\frac1{2\varkappa} |\nabla w_\eps|^2$, multiplying by
$\chi''(w_\eps)\in [0,1]$, and integrating over $Q$ we employ
\eqref{eq:nabla.weps} and arrive at 
\[
\frac1{\varkappa K} \int_Q  \chi''_\eta(w_\eps)|\nabla w_\eps|^2 \d x \d t
\leq \int_Q \chi''(w_\eps)\Big(\nabla w_\eps \cdot \widetilde{\mathcal
  K}_\eps \nabla \theta_\eps + C_* |\nabla^2 y_\eps|^2\Big) \d x \d t
\leq C_3,
\]
where the last integrand is bounded by \eqref{est-u+} and $p\geq 2$.  

For $r\in {[1,2[}$ we set $p=2/(2{-}r)$, $p'=2/r$, and $q=(1{+}\eta)
r/2$ and employ H\"older's estimate to obtain 
\begin{align}
\nonumber
& \| \nabla w_\eps\|_{L^r(Q)}^r =  \int_Q
  (1{+}w_\eps)^q\:\frac{|\nabla w_\eps|^r}{(1{+}w_\eps)^q} \d x \d t 
  \; \leq \; \|(1{+}w_\eps)^q\|_{L^p(Q)} \Big\|\frac{|\nabla
  w_\eps|^r}{(1{+}w_\eps)^q}  \Big\|_{L^{p'}(Q)}  
\\
\label{eq:n.weps:Holder}
& = \|1{+}w_\eps\|^q_{L^{qp}(Q)} \Big(\int_Q \frac{|\nabla
  w_\eps|^2}{(1{+}w_\eps)^{1+\eta}} \d x \d t\Big)^{1/p'}   \leq  \; 
\|1{+}w_\eps\|^q_{L^{qp}(Q)} \big(\varkappa K C_3/\eta\big)^{1/p'},
\end{align}
where crucially relied on  $p'=2/r$, $\chi''(w)=\eta/(1{+}w)^{1+\eta}$,
and the previous estimate. 
Using the a priori estimate $\| 1{+}w_\eps\|_{L^\infty(I;L^1(\varOmega))}
\leq T|\varOmega| + K=:K_1$ from \eqref{est-theta+}  we can now use the
anisotropic Gagliardo-Nirenberg interpolation (see e.g.\
\cite[Lem.\,4.2]{MieNau18?EWSK}) giving 
\begin{align*}
&\|1{+}w_\eps\|_{L^{r/\lambda}(Q)} \leq C \|
1{+}w_\eps\|_{L^\infty(I;L^1(\varOmega))}^{1-\lambda} \big(\|
1{+}w_\eps\|_{L^\infty(I;L^1(\varOmega))} + \|\nabla
w_\eps\|_{L^r(Q)}\big)^\lambda \ \text{ with }  \lambda =\frac d{d{+}1}. 
\end{align*}
For inserting this into \eqref{eq:n.weps:Holder} we need $qp \leq
r/\lambda$ which gives the restriction $r\leq 2-(1{+}\eta) \lambda$. 

Thus, for all $r\in {[1,(d{+}2)/(d{+}1)[}$ we find an $\eta=\eta_r \in
{]0,1[}$ such that the above estimates give 
\[
\| \nabla w_\eps\|_{L^r(Q)}^r  \leq C_r\,\big(1+ \| \nabla
w_\eps\|_{L^r(Q)}^{q\lambda}\big),  
\]
and $q_r\lambda < q_r=(1{+}\eta_r)r/2 < r$ provide $\| \nabla
w_\eps\|_{L^r(Q)} \leq K_r$. Using \eqref{eq:na.w.na.theta} and
$\pl_\theta \mathfrakw \geq \Epsilon > 0$ we easily find $\| \nabla
\theta_\eps\|_{L^r(Q)} \leq K_r$ and \eqref{eq:nabla.w.theta} is
established. 

Applying the Gagliardo-Nirenberg interpolation once again 
gives assertion \eqref{eq:est.w.theta}. 
 
Eventually, the a priori estimate \eqref{est-dot-w} is obtained
estimating all other terms in \eqref{heat-equation-reg}, 
when realizing that always
$H^{(d+3)/2}(\varOmega)\subset W^{1,\infty}(\varOmega)$.
\end{proof}

We are now in the position to pass to the limit $\eps \to 0$ in the
regularized system  \eqref{system-reg}-\eqref{BC-reg}, and thus
provide the proof of our main existence result presented in Theorem
\ref{thm:MainExist}. The approach is close to the convergence result
presented in Proposition \ref{pr:tau-0}: first we extract converging
subsequences and then pass to the limit in the mechanical momentum
balance. This also provides the necessary strong convergence of the
the strain rates that is needed to eventually pass to the limit in the
heat heat equation. 

\begin{proposition}[Convergence for $\eps \to0$]
 \label{prop2} 
 Let again \eqref{ass} and \eqref{eq:g.f.strong} hold. 
 Then, considering the sequence of time
 steps $\eps \to0$, there is a subsequence $(y_\eps ,\theta_\eps )$ of
 weak solutions to the regularized system
 \eqref{system-reg}-\eqref{BC-reg} obtained in
 Proposition~\ref{pr:tau-0} such that, for some $(y,\theta)$, it holds
\begin{subequations}\label{conv}\begin{align}\label{conv-y}
&y_\eps \to y&&\text{weakly* in }\ L^\infty(I;W^{2,p}(\varOmega;\R^d))\cap H^1(I;L^2(\varOmega;\R^d))\ \ \text{ and }&&
\\&\label{conv-theta}
\theta_\eps \to\theta&&\text{weakly in }\ L^r(I;W^{1,r}(\varOmega))
\ \ \text{ for all }\ 1\le r<(d{+}2)/(d{+}1).&&
\end{align}\end{subequations}
Moreover, every couple $(y,\theta)$ obtained in such a way is a weak
solution, according Definition~\ref{def}, of the boundary-value
problem \eqref{system}--\eqref{BC} satisfying the initial values
\eqref{IC}.%
\end{proposition}
\begin{proof} The proof follows the lines of the proof of
Proposition  \ref{pr:tau-0}, so we do not repeat all details of the
arguments. 


\STEP{Step 1: Extraction of converging subsequences.}  Using the
a priori estimates \eqref{est+} and \eqref{est++}, Banach's selection
principle allows us to choose a subsequence and some $(y,\theta)$ such
that \eqref{conv} holds. By the Aubin-Lions' theorem interpolated with
the estimates \eqref{est-u} and \eqref{est-theta}, we have also
\begin{subequations}
  \label{conv-strong+}
 \begin{align}\label{conv-y-strong++}
   \nabla y_\eps &\to\nabla y&&\text{strongly in }\
   L^\infty(Q;\R^{d\times d})\ \ \text{ and } 
\\
 \label{conv-w-strong++} 
  w_\eps &\to w&&\text{strongly in }\ L^p(Q)
    \ \ \text{ with any }\ 1\le p<1+2/d,\\
 \label{conv-theta-strong++} 
  \theta_\eps &\to\theta&&\text{strongly in }\ L^p(Q)
    \ \ \text{ with any }\ 1\le p<1+2/d.
\end{align}
\end{subequations}
The proof of \eqref{conv-y-strong++} is similar to
\eqref{conv-y-strong}. For \eqref{conv-w-strong++} we proceed as for
\eqref{conv-theta-strong} by using the estimates on $w_\eps$ given in
\eqref{est++}. Using the relation $w_\eps = \mathfrakw(\nabla
y_\eps,\theta_\eps)$ we also obtain the strong convergence
\eqref{conv-theta-strong++}.
%


\STEP{Step 2: Convergence in the mechanical equation.} 
The limit passage in the momentum balance
\eqref{momentum-eq-reg}-\eqref{BC-reg} works as before, again using
the Minty trick \eqref{eq:Minty}. Of course, the additional
regularizing viscosity  term $\eps \nabla \DT y_\eps$ vanishes because
of our a priori bound \eqref{est-F-dot+}:
\[
\bigg|\int_Q\eps \nabla\DT y_\eps {:}\nabla z \,\d x\d t\bigg|
\le\eps \big\|\nabla\DT y_\eps \big\|_{L^2(Q;\R^{d\times
    d})}\big\|\nabla z\big\|_{L^2(Q;\R^{d\times
    d})}=   C \eps \to 0. 
\]  


\STEP{Step 3: Balance of mechanical energy.} As in the proof of
Proposition \ref{pr:tau-0} we derive from the property that the limit
couple $(y,\theta)$ solves the mechanical equation that the following
mechanical energy relation holds: 
\begin{align}
 &   \calM(y(T))+ \int_0^T\!\! 2 \calR(y,\DT
  y, \theta)  \d t   = \calM(y_0) +\int_0^T  \langle \ell,\DT
  y \rangle \d t - \int_Q \pl_F\phi(\nabla y,\theta) {:} \nabla
  \DT y \d x \d t . 
  \label{eq:BalMechEnerg++}
\end{align}


\STEP{Step 4: Strong convergence of the symmetric strain rates.}
We can pass to 
the limit $\eps\to 0$ in the mechanical energy relation
\eqref{eq:BalMechEnerg}. Comparing the result with
\eqref{eq:BalMechEnerg++} we obtain 
\begin{equation}
  \label{eq:Diss.Cvg}
  \lim_{\eps\to 0} \int_Q
  \xi(\nabla y_\eps,\nabla \DT y_\eps,\theta_\eps) \d x \d t = \int_Q
  \xi(\nabla y,\nabla \DT y,\theta) \d x \d t. 
\end{equation}
To conclude strong convergence we use the special form
\eqref{def-of-xi}, namely $\xi(F,\DT F,\theta)= 2\hat\zeta(F^\top F,
F^\top \DT F {+}\DT F^\top F,\theta)$. From the pointwise convergence
$\theta_\eps \to \theta$, the uniform convergence $F_\eps :=\nabla
y_\eps \to F=\nabla y$, and the weak convergence $\DT F_\eps :=\nabla
\DT y_\eps \weak  \DT F$ in $L^2(Q)$ we obtain 
\[
V_\eps:= F_\eps^\top \DT F_\eps {+} \DT F_\eps^\top F_\eps \weak 
F^\top \DT F{+} \DT F^\top F =:V \ \text{ in } L^2(Q).
\]
With the coercive and quadratic structure of $\hat \zeta$ assumed in
\eqref{ass-hat-zeta} we proceed as follows:
\begin{align*}
2\alpha\|V_\eps{-}V\|_{L^2(Q)}^2  &\leq  \int_Q 2\hat\zeta(
C_\eps,V_\eps{-}V,\theta_\eps) \d x \d t \\
& = \int_Q\Big( 2\hat\zeta( C_\eps,V_\eps,\theta_\eps) 
  - 2V_\eps : \mathbb D(C_\eps,\theta_\eps) V+
  2\hat\zeta( C_\eps,V,\theta_\eps) \Big)  \d x \d t
\\
&= \int_Q\Big(  \xi(F_\eps,\DT F_\eps,\theta_\eps) 
  - 2V_\eps : \mathbb D(C_\eps,\theta_\eps) V+
   \xi(F,\DT F,\theta) \Big)  \d x \d t + \delta(\eps),
\\
& \hspace{7em} \text{with } \delta(\eps)= \int_\varOmega 2 V:\big( 
 \mathbb D (C_\eps,\theta_\eps) {-} \mathbb D(C,\theta)\big) V  \d x \d t.  
\end{align*}
We see that the first term converges by \eqref{eq:Diss.Cvg}, while the
second term converges by the weak convergence $V_\eps \weak V$ and the
strong convergence $ \mathbb D(C_\eps,\theta_\eps) V \to  \mathbb
D(C,\theta) V$ (as $\mathbb D$ is bounded and the arguments converge
pointwise). Similarly, $\delta(\eps)\to 0$ by Lebesgue's dominated
convergence theorem, and thus we conclude the strong convergence 
$\|V_\eps{-}V\|_{L^2(Q)} \to 0$.  


\STEP{Step 5: Limit passage in the heat equation.} Testing the
regularized heat equation \eqref{heat-equation-reg} with boundary
conditions \eqref{BC2-reg} by smooth function $v$ with
$V(T,\cdot)\equiv 0$ we find 
\begin{align}
\nonumber
&\int_Q \Big(\nabla \theta_\eps \cdot \mathcal K(\nabla y_\eps, \theta_\eps)
   \nabla v  -\big(\xi^\mathrm{reg}_\eps(\nabla y_\eps,\nabla \DT
   y_\eps, \theta_\eps) {+}\pl_F \phi(\nabla y_\eps, \theta_\eps){:} 
   \nabla \DT y_\eps \big)v -w_\eps \DT v \Big) \d x \d t \quad \\
& \label{eq:HeatEq.eps.weak} \hspace{8em}
  + \int_\varSigma \kappa \theta_\eps v \d S \d t = 
  \int_\varSigma \kappa \theta_{\flat,\eps} v \d S \d t + \int_\varOmega
  \mathfrakw(\nabla y_0,\theta_{0,\eps}) v(0) \d x. 
\end{align} 
Here the first term passes to the limit by $\nabla\theta_\eps \weak
\nabla \theta$ and $\mathcal K(\nabla y_\eps, \theta_\eps)
   \nabla v \to \mathcal K(\nabla y, \theta)$. In the second term we
   use 
\[
\xi^\mathrm{reg}_\eps (\nabla y_\eps,\nabla \DT y_\eps,\theta_\eps)=
\frac{ \xi(\nabla y_\eps,\nabla \DT y_\eps,\theta_\eps)}
{1{+}\eps \xi(\nabla y_\eps,\nabla \DT y_\eps,\theta_\eps)} = 
  \frac{2\hat\zeta(C_\eps,V_\eps,\theta_\eps)} 
   {1 + 2\eps \hat\zeta (C_\eps,V_\eps,\theta_\eps)} \leq
   2K|V_\eps|^2=:g_\eps. 
\]
Because of Step~4, we know $V_\eps \to V$ strongly in $L^2( Q; \R^{d\times
  d}_\mathrm{sym} )$. Hence, we have $g_\eps \to g:=K|V|^2$ in
$L^1(Q)$ and may assume, after extracting another subsequence, 
$ V_\eps(t,x)\to V(t,x)$ a.e.\ in $Q$. By the uniform/pointwise convergence of
$C_\eps$ and $\theta_\eps$ for any $v\in C^0(\ol Q)$ we obtain 
\[
g_\eps\|v\|_{L^\infty(Q)} \geq \xi^\mathrm{reg}_\eps (\nabla y_\eps,\nabla \DT
y_\eps,\theta_\eps)v \to \xi(\nabla y,\nabla \DT y,\theta) v \ \ 
\text{ a.e.\ in } Q. 
\]
As the majorants $g_\eps\|v\|_{L^\infty(Q)}$ converge to
$g\|v\|_{L^\infty(Q)}$ in $L^1(Q)$ the generalized dominated
convergence theorem implies convergence of the second term in
\eqref{eq:HeatEq.eps.weak}. 

In the third term we have weak convergence of $\nabla \DT y_\eps$ and
strong convergence of $v\pl_F \phi(\nabla y_\eps,
\theta_\eps)$. Similarly, the remaining four terms converge to the
desired limits. Thus, we have shown that $(y,\theta)$ satisfy
\eqref{heat-weak}, which finishes the proof of Proposition \ref{prop2}.
\end{proof}

\begin{remark}[{\it Strong convergence of $y_{\eps\tau}$ and $y_\eps$}]
\upshape
Strengthening monotonicity of $\mfhel$, cf.\ \eqref{ass-H1}, for the 
strict monotonicity
\[
\forall\,G_1,G_2\in\R^{d\times d\times d}:
\qquad(\mfhel(G_1){-}\mfhel(G_2))\Vdots(G_1{-}G_2)
\ge c_0|G_1{-}G_2|^p,
\]
we use the argumentation after \eqref{eq:liminf.Ij} to show $\ol
y_{\eps\tau}(t) \to y_\eps(t) $ strongly in $W^{2,p}(\varOmega;\R^d)$ for
all $t\in [0,T]$. Similarly, in Proposition \ref{prop2} one can show 
$y_{\eps}(t) \to y(t) $ strongly in $W^{2,p}(\varOmega;\R^d)$. Together
with the $L^\infty$-estimate \eqref{est-u}, we can also strengthen the
weak* convergence \eqref{conv-y-} in
$L^\infty(I;W^{2,p}(\varOmega;\R^d))$ to a strong convergence 
in $L^q(I;W^{2,p}(\varOmega;\R^d))$ for all $q\in {[1,\infty[}$.
The same applies to \eqref{conv-y}.
\end{remark}



\begin{remark}[{\it Dynamical problems}]\label{rem-Galerkin}
  \upshape 
  Introducing the \emph{kinetic energy} $\frac12\varrho|\DT y|^2$ with
  a mass density $\varrho=\varrho(x)>0$ leads to an \emph{inertial
    force} $\varrho\DDT y$ in the momentum equation
  \eqref{momentum-eq}, which would make the nonlinear problem
  hyperbolic.  It is generally recognized as analytically very
  troublesome.  Here, it would work for isothermal situation like in
  Corollary~\ref{Cor:ViscoElast} if we would be able to 
  work with weak convergence, i.e.\ $\mathscr{H}$ needs to be 
  quadratic ($p=2$). Staying with $\calH$ depending on the second
  gradient $\nabla^2 y$ we would be forced to give up the
  determinant constraint $\det \nabla y>0$, which is indeed possible 
  if heat conduction is not considered. Alternatively, one may take
  $\calH$ quadratic but coercive in Hilbert space norms $H^s(\varOmega)$
  with $s>1+d/2$, such that $H^s(\varOmega)$ still embeds into
  $C^{1,\alpha}$ for some $\alpha>0$, cf.\ also
  \cite[Ch.\,9.3]{KruRou19MMCM}.
  In the anisothermal situation, it seems difficult to ensure that the
  acceleration $\DDT y\in L^2(I;H^{1+\kappa}(\varOmega;\R^{d\times d})$
  stays in duality with the velocity $\DT y$. The regularity seems
  difficult and the higher-order viscosity is inevitably very
  nonlinear to comply with frame-indifference while the corresponding
  generalization of Korn's inequality does not seem available.
\end{remark}

\begin{remark}[{\it Other transport processes: flow in porous media}]
  \upshape Beside heat transport, one can also consider other
  transport processes in a similar way. The transport coefficients can
  be pulled back as in \eqref{K-pull-back}.  For example, 
  considering mass transport for a concentration $c$ one has to make
  the free energy $\psi$ also $c$-dependent and to augmenting it by a
  capillarity-like gradient term $\frac12\varkappa|\nabla c|^2$. The 
  dissipation potential $\calR$ will then be augmented by the nonlocal
  term $\frac12 \| \mathcal{M} ( \nabla
  y,c)^{1/2}\nabla \Delta_{\mathcal{M}(\nabla y,c)}^{-1} \DT
  c \|^2_{L^2(\varOmega)}$ with 
  $\Delta_{\mathcal{M}}^{-1}: f \mapsto \mu$ denoting the linear operator
  $H^1(\varOmega)^*\to H^1(\varOmega)$ defined by the weak solution $\mu$
  to the equation $\DIV(\mathcal{M}\nabla\mu)=f$.  Considering the
  mobility tensor $\bbM=\bbM(x,c)$, we can define the pulled-back
  tensor ${\mathcal M}(x,F,c):=(\Cof F^\top)\bbM(x,c)\Cof F/\det F$
  and augment the system for the diffusion equation of the
  Cahn-Hilliard type:
\begin{subequations}\label{system-CH}
\begin{align}
&\DIV \big(\sigma_\mathrm{vi}(\nabla y,\nabla\DT y,\theta)
+
\pl_F^{}\psi(F,c,\theta)-\DIV \mfhel(\nabla^2y)\big)+g
=0 
\,,\label{momentum-eq-CH}
\\
&\DT c-\DIV \big(\mathcal{M}(\nabla y,c)\nabla\mu\big)=0
\ \ \ \text{ with }\ \ \mu=\pl_c\psi(\nabla y,c,\theta)-\varkappa\Delta c,
\label{CH}
\\&\nonumber
c_\mathrm{v}(\nabla y,c,\theta)\DT\theta
 -\DIV \big(\mathcal{K}(\nabla y,\theta)\nabla\theta\big)
 =\xi(\nabla y,\nabla\DT y,\theta)
\\&\qquad\qquad\qquad\qquad
+\theta\pl_{F\theta}^2\psi(\nabla y,c,\theta){:}\nabla\DT y
+\nabla\mu \cdot {\mathcal M}(x,\nabla y,c)\nabla\mu
\label{heat-equation-CH}
\end{align}\end{subequations}
with $\sigma_\mathrm{vi}$ as in \eqref{momentum-eq},
$c_\mathrm{v}(F,c,\theta)=-\theta\pl_{\theta\theta}^2\psi(F,c,\theta)$, and 
$\xi$ from \eqref{def-of-xi}. In \eqref{CH}, the variable $\mu$ is called 
a \emph{chemical potential}. 
One can also 
augment the model by some inelastic 
(plastic or creep-type) strain like
in \cite{RouSteTEPR} where also the inertial forces have been involved 
  and the viscosity ignored but the concept of small elastic strains imposed 
as a modeling assumption.
\end{remark}
\bigskip

\noindent{\it Acknowledgments.} A.M. is grateful for the hospitality
and support of Charles University and for partial support by Deutsche
Forschungsgemeinschaft (DFG) via the SFB\,1114 \emph{Scaling Cascades
  in Complex Systems} (subproject B01 ``Fault Networks and Scaling
Properties of Deformation Accumulation''). T.R.  is thankful for
hospitality and support of the Weierstra\ss-Institut Berlin.
Also
the partial support of the Czech Science Foundation projects
17-04301S   (as for the focuse to dissipative evolutionary systems),
18-03834S (as for the application in modeling of shape-memory alloys), and
19-29646L (as for the focuse to large strains), 
as well as through the institutional support
RVO:\,61388998 (\v CR).


\bibliographystyle{my_alpha}

\bibliography{trber7}

\newcommand{\etalchar}[1]{$^{#1}$}
\begin{thebibliography}{11}\itemsep0.1em

\bibitem[Ant98]{Antm98PUVS}
{\scshape S.~S.~Antman}.
\newblock Physically unacceptable viscous stresses.
\newblock {\em Zeitschrift angew. Math. Physik}, 49, 980--988, 1998.

\bibitem[BaC11]{BalCro11LMPI}
{\scshape J.~M.~Ball {\upshape and} E.~C.~M.~Crooks}.
\newblock Local minimizers and planar interfaces in a phase-transition model
  with interfacial energy.
\newblock {\em Calc. Var.}, 40, 501--538, 2011.

\bibitem[Bal77]{Ball77CIET}
{\scshape J.~Ball}.
\newblock Constitutive inequalities and existence theorems in nonlinear
  elastostatics.
\newblock In {\em Nonlinear Analysis and Mechanics: Heriot-Watt Symposium
  (Edinburgh, 1976), Vol. I}, pages 187--241. Res. Notes in Math., No. 17.
  Pitman, London, 1977.

\bibitem[Bal02]{Ball02SOPE}
{\scshape J.~M.~Ball}.
\newblock Some open problems in elasticity.
\newblock In P.~Newton, P.~Holmes, {\upshape and} A.~Weinstein, editors, {\em
  Geometry, Mechanics, and Dynamics}, pages 3--59. Springer, New York, 2002.

\bibitem[Bal10]{Ball10PPNE}
{\scshape J.~M.~Ball}.
\newblock Progress and puzzles in nonlinear elasticity.
\newblock In J.~Schr\"oder {\upshape and} P.~Neff, editors, {\em Poly-, Quasi-
  and Rank-One Convexity in Applied Mechanics}, volume 516 of {\em CISM Intl.
  Centre for Mech. Sci.}, pages 1--15. Springer, Wien, 2010.

\bibitem[Bar10]{Barb10NDEM}
{\scshape V.~Barbu}.
\newblock {\em Nonlinear Differential Equations of Monotone Types in Banach
  Spaces}.
\newblock Springer, New York, 2010.

\bibitem[Bat76]{Batr76TNEM}
{\scshape R.~C.~Batra}.
\newblock Thermodynamics of non-simple elastic materials.
\newblock {\em J. Elasticity}, 6, 451--456, 1976.

\bibitem[BD{\etalchar{*}}97]{BDGO97NPDE}
{\scshape L.~Boccardo, A.~Dall'aglio, T.~Gallou{\"{e}}t, {\upshape and}
  L.~Orsina}.
\newblock Nonlinear parabolic equations with measure data.
\newblock {\em J. Funct. Anal.}, 147, 237--258, 1997.

\bibitem[Bet86]{Beto86KSMC}
{\scshape D.~E.~Betounes}.
\newblock Kinematics of submanifolds and the mean curvature normal.
\newblock {\em Arch. Rational Mech. Anal.}, 96, 1--27, 1986.

\bibitem[BlG00]{BlaGui00ESNS}
{\scshape D.~Blanchard {\upshape and} O.~Guib\'e}.
\newblock Existence of a solution for a nonlinear system in
  thermoviscoelasticity.
\newblock {\em Adv.\ Diff.\ Eq.}, 5, 1221--1252, 2000.

\bibitem[BoB03]{BonBon03EUS3}
{\scshape E.~Bonetti {\upshape and} G.~Bonfanti}.
\newblock Existence and uniqueness of the solution to a 3{D} thermoelastic
  system.
\newblock {\em Electronic J. Diff. Eqs.}, 50, 1--15, 2003.

\bibitem[BoG89]{BocGal89NEPE}
{\scshape L.~Boccardo {\upshape and} T.~Gallou{\"e}t}.
\newblock Non-linear elliptic and parabolic equations involving measure data.
\newblock {\em J. Funct. Anal.}, 87, 149--169, 1989.

\bibitem[Br{\'e}73]{Brez73OMMS}
{\scshape H.~Br{\'e}zis}.
\newblock {\em Operateur Maximaux Monotones et Semi-groupes de Contractions
  dans les Espaces de Hilbert}.
\newblock North-Holland, Amsterdam, 1973.

\bibitem[Daf82]{Dafe82GSSI}
{\scshape C.~M.~Dafermos}.
\newblock Global smooth solutions to the initial boundary value problem for the
  equations of one-dimensional thermoviscoelasticity.
\newblock {\em SIAM J. Math. Anal.}, 13, 397--408, 1982.

\bibitem[{Da}L10]{DalLaz10QCGF}
{\scshape G.~{Dal Maso} {\upshape and} G.~Lazzaroni}.
\newblock Quasistatic crack growth in finite elasticity with
  non-interpenetration.
\newblock {\em Ann. Inst. H. Poinc. Anal. Non Lin.}, 27(1), 257--290, 2010.

\bibitem[Dem00]{Demo00WSCN}
{\scshape S.~Demoulini}.
\newblock Weak solutions for a class of nonlinear systems of viscoelasticity.
\newblock {\em Arch. Rat. Mech. Anal.}, 155, 299--334, 2000.

\bibitem[DSF10]{DuSoFi10TSMF}
{\scshape F.~P.~Duda, A.~C.~Souza, {\upshape and} E.~Fried}.
\newblock A theory for species migration in a finitely strained solid with
  application to polymer network swelling.
\newblock {\em J. Mech. Phys. Solids}, 58, 515--529, 2010.

\bibitem[DST01]{DeStTz01VAST}
{\scshape S.~Demoulini, D.~Stuart, {\upshape and} A.~Tzavaras}.
\newblock A variational approximation scheme for three dimensional
  elastodynamics with polyconvex energy.
\newblock {\em Arch. Rat. Mech. Anal.}, 157, 325--344, 2001.

\bibitem[FeM06]{FeiMal06NSET}
{\scshape E.~Feireisl {\upshape and} J.~M\'alek}.
\newblock On the {N}avier-{S}tokes equations with temperature-dependent
  transport coefficients.
\newblock {\em Diff. Equations Nonlin. Mech.}, pages 14pp.(electronic), Art.ID
  90616, 2006.

\bibitem[FoL07]{FonLeo07MMCV}
{\scshape I.~Fonseca {\upshape and} G.~Leoni}.
\newblock {\em Modern Methods in the Calculus of Variations: {$L^p$} spaces}.
\newblock Springer, 2007.

\bibitem[FrG06]{FriGur06TBBC}
{\scshape E.~Fried {\upshape and} M.~E.~Gurtin}.
\newblock Tractions, balances, and boundary conditions for nonsimple materials
  with application to liquid flow at small-lenght scales.
\newblock {\em Arch.\ Rational Mech.\ Anal.}, 182, 513--554, 2006.

\bibitem[FrK18]{FriKru18PNLV}
{\scshape M.~Friedrich {\upshape and} M.~Kru\v{z}\'{\i}k}.
\newblock On the passage from nonlinear to linearized viscoelasticity.
\newblock {\em SIAM J. Math. Anal.}, 50, 4426--4456, 2018.

\bibitem[GoS93]{GovSim93CSD2}
{\scshape S.~Govindjee {\upshape and} J.~C.~Simo}.
\newblock Coupled stress-diffusion: case {II}.
\newblock {\em J. Mech. Phys. Solids}, 41, 863--887, 1993.

\bibitem[HeK09]{HeaKro09IWSS}
{\scshape T.~J.~Healey {\upshape and} S.~Kr\"omer}.
\newblock Injective weak solutions in second-gradient nonlinear elasticity.
\newblock {\em ESAIM: Control, Optim. \& Cal. Var.}, 15, 863--871, 2009.

\bibitem[KPS19]{KrPeSc19?GPEM}
{\scshape M.~Kru{\v{z}}{\'\i}k, P.~Pelech, {\upshape and}
  A.~Schl{\"o}merkemper}.
\newblock Gradient polyconvexity in evolutionary models of shape-memory alloys.
\newblock {\em J. Optim. Theo. Appl.}, 2019.
\newblock Online first.

\bibitem[KrR19]{KruRou19MMCM}
{\scshape M.~Kru\v{z}{\'\i}k {\upshape and} T.~Roub{\'\i}{\v{c}}ek}.
\newblock {\em Mathematical Methods in Continuum Mechanics of Solids}.
\newblock Springer, Switzerland, 2019.

\bibitem[LeM13]{LewMuch13LERS}
{\scshape M.~Lewicka {\upshape and} P.~B.~Mucha}.
\newblock A local existence result for system of viscoelasticity with physical
  viscosity.
\newblock {\em Evolution Equations \& Control Theory}, 2, 337--353, 2013.

\bibitem[MaM09]{MaiMie09GERI}
{\scshape A.~Mainik {\upshape and} A.~Mielke}.
\newblock Global existence for rate-independent gradient plasticity at finite
  strain.
\newblock {\em J. Nonlin. Science}, 19(3), 221--248, 2009.

\bibitem[MiE68]{MinEsh68FSTL}
{\scshape R.~Mindlin {\upshape and} N.~Eshel}.
\newblock On first strain-gradient theories in linear elasticity.
\newblock {\em Intl. J. Solid Structures}, 4, 109--124, 1968.

\bibitem[Mie11]{Miel11FTDM}
{\scshape A.~Mielke}.
\newblock Formulation of thermoelastic dissipative material behavior using
  {GENERIC}.
\newblock {\em Contin. Mech. Thermodyn.}, 23(3), 233--256, 2011.

\bibitem[Mie13]{Miel13TMER}
{\scshape A.~Mielke}.
\newblock Thermomechanical modeling of energy-reaction-diffusion systems,
  including bulk-interface interactions.
\newblock {\em Discr. Contin. Dynam. Systems Ser.~S}, 6(2), 479--499, 2013.

\bibitem[MiM18]{MieMit18CEER}
{\scshape A.~Mielke {\upshape and} M.~Mittnenzweig}.
\newblock Convergence to equilibrium in energy-reaction-diffusion systems using
  vector-valued functional inequalities.
\newblock {\em J. Nonlin. Science}, 28(2), 765--806, 2018.

\bibitem[MiN18]{MieNau18?EWSK}
{\scshape A.~Mielke {\upshape and} J.~Naumann}.
\newblock On the existence of global-in-time weak solutions and scaling laws
  for {Kolmogorov}'s two-equation model of turbulence.
\newblock {\em Arch. Rational Mech. Anal.}, 2018.
\newblock Submitted. WIAS preprint 2545, arXiv:1801.02039.

\bibitem[MiR15]{MieRou15RIST}
{\scshape A.~Mielke {\upshape and} T.~Roub{\'\i}{\v{c}}ek}.
\newblock {\em Rate-Independent Systems -- Theory and Application}.
\newblock Springer, New York, 2015.

\bibitem[MiR16]{MieRou16RIEF}
{\scshape A.~Mielke {\upshape and} T.~Roub{\'\i}{\v{c}}ek}.
\newblock Rate-independent elastoplasticity at finite strains and its numerical
  approximation.
\newblock {\em Math. Models Meth. Appl. Sci.}, 6, 2203--2236, 2016.

\bibitem[MOS13]{MiOrSe14ANVM}
{\scshape A.~Mielke, C.~Ortner, {\upshape and} Y.~Seng\"ul}.
\newblock An approach to nonlinear viscoelasticity via metric gradient flows.
\newblock {\em SIAM J. Math. Anal.}, 46, 1317--1347, 2013.

\bibitem[MRS13]{MiRoSa13NADN}
{\scshape A.~Mielke, R.~Rossi, {\upshape and} G.~Savar\'{e}}.
\newblock Nonsmooth analysis of doubly nonlinear evolution equations.
\newblock {\em Calc. Var PDE}, 46(1-2), 253--310, 2013.

\bibitem[MRS18]{MiRoSa18GERV}
{\scshape A.~Mielke, R.~Rossi, {\upshape and} G.~Savar\'{e}}.
\newblock Global existence results for viscoplasticity at finite strain.
\newblock {\em Arch. Rational Mech. Anal.}, 227(1), 423--475, 2018.

\bibitem[Nef02]{Neff02KFIN}
{\scshape P.~Neff}.
\newblock On {K}orn's first inequality with non-constant coefficients.
\newblock {\em Proc. Royal Soc. Edinburgh}, 132A, 221--243, 2002.

\bibitem[Pod02]{Podi02CISM}
{\scshape P.~Podio-Guidugli}.
\newblock Contact interactions, stress, and material symmetry, for nonsimple
  elastic materials.
\newblock {\em Theor. Appl. Mech.}, 28-29, 261--276, 2002.

\bibitem[Pom03]{Pomp03KFIV}
{\scshape W.~Pompe}.
\newblock Korn's {F}irst {I}nequality with variable coefficients and its
  generalization.
\newblock {\em Comment. Math. Univ. Carolinae}, 44, 57--70, 2003.

\bibitem[RoS18]{RouSteTEPR}
{\scshape T.~Roub{\'\i}{\v{c}}ek {\upshape and} U.~Stefanelli}.
\newblock Thermodynamics of elastoplastic porous rocks at large strains,
  earthquakes, and seismic waves.
\newblock {\em SIAM J. Appl. Math.}, 78, 2597--2625, 2018.

\bibitem[Rou04]{Roub04MMES}
{\scshape T.~Roub{\'{\i}}{\v{c}}ek}.
\newblock Models of microstructure evolution in shape memory alloys.
\newblock In P.~P.~Castaneda, J.~Telega, {\upshape and} B.~Gambin, editors,
  {\em Nonlin. Homogenization and its Appl. to Composites, Polycrystals and
  Smart Mater.}, pages 269--304. Kluwer, Dordrecht, 2004.

\bibitem[Rou09]{Roub09TVES}
{\scshape T.~Roub{\'{\i}}{\v{c}}ek}.
\newblock Thermo-visco-elasticity at small strains with ${L}^1$-data.
\newblock {\em Quarterly Appl. Math.}, 67, 47--71, 2009.

\bibitem[Rou13]{Roub13NPDE}
{\scshape T.~Roub{\'\i}{\v{c}}ek}.
\newblock {\em Nonlinear Partial Differential Equations with Applications}.
\newblock Birkh\"auser, Basel, 2nd edition, 2013.

\bibitem[{\v{S}}il85]{Silh88PTNB}
{\scshape M.~{\v{S}}ilhav\'{y}}.
\newblock Phase transitions in non-simple bodies.
\newblock {\em Arch. Rat. Mech. Anal.}, 88, 135--161, 1985.

\bibitem[Ste15]{Stei15GFCM}
{\scshape P.~Steinmann}.
\newblock {\em Geometric Foundations of Continuum Mechanics}.
\newblock Springer, 2015.
\newblock LAMM Vol.\,2.

\bibitem[Tou62]{Toup62EMCS}
{\scshape R.~Toupin}.
\newblock Elastic materials with couple stresses.
\newblock {\em Arch. Rat. Mech. Anal.}, 11, 385--414, 1962.

\bibitem[TrA86]{TriAif86GALD}
{\scshape N.~Triantafyllidis {\upshape and} E.~Aifantis}.
\newblock A gradient approach to localization of deformation. {I}.
  {H}yperelastic materials.
\newblock {\em J. Elast.}, 16, 225--237, 1986.

\bibitem[Tve08]{Tved08QEVS}
{\scshape B.~Tvedt}.
\newblock Quasilinear equations of viscoelasticity of strain-rate type.
\newblock {\em Arch. Rational Mech. Anal.}, 189, 237--281, 2008.

\bibitem[Vis96]{Visi96MPT}
{\scshape A.~Visintin}.
\newblock {\em Models of Phase Transitions}.
\newblock Birkh{\"a}user, Boston, 1996.

\end{thebibliography}

\end{sloppypar}
\end{document}